\input{style/arxiv-general.cfg}
\documentclass[MSNbibl,number,citesort,seceqn,rotating,dvips]{arxbj}
\makeatletter
   \@ifpackageloaded{graphicx}{}{\usepackage{graphicx}}
\makeatother
\usepackage{upgreek}

\aid{0}
\volume{21}
\issue{2}
\pubyear{2015}
\firstpage{781}
\lastpage{831}
\doi{10.3150/13-BEJ587} % kopijuoti is 'New paper accepted'

\makeatletter

\newcommand{\rrvert}{\vert}
\newcommand{\llvert}{\vert}
\newcommand{\eqref}[1]{(\ref{#1})}

\newcommand{\C}{\mathbb C}
\newcommand{\R}{\mathbb R}
\newcommand{\rr}{\mathbb R}
\newcommand{\Z}{\mathbb{Z}}
\newcommand{\n}{^{(n)}}

\newcommand{\varthetab}{{\bolds\vartheta}}

\newcommand{\pr}{^{\prime}}

\newcommand{\uut}[1] {
\begin{array}[t]{c}{#1}\vspace*{-2pt}\\[-.06mm]
\widetilde{}\vspace{-3.1mm}% espace a regler ici
\end{array}
}

\newcommand{\IC}{\mathbb{C}}
\newcommand{\IG}{\mathbb{G}}
\newcommand{\IN}{\mathbb{N}}
\newcommand{\IP}{\mathrm{P}}
\newcommand{\IR}{\mathbb{R}}
\newcommand{\IS}{\mathbb{S}}
\newcommand{\IZ}{\mathbb{Z}}

\newcommand{\E}{\mathrm{E}}
\newcommand{\cF}{\mathcal{F}}

\newcommand{\Cov}{\operatorname{Cov}}
\newcommand{\Var}{\operatorname{Var}}
\newcommand{\cum}{\operatorname{cum}}

\newcommand{\bdelta}{{\bolds\delta}}

\newcommand{\bzeta}{{\bolds\zeta}}
\newcommand{\tr}{\prime}

\newcommand{\gammar}{\gamma^{U}}

\newcommand{\bQy}{\mathbf{Q}}
\newcommand{\bQu}{\mathbf{Q}^{U}}
\newcommand{\bzetay}{\bzeta}
\newcommand{\bzetau}{\bzeta^{U}}

\newcommand{\ef}[1]{\frak{f}_{#1}}
\newcommand{\efh}[1]{\hat{\frak{f}}_{#1}}

\newcommand{\efs}[1]{\!\stackrel{\scriptscriptstyle\circ}{\frak{f}}_{#1}\!}

\newcommand{\efth}[1]{\uut{\hat{\frak{f}}}{}_{#1}}

\newcommand{\Lr}[1]{{L^{U}_{#1}}}
\newcommand{\Lnr}[1]{L_{#1}}
\newcommand{\Lrh}[1]{\uut{\hat{L}}{}_{#1}}
\newcommand{\Lnryh}[1]{\hat{L}_{#1}}
\newcommand{\Lnc}[1]{L^{U}_{#1}}
\newcommand{\Lnch}[1]{\hat{L}^{U}_{#1}}

\newcommand{\Zx}[1]{Z^{X}_{#1}}
\newcommand{\Zy}[1]{Z_{#1}}
\newcommand{\Zu}[1]{Z^{U}_{#1}}
\newcommand{\Zyh}[1]{\hat{Z}_{#1}}
\newcommand{\Zuh}[1]{\hat{Z}^{U}_{#1}}

\newcommand{\Zyrh}[1]{\uut{\hat{Z}}{}_{#1}}
\newcommand{\dex}[1]{\bdelta^{X}_{#1}}
\newcommand{\dey}[1]{\bdelta_{#1}}
\newcommand{\deu}[1]{\bdelta^{U}_{#1}}
\newcommand{\deyh}[1]{\hat{\bdelta}_{#1}}
\newcommand{\deyrh}[1]{\uut{\hat{\bdelta}}{}_{#1}}

\newcommand{\aeyh}[1]{\hat{a}_{#1}}

\newcommand{\bey}[1]{\mathbf{b}_{#1}}

\newcommand{\beyhpr}[1]{\hat{\mathbf{b}}^{\prime}_{#1}}

\newcommand{\beyh}[1]{\hat{\mathbf{b}}_{#1}}
\newcommand{\beuh}[1]{{\mathbf{b}}^{U}_{#1}}
\newcommand{\beyrh}[1]{\uut{\hat{\mathbf{b}}}{}_{#1}}

\newtheorem{Prop}{Proposition}[section]
\newtheorem{theorem}{Theorem}[section]
\newtheorem{Corol}{Corollary}[section]
\newproclaim{bem}{Remark}[section]
\newproclaim{Assumption}{Assumption}%[section]
\newtheorem{lemma}{Lemma}[section]

\makeatother

\begin{document}
\begin{frontmatter}

\title{Of copulas, quantiles, ranks and spectra: An~$L_1$-approach to
spectral analysis}
\runtitle{$L_1$-approach to
spectral analysis}

\begin{aug}
%%%% inicialai - be tarpu
\author[A]{\inits{H.}\fnms{Holger} \snm{Dette}\thanksref{A}\ead[label=e1]{holger.dette@rub.de}},
\author[B,C]{\inits{M.}\fnms{Marc} \snm{Hallin}\corref{}\thanksref{B,C}\ead[label=e2]{mhallin@ulb.ac.be}},
\author[A]{\inits{T.}\fnms{Tobias} \snm{Kley}\thanksref{A}\ead[label=e3]{tobias.kley@rub.de}}
\and\\
\author[A]{\inits{S.}\fnms{Stanislav}~\snm{Volgushev}\thanksref{A}\ead[label=e4]{stanislav.volgushev@rub.de}}
%%\runauthor{} %% auto
\address[A]{Fakult\"{a}t f\"{u}r Mathematik, Ruhr-Universit\"{a}t Bochum, 44780 Bochum, Germany.\\
\printead{e1,e3,e4}}
\address[B]{ECARES, Universit\'{e} Libre de Bruxelles, 1050 Brussels, Belgium.
\printead{e2}}
\address[C]{ORFE, Princeton University, Princeton, NJ 08544, USA}
\end{aug}

% HISTORY:
\received{\smonth{5} \syear{2012}}
\revised{\smonth{4} \syear{2013}}

% ABSTRACT
%
\begin{abstract}
In this paper, we present an alternative method for the spectral
analysis of
{a univariate, strictly stationary} time series
$\{Y_t\}_{t\in\Z}$.
We define a ``new'' spectrum as the Fourier transform of the
differences between copulas
of the pairs $(Y_t,Y_{t-k})$ and the independence copula.
This object is called a {\it copula spectral density kernel}
and allows to separate the marginal and serial aspects of a time series.
We show that this spectrum is closely related to the concept of
quantile regression.
Like quantile regression, which provides much more information about
conditional distributions than classical location-scale regression
models, copula spectral density kernels are
more informative than traditional spectral densities obtained from
classical autocovariances.
In particular, copula spectral density kernels, in their population
versions, provide (asymptotically provide, in their sample versions)
a complete description of the {copulas}
of all pairs $(Y_t,  Y_{t-k})$. Moreover, they inherit the robustness
properties of classical quantile regression,
and do not require any distributional assumptions
such as the existence of finite moments.
In order to estimate the copula spectral density kernel,
we introduce rank-based Laplace periodograms
which are calculated as bilinear forms of weighted $L_1$-projections of
the ranks
of the observed time series onto a harmonic regression model. We
establish the asymptotic
distribution of those periodograms, and the consistency of adequately
smoothed versions.
The finite-sample properties of the new methodology, and its potential for
applications are briefly investigated by simulations and a short
empirical example.
\end{abstract}

% KEYWORDS
% visi is mazosios raides ir pagal abecele
\begin{keyword}
\kwd{copulas}
\kwd{periodogram}
\kwd{quantile regression}
\kwd{ranks}
\kwd{spectral analysis}
\kwd{time reversibility}
\kwd{time series}
\end{keyword}

\end{frontmatter}

%s1 #&#
\section{Introduction}\label{introintro}

%s1.1 #&#
\subsection{The location-scale paradigm}\label{intro1}

Whether linear or not, most traditional time series models are of the
conditional location/scale type: {conditionally} on past values
$Y_{t-1},Y_{t-2}, \ldots$\,, the random variable $Y_t$ is of the form
%
%e1.1 #&#
\begin{equation}
\label{locscale} Y_t = \psi %_{\bolds\vartheta}
(Y_{t-1},Y_{t-2},
\ldots ) + \sigma % _{\bolds\vartheta}
(Y_{t-1},Y_{t-2}, \ldots )
\varepsilon_t,\quad\quad t\in\Z,
\end{equation}
where $\{\varepsilon_t \}_{t \in\IZ}$ is white noise (either strong
or weak, depending on the authors -- here, by white noise we
throughout mean strong, i.e., independent white noise), and
$\varepsilon_t$ is independent of (in the case of weak white noise,
orthogonal to) $Y_{t-1},Y_{t-2}, \ldots $\,. The $(Y_{t-1},Y_{t-2},
\ldots)$-measurable functions $\psi
%_\varthetab
$ and $\sigma
%_\varthetab
$ are (conditional) location and scale functions,
possibly para\-metrized by some $\varthetab$. Equation (\ref
{locscale}) may characterize a data-generating process -- in which case
``$=$'' in (\ref{locscale}) is to be considered as ``almost sure
equality'' -- or, more generally, it simply describes $Y_t$'s
conditional (on $Y_{t-1},Y_{t-2}, \ldots$) distribution -- and
``$=$'' is to be interpreted as ``equality in (conditional)
distribution''. Such distinction is, however, irrelevant from a
statistical point of view, as it has no impact on likelihoods.

In model \eqref{locscale}, the distribution of $Y_t $ conditional on
$Y_{t-1},Y_{t-2}, \ldots$ is nothing but the distribution of
$\varepsilon_t$, rescaled by the conditional scale parameter $ \sigma
(Y_{t-1},Y_{t-2}, \ldots)$ and shifted by the conditional location
parameter $\psi(Y_{t-1},Y_{t-2}, \ldots)$. Sophisticated as they
may be, the mappings
\[
(Y_{t-1},Y_{t-2}, \ldots) \mapsto\bigl(\psi(Y_{t-1},Y_{t-2},
\ldots), \sigma(Y_{t-1},Y_{t-2}, \ldots )\bigr)
\]
only can account for a very limited type of dynamics for the process $\{
Y_t\}_{t \in\IZ} $. The volatility dynamics for such models, for
instance, are quite poor, being of a pure rescaling nature. In
particular, no impact of past values on skewness, kurtosis, tails, can
be reflected.
All standardized conditional distributions strictly coincide with that
of $\varepsilon$, and all conditional $\tau$-quantiles, hence all
values at risk, follow, {irrespectively} of $\tau$, from those of
$\varepsilon$ via one single linear transformation.

Note that the interpretation of $\psi$ and $\sigma$ depends on the
identification constraints on $\varepsilon$: if $\varepsilon$ is
assumed to have mean zero and variance one, then $\psi$ and $\sigma$
are a conditional mean and a conditional standard error, respectively.
In this case, models of the form (\ref{locscale}) clearly belong to
the $L_2$-Gaussian legacy. If $\varepsilon$ is assumed to have median
zero and
expected absolute deviation or median absolute deviation one, $\psi$
and $\sigma$ are a conditional median and a conditional expected or
median absolute deviation.

On the basis of these ``remarks'', the following questions naturally arise:
Can we do better? Can we go beyond that (conditional) ``location-scale
paradigm''? Can we model richer dynamics under which the conditional
quantiles of $Y$ are not just a shifted and rescaled version of those
of $\varepsilon$, and under which the whole conditional distribution
of $Y_t$%conditional on $(Y_{t-1},Y_{t-2}, \ldots )$
, not just its location and scale, can be affected by the past? And,
can we achieve this in a statistically tractable way?

%structural assumptions, however, such as $\{Y_t\}$ stationary, or $
%

%s1.2 #&#
\subsection{Marginal and serial features}\label{intro2}

Another drawback of models of the form (\ref{locscale}) is their
sensitivity to nonlinear marginal transformations. If two statisticians
observe the same time series, but measure it on different scales, $Y_t$
and $Y_t^3$ or $\mathrm{e}^{Y_t}$, for instance, and both adjust a model of the
form (\ref{locscale}) to their measurements, they will end up with
drastically different analyses and predictions. The only way to avoid
this problem
consists in disentangling the marginal (viz., related to the scale of
measurement) aspects of the series under study from its serial aspects,
that is, basing the description of serial dependence features on
quantities such as the $F(Y_t)$'s, where $F$ is $Y_t$'s marginal
distribution function. Those quantities do not depend on the
measurement scale since they are invariant under {continuous strictly
increasing} transformations.

This point of view is closely related to the concept of copulas (see
Nelsen \cite{Nelsen2006} or Genest and Favre \cite{GenestFavre2007}).
Consider, for instance, a strictly stationary Markovian process $\{Y_t\}
_{t \in\IZ}$
of order one. This process is fully characterized by the joint
distribution of $(Y_t, Y_{t-1})$ or, equivalently, by the marginal
distribution function $F$ (equivalently, the quantile
function $F^{-1}$) of $Y_t$, along with the joint distribution of
$(U_t,U_{t-1}) :=(F(Y_t),  F(Y_{t-1}))$, a ``serial copula of order
one''. In such a description, the marginal features of the process $\{
Y_t\}_{t \in\IZ}$ are entirely described by $F$, independently of the
serial features, that are accounted for by the serial copula. Two
statisticians observing the same phenomenon but recording $Y_t$ and
$\mathrm{e}^{Y_t}$, respectively, would use distinct quantile functions, but
they would agree on serial features.

In more general cases, serial copulas of order one are not sufficient,
and higher-order or multiple copulas may be needed. Note that the
description of the model in this context is clearly ``in distribution'':
$U_t$ is not related to $U_{t-1}$ through any direct interpretable
``almost sure relation'' reflecting some ``physical'' data-generating mechanism.

%s1.3 #&#
\subsection{A new nonparametric approach}\label{intronew}

The objective of this paper is to show how to overcome the limitations
of conditional location-scale modelling described in Sections~\ref{intro1} and \ref{intro2}, and to provide statistical tools for a
fully general approach to time series modelling. Not surprisingly,
those tools are essentially related to copulas, quantiles and ranks.
The traditional nonparametric techniques, such as spectral analysis (in
its usual $L_2$-form), which only account for second-order serial
features, cannot handle such objects, and we therefore propose and
develop an original, flexible and fully nonparametric $L_1$-spectral
analysis method.

While classical spectral densities are obtained as Fourier transforms
of classical covariance functions, we rather
define spectral density \textit{kernels}, associated with covariance
\textit{kernels} of the form (for $(\tau_1,\tau_2) \in(0,1)^2$)
%
%e1.2 #&#
\begin{equation}
\label{lcov} \gamma_k(x_1 , x_2):=
\operatorname{Cov}\bigl(I\{ Y_t\leq x_1\}, I\{
Y_{t-k}\leq x_2\}\bigr)
\end{equation}
(Laplace {cross-}covariance kernels) or
%
%e1.3 #&#
\begin{equation}
\label{ccov} \gammar_k(\tau_1,\tau_2):=
\operatorname{Cov}\bigl(I\{ U_t\leq\tau_1\}, I\{
U_{t-k}\leq\tau_2\}\bigr)
\end{equation}
(copula {cross-}covariance kernels), where $U_t := F(Y_t)$ and $F$
denotes the marginal distribution of
the strictly stationary process $\{Y_t\}_{t \in\IZ}$ and $I \{A\}$
stands for the indicator function of $A$.
Contrary to covariance functions, the
{\it kernels} $\{\gamma_k(x_1 , x_2) | x_1,x_2 \in\R\}$ and $ \{
\gammar_k(\tau_1,\tau_2)| \tau_1,\tau_2 \in(0,1)\}$
allow for a complete description of arbitrary bivariate distributions
for the couples $(Y_t, Y_{t-k})$ and
the corresponding copulas, respectively, and thus escape the
conditional location-scale
paradigm discussed in Section~\ref{intro1}. They are able to account
for sophisticated dependence features that covariance-based
methods are unable to detect, such as time-irreversibility, tail
dependence, varying conditional skewness or kurtosis, etc. And, in view
of the desired separation between marginal and serial features
expressed in Section~\ref{intro2}, special virtues, such as
invariance/equivariance (with respect to continuous order-preserving
marginal transformations), can be expected from the copula covariance
kernels defined in \eqref{ccov}.

Classical nonparametric spectral-based inference methods have proven
quite effective (see, e.g., Granger \cite{Granger1964}, Bloomfield \cite{Bloomfield1976}),
essentially in a Gaussian context, where dependencies are fully
characterized by autocovariance functions. Therefore, it can be
anticipated that similar methods, based on estimated versions of
Laplace or copula spectral kernels (associated with Laplace and copula
covariance kernels, respectively) would be quite useful in the study of
series exhibiting those features that classical covariance-related
spectra cannot account for.

Estimation of Laplace and copula spectral kernels, however, calls for a
substitute to the \textit{ordinary periodogram} concept considered in
the classical approach. We therefore introduce Laplace and copula
\textit{periodogram kernels}. While ordinary periodograms are defined via
least squares regression of the observations on the sines and cosines
of the harmonic basis, our periodogram kernels are obtained via
quantile regression in the Koenker and
  Bassett \cite{KoenkerBassett1978} sense. A study
of their asymptotic properties shows that,
just as ordinary periodograms, they produce asymptotically unbiased
estimates (more precisely, the mean of their asymptotic distribution is
$2\uppi$ times the corresponding spectrum), and we therefore also
consider smoothed versions that yield consistency. Asymptotic results
show that copula periodograms, as anticipated, are preferable to the
Laplace ones, as their asymptotic behavior only depends on the
bivariate copulas of the pairs $(U_t, U_{t-k})$, not on the (in general
unknown) marginal distribution $F$ of the $Y_t$'s.

Unfortunately, copula periodogram kernels are not statistics, since
their definition involves the transformation of $Y_t$ into $U_t$, hence
the knowledge
of the marginal distribution function $F$. We therefore introduce a
third periodogram kernel, based on the empirical version $\hat{F}_n$
of $F$, that is, on the
\textit{ranks} of the random variables $Y_1,\ldots, Y_n$, and
establish, under mild assumptions, the asymptotic equivalence of that
rank-based Laplace periodogram with the copula one. Smoothed rank-based
Laplace periodogram kernels, accordingly, seem to be the adequate tools
in this context. We conclude with a brief numerical illustration --
simulations and an empirical application -- of their potential use in
practical problems. % more can be be found in a companion paper (Dette
%et al. 2012) \textbf{Soll das hier so stehen bleiben??}.

%s1.4 #&#
\subsection{Review of related literature}

Quantities of the form (\ref{lcov}) and (\ref{ccov})
naturally come into the picture when the {\it clipped processes} $(I\{
Y_t \leq x\})_{t \in\IZ}$ and $(I\{U_t \leq\tau\})_{t \in\IZ}$
are investigated.
Such clipped processes have been
{considered} earlier in the literature (see, for instance, Kedem \cite{Kedem1980}).
In the field of signal processing, the idea to replace the quadratic
loss by other loss functions
has been discussed by Katkovnik \cite{Katkovnik1998}, who proposes using
$L_p$-distances and analyzes the properties of the resulting
\textit{M-periodograms}.
Hong \cite{Hong2000} used the Laplace covariances corresponding to
positive lags
to construct a test for serial dependence.
Linton and Whang \cite{LintonWhang2007} considered sequences of Laplace
autocorrelations $\gamma_k(\tau,\tau) / \gamma_0(\tau,\tau)$
(called {\it quantilogram} by these authors) in order to
test for directional predictability.
{Mikosch and
  Zhao \cite{MikoschZhao2013} define a periodogram generated from a
suitable sequence of indicator functions of rare events.}

In a pioneering paper, Li \cite{Li2008} suggested least absolute
deviation estimators in a harmonic regression model assuming
that the median of the random variables $Y_t$ is zero. The focus of
this author is on the quantities of the form (for $\omega\in(0,\uppi
)$; throughout, i stands for the root of $-1$)
\[
f_{0,0}(\omega) = \frac{1}{2\uppi} \sum_{k \in\IZ}
\gamma_k(0,0) \exp(\mathrm{i} k \omega),\quad\quad \omega\in(0,\uppi),
\]
the collection of which he calls the \textit{Laplace spectrum}. He
constructs an asymptotically unbiased estimator for a quantity which differs
from $f_{0,0}(\omega_j)$ ($\omega_j$ the $j$th Fourier frequency) by
a factor involving $1/(F'(0))^2$ and, in Li \cite{Li2011},
extends his results to arbitrary quantiles. An important drawback of
Li's method is that it requires estimates of
the quantity $F'(0)$ in order to obtain an estimate of the Laplace
spectrum; moreover, the consistency of
a smoothed version of his estimates is not established.
{More recently}, Hagemann \cite{Hagemann2011} proposed an alternative method to
estimate the Laplace spectrum (called {\it quantile spectrum} by this author),
which is based on the estimation of a linearization of Li \cite
{Li2008}'s statistic. This approach does
not suffer from the drawbacks of Li's method, and yields consistent
estimates avoiding estimation of the {marginal} density; on the other
hand, it does not allow a direct interpretation in terms of (weighted)
absolute deviation estimates.

In order to obtain a complete description of the two-dimensional
distributions at lag $k$, Hong \cite{Hong1999} introduced a {\it
generalized spectrum} defined
as the covariance $\operatorname{Cov} (\mathrm{e}^{{\rm i}x_1Y_t}, \mathrm{e}^{{\rm
i}x_2Y_{t+k}})$; this concept
was used by Chung and Hong \cite{ChungHong2007}
to test for directional predictability.
Recently, Lee and Rao \cite{LeeRao2011} considered a
Fourier transform of the differences between the joint density of the
pairs $(Y_t, Y_{t-k})$ and the product of their marginal densities to
investigate
serial dependence.
Unlike ours, these methods are not invariant with respect to
transformations of the marginal
distributions.

Finally, there exist some recent proposals using pair-copula
constructions to describe dependencies in the time-domain.
Domma, Giordano and
  Perri \cite{DommaEtAl2009} assume first-order Markov dependence, so that only
distributions of pairs $(Y_t, Y_{t+1})$ at lag $k=1$ need to be considered.
Smith \textit{et~al.} \cite{SmithEtAl2010} decompose the distribution at a point in time,
conditional upon the past, into the product
of a sequence of bivariate copula densities and the marginal density,
known as D-vine (Bedford and Cooke \cite{BedfordCooke2002}).

The approach presented in this paper differs from these references in
many important aspects.
Essentially, it combines their attractive features while avoiding some
of their drawbacks.
It shares the quantile-based flavor of
Kedem \cite{Kedem1980}, Linton and Whang \cite{LintonWhang2007}, Li \cite{Li2008,Li2011}
and Hagemann \cite{Hagemann2011}.
In contrast to the latter, however, we do not focus on a particular
quantile, and consider copula \textit{cross-}covariances $\gamma_k^U
(\tau_1,\tau_2)$ for
{\it all} pairs $(\tau_1, \tau_2)$,
while Li \cite{Li2008,Li2011} and Hagemann \cite{Hagemann2011}
restrict to the case $\tau_1 = \tau_2$. As a consequence,
we obtain, as in the characteristic function approach of Hong \cite{Hong1999},
a complete characterization of the dependencies among the pairs $(Y_t,
Y_{t-k})$.
This allows to address such important features as time reversibility
(see Proposition~\ref{timerev}) or tail dependence in general.
By replacing the original observations with their ranks, we furthermore
achieve an
attractive invariance property with respect to modifications of
marginal distributions, which is not satisfied
in the case of Hong \cite{Hong1999}'s method.
Moreover, we also avoid the scaling problem of Li's estimates and
establish the consistency of a smoothed version of periodograms.
Finally, because our method is
related to the concept of copulas, it allows to separate the marginal
and serial aspects of a time series,
which should make it attractive for practitioners.

%s1.5 #&#
\subsection{Outline of the paper}\label{introutline}
The paper is organized as follows. In Section~\ref{secLaplace}, we
introduce the concepts of {\it Laplace} and {\it copula
cross-covariance kernels} which, in this quantile-based approach, are
to replace the ordinary autocovariance function. The corresponding
spectra and periodograms are introduced in Sections~\ref{secLaplaceSpectra} and \ref{secperiod}, respectively. Section~\ref{321315458} deals with the asymptotic properties of the Laplace, copula, and
rank-based
Laplace periodograms. In Section~\ref{sec4}, smoothed periodograms are
considered, and the smoothed rank-based
Laplace periodogram kernel is shown to be
a consistent estimator of the copula spectral density. Some numerical
illustration is provided in Section~\ref{sec5}, and most of the
technical details are concentrated in an appendix.\vspace*{-2pt}

%s2 #&#
\section{An $L_1$-approach to spectral analysis}

%s2.1 #&#
\subsection{The Laplace and copula cross-covariance kernels}\label{secLaplace}
Covariances clearly are not sufficient for describing a serial copula.
We therefore introduce the following concept, which will be convenient
for that purpose.
Let $\{Y_t\}_{t\in\Z}$ be a strictly stationary process
%, with absolutely continuous distribution and nonvanishing densities,
%hence monotone marginal distribution function $F_Y$,
and define its {\it copula cross-covariance kernel} of lag $ k\in\IZ$
of $\{Y_t\}_{t \in\IZ}$ as\vspace*{-2pt}
\[
\gammar_k := \bigl\{ \gammar_k(\tau_1,
\tau_2) | (\tau_1, {\tau _2}) \in(0,
1)^2 \bigr\},
\]
where $\gammar_k(\tau_1,\tau_2)$ is defined in \eqref{ccov}.
Similarly, define the {\it Laplace cross-covariance kernel} of
lag $k\in\IZ$ of $\{Y_t\}_{t \in\IZ}$ as\vspace*{-2pt}
\[
\gamma_k := \bigl\{ \gamma_k(x_1 ,
x_2) | (x_1 , x_2) \in\R^2 \bigr
\},
\]
where $ \gamma_k(x_1 , x_2)$ is defined in \eqref{lcov}.
Contrary to traditional cross-covariances, copula and Laplace
cross-covariance kernels exist for all $k$ (no finite variance
assumption needed). The words ``covariance'' and ``cross-covariance''
are used out of time series classical terminology; but we only consider
covariances of indicators, which always exist, and provide a canonical
description of their joint distributions. The copula cross-covariance
kernel of order $k$ indeed entirely characterizes the joint
distribution of $(U_t, U_{t-k})$, and conversely, without requiring any
information on the distribution function $F$ of $Y_t$. Along with $F$,
the copula cross-covariance kernel of order $k$ entirely characterizes
the Laplace cross-covariance kernel of order $k$ and the joint
distribution of $(Y_t, Y_{t-k})$, and conversely. If $\int x^2\,{\rm d}F
<\infty$, the distribution function $F$ of $Y_t$ and the collection of
copula cross-covariance kernels of all orders jointly characterize the
autocovariance function of $\{Y_t\}_{t\in\Z}$.\vspace*{-2pt}

%Cross-covariance kernels only involve bivariate joint distributions
%and bivariate copulas, but naturally extend to higher-dimensional
%distributions as follows: the {\it copula} and {\it Laplace cumulant
%kernels } of order $m$ $\{Y_t\}$ are the functions, from $(0,1)^{m+1}$
%and $\R^{m+1}$ to $[ -1/2^{m+1}, 1/2^{m+1}]$, respectively,
%(\tau_0,\tau_1,\ldots,\tau_m) \mapsto\gammar_{k_1\ldots k_m}(\tau_0,
%:= {\rm E}\left[(I\{ U_t\leq\tau_0\} - \tau_0) (I\{ U_{t-k_1}\leq\tau_1
% and
%$$(x_0,\ldots, x_m)\mapsto\gamma_{k_1\ldots k_m} (x_0,\ldots,
%x_m):= \gammar_{k_1\ldots k_m}(F_Y(x_0),\ldots, F_Y(x_m)).
%$$
%
%We then have the following result (see the Appendix for a proof).
%kernels} of order less than or equal to $q$ jointly determine all
%joint distributions and all {\em serial copulas}, respectively, up to
%dimension $q$, and vice-versa. }
%In the sequel, however, we concentrate on cross-covariance kernels.
%
% \textbf{Should we make a reference to the intended companion paper
%here?}

%s2.2 #&#
\subsection{The Laplace and copula spectral density kernels}\label{secLaplaceSpectra}

Assume that the Laplace cross-covariance kernels $\gamma_k$
(equivalently, the copula cross-covariance kernels $\gammar_k$), $k
\in\mathbb{Z}$ are absolutely summable, that is, assume that they
satisfy\vspace*{-2pt}
\[
\sum_{k=-\infty}^\infty\bigl |\gamma_k (x_1 ,  x_2 )
\bigr |
< \infty \qquad\mbox{for all }(x_1,x_2 ) \in\rr^2
.
\]
Then, $\gamma_k$ admits the representation\vspace*{-2pt}
\[
\gamma_k (x_1 , x_2 ) = \int
_{-\uppi}^\uppi \mathrm{e}^{{\rm i}k\omega}\ef{x_1 ,
x_2 }(\omega)\,{\rm d}\omega,\quad\quad (x_1, x_2 ) \in
\rr^2
\]
with\vspace*{-2pt}
%
%e2.1 #&#
\begin{equation}
\label{sd1} \ef{x_1 , x_2 }(\omega) :=
\frac{1}{2\uppi}\sum_{k=-\infty}^\infty
\gamma_k(x_1 , x_2 )\mathrm{e}^{-{\rm i}k\omega} ,\quad\quad
(x_1,x_2 ) \in\rr^2.
\end{equation}
The collection $\{\omega\mapsto\ef{x_1 , x_2 }(\omega)|
(x_1,x_2 ) \in\rr^2\}$, call it the {\it Laplace spectral density
kernel}, is such that
each mapping $\omega\in(-\uppi,\uppi] \mapsto\ef{x_1 , x_2 }(\omega
)$, $(x_1 , x_2 ) \in\rr^2$, is continuous and satisfies (writing
$\bar z$ for the complex conjugate of $z\in\C$)
%
%e2.2 #&#
\begin{equation}
\label{transpo} \ef{x_1, x_2} (-\omega) =
\ef{x_2, x_1} (\omega) = {\overline{\ef
{x_1, x_2} (\omega)}}.
\end{equation}
Similarly define the {\it copula spectral density kernel} as
%
%e2.3 #&#
\begin{equation}
\label{sd2} \ef{q_{\tau_1} , q_{\tau_2}}(\omega) =
\frac{1}{2\uppi}\sum_{k=-\infty}^\infty
\gammar_k(\tau_1,\tau_2)\mathrm{e}^{-{\rm i}k\omega},\quad\quad (
\tau_1, \tau_2 ) \in(0,1)^2,
\end{equation}
where $q_{\tau_i} := F^{-1}(\tau_i)$ ($i=1,2$).
Note that $\ef{q_{\tau_1} , q_{\tau_2}}$ is the Fourier transform of the
differences between copulas
of the pairs $(Y_t,Y_{t-k})$ and the independence copula.
Clearly, the same identity (\ref{transpo}) holds for $\ef{q_{\tau_1}
, q_{\tau_2}}(\omega)$ as for $\ef{x_1, x_2} (\omega)$.

Throughout this paper, we denote by $\stackrel{d}{=}$ equality in
distribution and
define $\Im z$ and $\Re z$ as the imaginary and real part of $z \in\IC
$, respectively. Obviously, we have
$\Im\ef{x_1, x_2} (\omega) = 0$ for all $\omega$ if and only if $
\gamma_k(x_1 , x_2 ) = \gamma_{-k}(x_1 , x_2 )$ for all $k$,
and we obtain the following result.

%pr2.1 #&#
\begin{Prop} \label{timerev}
The following statements are
equivalent: %
\begin{enumerate}[(3)]
\item[(1)] $(Y_t, Y_{t+k})\stackrel{d}{=}(Y_t, Y_{t-k})$ for all
$k\in\mathbb{Z}$ (pairwise time-reversibility);
\item[(2)]
$\Im\ef{x_1, x_2} (\omega) = 0$ for all $\omega\in(0,\uppi)$ and
$(x_1, x_2)\in\mathbb{R}^2$;
\item[(3)]
$\Im\ef{q_{\tau_1}, q_{\tau_2}} (\omega) = 0$ for all $\omega\in
(0,\uppi)$ and $(\tau_1, \tau_2)\in(0,1)^2$.
\end{enumerate}
\end{Prop}

%s2.3 #&#
\subsection{The Laplace, copula  and rank-based Laplace periodogram
kernels}\label{secperiod}

The copula cross-covariance kernels describe the serial behavior of
$Y_t$'s quantiles.
%Earlier attempts have been made in that direction, the most typical of
%which is perhaps Koenker and Xiao (2004, 2006)'s definition of {\it
%quantile AR($p$) models}, of the form (for $p=1$)
% $
%
%%\begin{equation}\label{KX}
%Y_t \stackrel{d}{=} \theta_0(U_t)+ \theta_{1}(U_t)Y_{t-1}
%$
%%\end{equation}
%where $\theta_0, \theta_1$ are unknown functions from $[0, 1]$ to $
%over $[0,1]$, and $\stackrel{d}{=} $ stands for equality in
%distribution. Provided that $u\mapsto\theta_0(u)+
%(this imposes constraints on $\theta_0$ and $\theta_1$), the $\tau$th
%conditional quantile function
%of $Y_t$ indeed can be written as
% $Q^{Y_t| Y_{t-1}=y}_\tau
%= \theta_0(\tau) +\theta_1(\tau)y.$
%The same model also can be reformulated as the
%random-coefficient AR(1) model
%Y_t \stackrel{d}{=} \mu_0 +\alpha_{t } Y_{t-1} + V_t ,
%where $\mu_0 = {\rm E}[{\theta_0}(U_t)]$, $V_t = {\theta_0}(U_t) -
%with distribution function $ v\mapsto\theta_0
%^{-1} (v + \mu_0)$. Under this latter form, it is clear that quantile
%AR($p$) models bring us back to the location-scale paradigm, be it
%with randomized rescaling, while marginal and serial features remain
%intertwined.
If quantiles are to be considered, it seems intuitively reasonable that
the traditional $L_2$-tools, which are closely related with the
concepts of mean and variance, be abandoned in favor of
quantile-related ones. In particular, traditional $L_2$-projections
should be replaced with (weighted) $L_1$-projections. Recall that, in
traditional spectral analysis, estimation is usually based on the {\it
ordinary periodogram}
\[
I_n(\omega_{j,n}):=\frac{1}{n} \Biggl|\sum
_{t=1}^n Y_t \mathrm{e}^{-{\rm i}t\omega_{j,n}}
\Biggr|^2,
\]
where { $ \omega_{j,n} = {2\uppi j}/{n}\in{\mathcal F}_n:=\{ {2\uppi
j}/{n}| j = 1, \ldots, \lfloor\frac{n-1}{2}\rfloor- 1,
\lfloor\frac{n-1}{2}\rfloor\}$ }
denote the positive {\it Fourier frequencies}. A straightforward
calculation shows that this can be expressed as
\[
I_n(\omega_{j,n}) = \frac{n}{4}\bigl\|\hat{\mathbf
b}_{n,\mathrm
{OLS}}(\omega_{j,n})\bigr\|^2 :=
\frac{n}{4} \hat{\mathbf b}^{\prime
}_{n,\mathrm{OLS}} (
\omega_{j,n})\pmatrix{1&{\rm i}
\cr
-{\rm i}&1 } \hat{\mathbf b}_{n,\mathrm{OLS} }(
\omega_{j,n}),
\]
where $\| \cdot\|$ denotes the euclidian norm, and
%
%e2.4 #&#
\begin{equation}
\label{L2regre} \bigl( \hat{a}_{n,\mathrm{OLS} }(\omega_{j,n}) , \hat{
\mathbf b}^{\prime
}_{n,\mathrm{OLS}}(\omega_{j,n})\bigr):=
\operatorname{Argmin}\limits_{(a,
{\mathbf b}\pr)\in\rr^3}\sum_{t=1}^n
\bigl(Y_t -\bigl(a,\mathbf {b}^\prime\bigr)
\mathbf{c}_t (\omega_{j,n}) \bigr)^2
\end{equation}
is the ordinary least squares estimator in the linear model with
{regressors} $\mathbf{c}_t(\omega_{j,n}):=(1, \cos(t\omega_{j,n}) ,
\sin(t\omega_{j,n}))\pr$, corresponding to an $L_2$-projection of
the observed series onto the harmonic basis.

If, instead of a representation of $Y_t$ itself, we are interested in a
representation, in terms of the harmonic basis, of $Y_t$'s quantile of
order $\tau$, it may seem natural to replace that ordinary
periodogram $I_n(\omega_{j,n})$ with
\[
\Lnryh{n,\tau} (\omega_{j,n}) := \frac{n}{4}\bigl\|\beyh{n,\tau }(
\omega_{j,n})\bigr\|^2 := \frac{n}{4} \beyhpr{n,\tau} (
\omega _{j,n})
\pmatrix{1&{\rm i}
\cr
-{\rm i}&1  } %
  \beyh{n,\tau}(\omega_{j,n}),
\]
where
%
%e2.5 #&#
\begin{equation}
\label{L1regre} \bigl(\aeyh{n,\tau}(\omega_{j,n}), \beyh{n,\tau}(
\omega_{j,n})\bigr) := \operatorname{Argmin}\limits_{(a, \mathbf{b}\pr)\in\rr^3}\sum
_{t=1}^n \rho_\tau
\bigl(Y_t -\bigl(a,\mathbf{b}^\prime\bigr)
\mathbf{c}_t (\omega_{j,n}) \bigr),
\end{equation}
and
\[
\rho_\tau(x):= x\bigl(\tau- I\{ x \leq0\}\bigr) = (1-\tau) | x
| I\{ x\leq0\} + \tau| x| I\{ x > 0\}, \quad\quad\tau\in(0,1),
\]
is the so-called {\it check function} (see Koenker \cite{Koenker2005}). In
definition (\ref{L1regre}), the $L_2$-loss
function, which yields the classical
definition (\ref{L2regre}), is thus replaced by Koenker and Bassett's
weighted $L_1$-loss which produces
quantile regression estimates -- see Koenker and
  Bassett \cite{KoenkerBassett1978}. That
this indeed is a sensible definition will follow from the asymptotic
results of Section~\ref{321315458}.
%Call the collection of $L^\tau_n(\omega_{j,n})$'s, for $\omega_{j,n}$
%ranging over the Fourier frequencies, the {\it Laplace periodogram}.

This $L_1$-approach has been taken by Li \cite{Li2008} for the
particular case $\tau= 1/2$, leading to a least absolute deviations
(LAD) regression coefficient $ \beyh{n,0.5}$ and
{later} by Li \cite{Li2011} {for arbitrary $\tau\in(0,1)$}.
More generally, for a given series $Y_1,\ldots, Y_n$,
define the {\it Laplace periodogram kernel} as the collection
%
%e2.6 #&#
\begin{equation}
\label{lperio} \hspace*{-10pt}\Lnryh{n,\tau_1,\tau_2}(
\omega_{j,n}) := \frac{n}{4} {\beyhpr{n,\tau_1}} (
\omega_{j,n})  %
\pmatrix{1&{\rm i}
\cr
-{\rm i}&1  } %
  \beyh{n,\tau_2} (
\omega_{j,n}),\quad\quad \omega_{j,n} \in{\mathcal F}_n, (
\tau_1, \tau_2) \in(0,1)^2.
\end{equation}
For any $(\tau_1,\tau_2, \omega_{j,n} )$, computation of $\Lnryh
{n,\tau_1,\tau_2}(\omega_{j,n})$ is immediate via the simplex
algorithm (as in classical Koenker--Bassett quantile regression,
see Koenker \cite{Koenker2005}).

Similarly, define the {\it copula periodogram kernel} as the Laplace
periodogram kernel  $\Lnch{n,\tau_1,\tau_2}(\omega_{j,n})$
associated with the series $U_1,\ldots, U_n$. This means that $\Lnch
{n,\tau_1,\tau_2}(\omega_{j,n})$ is obtained from \eqref{lperio}
by replacing the estimate $ \beyh{n,\tau} $ by {the second and third
components of the vector}
\[
\bigl(\hat a, \bigl( \mathbf{\hat b}^U\bigr)\pr\bigr) :=
\operatorname{Argmin}\limits_{(a, \mathbf{b}\pr)\in\rr^3}\sum_{t=1}^n
\rho_\tau \bigl(U_t -\bigl(a,\mathbf{b}^\prime
\bigr) \mathbf{c}_t (\omega_{j,n}) \bigr).
\]
Finally, because the %marginal
distribution function $F$ required for the calculation of $U_t=F(Y_t)$
is not
known, we introduce
the {\it empirical} or {\it rank-based %copula %\utL_n^{\tau_1\tau_2}
Laplace periodogram kernel} as the Laplace periodogram kernel $\Lrh
{n,\tau_1,\tau_2}(\omega_{j,n})$ associated
with the series $n^{-1}R\n_1,\ldots, n^{-1}R\n_n$, where $R\n_t$
denotes the rank of $Y_t$ among $Y_1,\ldots, Y_n$.
In other words, $\Lrh{n,\tau_1,\tau_2}(\omega_{j,n})$
is obtained from \eqref{lperio}
by replacing the estimate $ \beyh{n,\tau} $ by {the second and third
components of the vector}
\[
(\hat a, {\uut{\mathbf{\hat b\pr}}} ) := \operatorname{Argmin}\limits_{(a, \mathbf{ b}\pr)\in\rr^3}
\sum_{t=1}^n \rho_\tau
\bigl(n^{-1}R\n_t -\bigl(a,\mathbf{b}^\prime\bigr)
\mathbf{c}_t (\omega_{j,n}) \bigr).
\]

A few remarks about the notation used in this paper are in order.
With $\hat T$, we usually denote a statistic obtained from the original
series $Y_1,\ldots,Y_n$, such as
$\Lnryh{n,\tau_1,\tau_2}$. The notation $\hat{T}^U$ means that
$\hat T$ has been computed from the probability integral transform
$U_1,\ldots,U_n$ of the data %(recall that $U_t=F_Y(Y_t)$)
-- a typical example is $\Lnch{n,\tau_1,\tau_2}$. Finally, the
notation $\uut{\hat{T}}$ reflects the fact that $\hat T$ has been
computed from the
normalized ranks $n^{-1}R\n_1,\ldots, n^{-1}R\n_n$ (see, for
instance, the rank-based Laplace
periodogram kernel $\Lrh{n,\tau_1,\tau_2}$).

%s3 #&#
\section{Asymptotic properties}\label{321315458}

%s3.1 #&#
\subsection{Asymptotics of Laplace and copula periodogram
kernels}\label{sec31}

We now proceed to deriving the asymptotic distributions of the Laplace
and rank-based
Laplace periodogram kernels, which, as we shall see, establishes their
relation to the spectral density kernels defined in \eqref{sd1} and
\eqref{sd2}. Throughout the rest of the paper, we make the following
basic assumptions.\def\theAssumption{(A1)}

%as1 #&#
\begin{Assumption}\label{(A1)}  The process $\{Y_t\}_{t\in\Z}$ is
strictly stationary and $\beta$-mixing, such that
\[
\beta(n) := \sup_{k \geq1} \mathrm{E} \sup_{B \in\cF
_{n+k}^{\infty}}
\bigl|\mathrm{P}\bigl(B|\cF_{-\infty}^k\bigr) - \mathrm{P}(B)\bigr| = \mathrm{O}
\bigl(n^{-\delta}\bigr), \qquad\delta> 1, \mbox{ as $n \to\infty$,}
\]
where $\cF_l^m := \sigma(Y_l, \ldots, Y_m)$ denotes the $\sigma
$-field generated by $Y_l, \ldots, Y_m$.
\end{Assumption}

The class of $\beta$-mixing processes is well studied, and contains a
wide range of linear and nonlinear processes, including (possibly,
under mild additional assumptions)
%There is a wide range of linear and nonlinear processes that are $
%Examples include
ARMA, general nonlinear scalar ARCH, threshold ARCH, and exponential
ARCH processes (see Liebscher \cite{Liebscher2005}), $\operatorname{GARCH}(p,q)$ processes with
moments (see Boussama \cite{Boussama1998}) and $\operatorname{GARCH}(1,1)$ processes with no
assumptions regarding the moments (see Francq and
  Zako\"{i}an \cite{FrancqZakodian2006}),
generalized polynomial random coefficient vector autoregressive
processes, and a family of generalized hidden Markov processes (Carrasco and Chen \cite
{CarrascoChen2002}) which includes stochastic volatility ones.\def\theAssumption{(A2)}

%as2 #&#
\begin{Assumption}\label{(A2)} The distribution function $F$ of
$Y_t$ and the joint distribution functions $F_k$ of $(Y_t, Y_{t+k})$
are twice continuously differentiable, with uniformly (with respect to
their arguments, and also with respect to $k$) bounded derivatives.
Moreover, there exists a subset $T$ of $[0,1]$ and, for every $\tau\in
T$, a positive real $d_{\tau}$, such that
$\inf_{| x-q_{\tau}|\leq d_{\tau}}f(x) > 0$, %\]
where $f$ and $q_{\tau} := F^{-1}(\tau)$ denote the density
and $\tau$-quantile corresponding to the distribution function $F$.
\end{Assumption}

Denote by $\Lnryh{n,\tau_1,\tau_2}$ and $\Lnch{n,\tau_1,\tau_2}$,
respectively, the Laplace and copula periodogram kernels associated
with a realization of length $n$. For each $(\tau_1 , \tau_2)\in
(0,1)^2$ and $\omega\in(0, \uppi)$, write
%
%e3.1 #&#
\begin{equation}
\label{scaled} \efs{\tau_1, \tau_2}(\omega):= {
\ef{q_{\tau_1}, q_{\tau
_2}}(\omega)}/\bigl(f(q_{\tau_1})f(q_{\tau_2})
\bigr)
\end{equation}
for the\vspace*{-3pt} \textit{scaled} version of the spectral density kernel $\ef
{q_{\tau_1}, q_{\tau_2}}(\omega) $ defined in \eqref{sd2}.
In the following two statements, $\stackrel{\mathcal
L}{\longrightarrow}$ stands for convergence in distribution, and $\chi
_k^2$ denotes
a chi-square distribution with $k$ degrees of freedom. We also
introduce the piecewise constant function (defined on
the interval $(0,\uppi)$)
%
%e3.2 #&#
\begin{equation}
\label{gfkt} g_n(\omega):= \omega_{j,n},
\end{equation}
where $\omega_{j,n}$ is the Fourier frequency closest to $\omega
$~-- more precisely, $\omega_{j,n}$ is such that $\omega$ belongs
to $ (\omega_{j,n}-\frac{2\uppi}{n},  \omega_{j,n}+\frac{2\uppi
}{n}]$. {The following
result is the key for understanding the asymptotic properties of the
Laplace periodogram kernel.}

%th3.1 #&#
\begin{theorem}\label{PropBAsymp}
Let $\Omega:= \{\omega_1, \ldots, \omega_\nu\} \subset(0,\uppi)$
and $T := \{\tau_1, \ldots, \tau_p\} \subset(0,1)$ denote distinct
frequencies and distinct quantile orders, respectively. Let Assumptions
\textup{\ref{(A1)}} and \textup{\ref{(A2)}} be satisfied with \textup{\ref{(A2)}} holding for every $\tau\in T$. Then
\[
\sqrt{n} \bigl( \beyh{n,\tau}\bigl(g_n(\omega)\bigr)
\bigr)_{\tau\in T,
\omega\in\Omega} \operatorname{\stackrel{\mathcal{L}} {\longrightarrow}}
\limits_{n
\rightarrow\infty}
\bigl( N_\tau(\omega) \bigr)_{\tau\in T,
\omega\in\Omega},
\]
where $( N_\tau(\omega) )_{\tau\in T,  \omega\in\Omega}$ denotes
a Gaussian random vector with mean zero and covariance
%
%e3.3 #&#
\begin{eqnarray}
\label{asymvar} M_{\tau_{1},\tau_{2}}^{\omega_1, \omega_2}& := &\Cov\bigl(N_{\tau
_1}(
\omega_1), N_{\tau_2}(\omega_2)\bigr) %:=\mathrm{E}[\tilde{
%))^\prime]
\nonumber
\\[-8pt]\\[-8pt]
&=& \cases{ 4\uppi \pmatrix{ \Re{\efs{\tau_1,
\tau_2}(\omega)} & \Im{\efs{\tau_1, \tau _2}(
\omega)}
\cr
-\Im{\efs{\tau_1, \tau_2}(\omega)} & \Re{
\efs{\tau_1, \tau _2}(\omega)} } & \quad\mbox{if $
\omega_1 = \omega_2 =: \omega$,}\vspace*{3pt}
\cr
\pmatrix{ 0 & 0
\cr
0
& 0 } & \quad\mbox{if $\omega_1 \neq\omega_2$.} }\nonumber
\end{eqnarray}
\end{theorem}

\begin{pf} The proof consists of two basic steps which we only
sketch here.
Details are provided in Appendix~\ref{appA}.

\noindent\textit{Step 1.} The first step consists of a linearization
of the estimate $\beyh{n,\tau}(\omega_{j,n})$
defined in \eqref{L1regre}. To be precise,
for any $\tau\in(0,1)$, $\omega\in(0,\uppi)$, and $\bolds\delta\in
\mathbb{R}^3$, let
%
%e3.4 #&#
\begin{equation}
\label{barZ} \Zyh{n,\tau, \omega}(\bdelta) := \sum
_{t=1}^n \bigl(\rho_{\tau}
\bigl(Y_{t} - q_\tau- n^{-1/2}
\mathbf{c}_{t}^\tr(\omega)\bdelta\bigr) -
\rho_{\tau}(Y_{t} - q_\tau ) \bigr),
\end{equation}
where $\mathbf{c}_{t}(\omega)% = (c_{t,1}(\omega), c_{t,2}(\omega),
%c_{t,3}(\omega))^\tr
:= (1, \cos(\omega t), \sin(\omega t))^\tr$, and $q_\tau$ denotes
the $\tau$-quantile of $F$.
Further define
\[
\Zy{n, \tau, \omega}(\bdelta) := -\bdelta^\tr\bzetay_{n,\tau,
\omega}
+ \tfrac{1}{2} \bdelta^\tr\bQy_{n,\tau,\omega} \bdelta,
\]
where
%
%e3.5 #&#
\begin{equation}
\label{xinx} \bzetay_{n,\tau,\omega}  :=  n^{-1/2} \sum
_{t=1}^n \mathbf {c}_{t}(\omega) \bigl(
\tau- I\{Y_t \leq q_{\tau}\}\bigr),
\end{equation}
and
%
%e3.6 #&#
\begin{equation}\label{qnx}
\bQy_{n, \tau, \omega}  :=  f(q_{\tau}) {n^{-1}}\sum
_{t=1}^n \mathbf{c}_{t}(\omega)
\mathbf{c}_{t}^\tr(\omega).
\end{equation}

We first show that the minimizers
%
%e3.7 #&#
\begin{equation}
\label{deltaquer} \deyh{n,\tau,\omega} := \arg\min_\bdelta\Zyh{n,
\tau,\omega }(\bdelta) \quad\mbox{and}\quad \dey{n,\tau,\omega} := \arg\min
_\bdelta\Zy{n,\tau,\omega }(\bdelta) = (\bQy_{n,\tau,\omega})^{-1}
\bzetay_{n,\tau,\omega}
\end{equation}
are close in probability (uniformly with respect to $\omega\in
\mathcal{F}_n$).
Note that, from the definition in \eqref{L1regre}, it follows that the
random variable $\sqrt{n} \beyh{n,\tau}(\omega_{j,n})$
coincides with the second and third components of the vector $\deyh
{n,\tau,\omega}$.
Moreover, for
$\omega_{j,n} =
{2\uppi j}/{n}$, we have
%
%e3.8 #&#
\begin{equation}
\label{propq} \bQy_{n,\tau,\omega_{j,n}} = f(q_{\tau}) \operatorname{diag}
(1,1/2,1/2),
\end{equation}
where $\operatorname{diag}(a_1,\ldots,a_k) $ denotes the diagonal matrix with
diagonal elements $a_1,\ldots,a_k$.
More precisely, we establish the following bound
%
%e3.9 #&#
\begin{eqnarray}
\label{thm:uniflinearworank}  \sup_{\omega\in\mathcal{F}_n}\|
\deyh{n,\tau,\omega}-\dey {n,
\tau,\omega}\|& =& \mathrm{O}_{\mathrm{P}} \bigl( r_n(\delta) \bigr),\nonumber
\\[-8pt]\\[-8pt]
r_n(\delta) &:=& \bigl(n^{-1/8} \log n\bigr) \vee
\bigl(n^{(1/4)  (1-\delta)
/(1+\delta) } (\log n)^{3/2}\bigr).\nonumber
\end{eqnarray}

This result is obtained from the following arguments, for which the
details are provided in Section~\ref{details}.
Roughly speaking, bounds of the type \eqref{thm:uniflinearworank}
can be obtained by showing that the corresponding functions $\Zyh{n,
\tau, \omega}$ and $ \Zy{n, \tau, \omega}$
are uniformly close in probability. A precise statement is given in
Lemma~\ref{lem:convrate} (see Section~\ref{subsub121}), where we show
that \eqref{thm:uniflinearworank} follows
if the bound
%
%e3.10 #&#
\begin{equation}
\label{conslemm1}  \sup_{\omega\in\mathcal{F}_n}\sup_{\|\bdelta- \dey{n,\tau
,\omega}\| \leq\epsilon} \bigl|
\Zyh{n, \tau, \omega}(\bdelta) - \Zy {n, \tau, \omega}(\bdelta)\bigr| =
\mathrm{O}_{\mathrm{P}} \bigl( r_n(\delta)^2 \bigr)
\end{equation}
can be established for some $\epsilon>0$.

Note that
\begin{eqnarray*}
&&\mathrm{P} \Bigl( \sup_{\omega\in\mathcal{F}_n}\sup_{\|\bdelta-
\dey{n,\tau,\omega}\| \leq\epsilon} \bigl|
\Zyh{n, \tau, \omega }(\bdelta) - \Zy{n, \tau, \omega}(\bdelta)\bigr| >
r_n(\delta)^2 \Bigr)
\\
&&\quad\leq\mathrm{P} \Bigl(\sup_{\omega\in\mathcal{F}_n} \sup_{\|
\bdelta\| \leq\epsilon+ \|\dey{n,\tau,\omega}\|}
\bigl|\Zyh{n, \tau, \omega}(\bdelta) - \Zy{n, \tau, \omega}(\bdelta)\bigr| >
r_n(\delta )^2, \sup_{\omega\in\mathcal{F}_n} \|\dey{n,
\tau,\omega}\| \leq A\sqrt{\log n} \Bigr)
\\
&&\quad\quad{}+ \mathrm{P} \Bigl( \sup_{\omega\in\mathcal{F}_n} \sup_{\|\bdelta
\| \leq\epsilon+ \|\dey{n,\tau,\omega}\|}
\bigl|\Zyh{n, \tau, \omega }(\bdelta) - \Zy{n, \tau, \omega}(\bdelta)\bigr| >
r_n(\delta)^2, \sup_{\omega\in\mathcal{F}_n} \|\dey{n,
\tau,\omega}\| > A\sqrt {\log n} \Bigr)
\\
&&\quad\leq\mathrm{P} \Bigl(\sup_{\omega\in\mathcal{F}_n}\sup_{\|\bdelta
\| \leq\epsilon+ A\sqrt{\log n}}
\bigl|\Zyh{n, \tau, \omega}(\bdelta) - \Zy{n, \tau, \omega}(\bdelta)\bigr| >
r_n(\delta)^2 \Bigr) + \mathrm {P} \Bigl(\sup
_{\omega\in\mathcal{F}_n} \|\dey{n,\tau,\omega}\| > A\sqrt{\log n} \Bigr).
\end{eqnarray*}
By application of Lemma~\ref{lem:boundelta}, it is therefore
sufficient to show that, for an enlarged $A$,
%
%e3.11 #&#
\begin{equation}
\label{conslemm2} \sup_{\omega\in\mathcal{F}_n}\sup_{\|\bdelta\| \leq A\sqrt{\log
n}} \bigl|
\Zyh{n, \tau, \omega}(\bdelta) - \Zy{n, \tau, \omega }(\bdelta)\bigr| =
\mathrm{O}_{\mathrm{P}} \bigl( r_n(\delta)^2 \bigr)
\end{equation}
and \eqref{conslemm1}, hence also, in view of Lemma~\ref{lem:convrate}, \eqref
{thm:uniflinearworank} is proved.
The proof of \eqref{conslemm2} is deferred to Section~\ref{subsub12x}.

\noindent\textit{Step 2.} {As we have discussed at the beginning of the first
step, the asymptotic properties of $ \sqrt{n} \beyh{n,\tau}(\omega
_{j,n})$ can be obtained from those of the random variables $\bdelta
_{n,\tau,\omega}$ for which an explicit expression is available.}
More precisely, for given sets $\Omega:= \{\omega_1, \ldots, \omega
_\nu\} \subset(0,\uppi)$ of  frequencies and
$T := \{\tau_1, \ldots, \tau_p\} \subset(0,1)$, consider the linear
combination with coefficients ${\bolds\lambda}_{ik} \in\IR^2$,
$i=1,\ldots, \nu$, $k=1,\ldots, p$
%
%e3.12 #&#
\begin{eqnarray}
\label{eqn:PropLlin} %\begin{split}
&&\sum_{k=1}^p
\sum_{i=1}^{\nu} {\bolds\lambda}^{\prime}_{ik}
\sqrt{n} \beyh{n,\tau_k}\bigl(g_n(\omega_i
)\bigr)\nonumber\\[-8pt]\\[-8pt] % \bar\bdelta_{n, \tau_k, g_n(\omega_i)}
&&\quad = \sum_{k=1}^p \sum
_{i=1}^{\nu} {\bolds\lambda}^{\prime}_{ik}
\sum_{t=1}^n \frac{2}{f(q_{\tau_k})}
\frac{\mathbf{v}_{tn}(\omega
_i)}{\sqrt{n}} \bigl(\tau_k - I\{Y_t \leq
q_{\tau_k}\}\bigr) + \mathrm{o}_{\mathrm{P}}(1)\nonumber, %\end{split}
\end{eqnarray}
where $\mathbf{v}_{tn}(\omega) := (\cos(g_n(\omega)t), \sin
(g_n(\omega)t))^{\prime}$.
The first equality is a consequence of \eqref{deltaquer}, \eqref
{propq} and \eqref{thm:uniflinearworank}.
Along the same lines as in the proof of Theorem~2 of Li \cite{Li2008},
and using the fact that \textup{\ref{(A1)}} implies $\sum_{k=-\infty}^{\infty}
|\gamma_k(q_{\tau_1}, q_{\tau_2})| \leq C \sum_{k=-\infty}^{\infty
} |k|^{-\delta} < \infty$, we obtain that
\[
\Cov\Biggl(\sum_{t=1}^n
\frac{2}{f(q_{\tau_{k_1}})} \frac{\mathbf
{v}_{tn}(\omega_{i_1})}{\sqrt{n}} \bigl(\tau_{k_1} - I
\{Y_t \leq q_{\tau
_{k_1}}\}\bigr), \sum
_{t=1}^n\frac{2}{f(q_{\tau_{k_2}})} \frac{\mathbf
{v}_{tn}(\omega_{i_2})}{\sqrt{n}}
\bigl(\tau_{k_2} - I\{Y_t \leq q_{\tau
_{k_2}}\}\bigr)
\Biggr)% \xrightarrow{} M_{\tau_{k_1},\tau_
%{k_2}}^{\omega_{i_1}, \omega_{i_2}},
\]
converges to $M_{\tau_{k_1},\tau_{k_2}}^{\omega_{i_1}, \omega
_{i_2}}$ defined in \eqref{asymvar}.
Hence, we have
\[
\Var \Biggl(\sum_{t=1}^n \sum
_{k=1}^p \sum_{i=1}^{\nu}
{\bolds \lambda}^{\prime}_{ik} \frac{2}{f(q_{\tau_k})}
\frac{\mathbf
{v}_{tn}(\omega_i)}{\sqrt{n}} \bigl(\tau_k - I\{Y_t \leq
q_{\tau_k}\}\bigr) \Biggr) \rightarrow\Var\Biggl(\sum
_{k=1}^p \sum_{i=1}^{\nu}
{\bolds\lambda}^{\prime}_{ik} N_{\tau
_k}(
\omega_i )\Biggr) .
\]

By an application of the central limit theorem for triangular arrays of
strongly mixing random variables in Francq and
  Zako\"{i}an \cite{FrancqZakoian2005}, with
$\kappa= 0$, $T_n = 0$, $r^{*} = (\delta-1)/(2+4\delta)$ and $\nu
^{*}=3/(\delta-1)$,
we deduce that
\[
\sum_{t=1}^n \sum
_{k=1}^p \sum_{i=1}^{\nu}
{\bolds\lambda }^{\prime}_{ik} \frac{2}{f(q_{\tau_k})}
\frac{\mathbf
{v}_{tn}(\omega_i)}{\sqrt{n}} \bigl(\tau_k - I\{Y_t \leq
q_{\tau_k}\}\bigr) \stackrel{\mathcal L} {\longrightarrow} \mathcal{N}
\Biggl( 0, \Var\Biggl(\sum_{k=1}^p \sum
_{i=1}^{\nu} {\bolds\lambda}^{\prime}_{ik}
N_{\tau_k}(\omega_i )\Biggr) \Biggr),
\]
where $( N_\tau(\omega) )_{\tau\in T,  \omega\in\Omega}$ denotes
a Gaussian random vector with mean zero and covariance matrix $\Cov
(N_{\tau_1}(\omega_1), N_{\tau_2}(\omega_2)) = M_{\tau_{k_1},\tau
_{k_2}}^{\omega_{i_1}, \omega_{i_2}}$.
Because of \eqref{eqn:PropLlin}, the quantity
\[
\sqrt{n} \sum_{k=1}^p \sum
_{i=1}^{\nu} {\bolds\lambda}^{\prime
}_{ik}
\beyh{\tau_k}\bigl(g_n(\omega_i )\bigr)
\]
converges in distribution to the same normal limit. Thus, it follows
from the traditional Cram\'er--Wold device
that
\[
\bigl(\sqrt{n} \beyh{n,\tau}\bigl(g_n(\omega)\bigr)
\bigr)_{\tau\in T,  \omega\in\Omega} \operatorname{\stackrel{\mathcal{L}} {\longrightarrow}}
\limits_{n
\rightarrow\infty}
\bigl( N_\tau(\omega) \bigr)_{\tau\in T,
\omega\in\Omega}.
\]
\upqed\end{pf}

 As an immediate consequence of the above result, the continuous
mapping theorem yields the asymptotic
distribution of a collection of Laplace periodogram kernels.

%th3.2 #&#
\begin{theorem}\label{PropLAsymp}
Under the assumptions of Theorem~\ref{PropBAsymp},
%
%e3.13 #&#
\begin{equation}
\label{LAsympt} \bigl(\Lnryh{n,\tau_1,\tau_2}
\bigl(g_n(\omega_1)\bigr),\ldots, \Lnryh{n,\tau
_1,\tau_2}\bigl(g_n(\omega_\nu)
\bigr)\bigr) \stackrel{\mathcal L} {\longrightarrow }\bigl(\Lnr{
\tau_1,\tau_2}(\omega_1),\ldots, \Lnr{
\tau_1,\tau _2}(\omega_\nu)\bigr),%
\end{equation}
where the random variables $\Lnr{\tau_1,\tau_2}$ associated with
distinct frequencies are mutually independent. Moreover,
%
%e3.14 #&#
\begin{equation}\label{tau1tau1}
\Lnr{\tau_1,\tau_2}(\omega)\sim \uppi\efs{
\tau_1, \tau_2}(\omega) \chi^2_2
\quad\quad\mbox{if } \tau_1=\tau_2,
\end{equation}
and
\[
\Lnr{\tau_1,\tau_2}(\omega)\stackrel{d} {=}
\frac{1}{4} (Z_{11}, Z_{12})
\pmatrix{ 1&{\rm i}
\cr
-{\rm i}&1 }\pmatrix{
 Z_{21}
\cr
Z_{22}  } %
 \quad\quad \mbox{if }
\tau_1\neq\tau_2,
\]
where
$(Z_{11}, Z_{12}, Z_{21}, Z_{22})\pr$
is a Gaussian vector with mean zero and covariance matrix
%
%e3.15 #&#
\begin{equation}
{\bolds\Sigma} _4(\omega):= 4\uppi\pmatrix{ {\efs{{\tau_1}, {\tau_1}}(\omega)}
&0&\Re{\efs{{\tau_1}, {\tau_2}}(\omega)} &\Im{\efs{{
\tau_1}, {\tau_2}}(\omega)}
\cr
0& {\efs{{\tau_1}, {\tau_1}}(\omega)} &-\Im{\efs{{
\tau_1}, {\tau_2}}(\omega)} &\Re{\efs{{
\tau_1}, {\tau_2}}(\omega)}
\cr
\Re{\efs{{\tau_1}, {\tau_2}}(\omega)} &-\Im{\efs{{
\tau_1}, {\tau_2}}(\omega)} & {\efs{{\tau_2},
{\tau_2}}(\omega)} &0
\cr
\Im{\efs{{\tau_1}, {\tau_2}}(\omega)} & \Re{\efs{{
\tau_1}, {\tau_2}}(\omega)} & 0 & {\efs{{
\tau_2}, {\tau_2}}(\omega)} }.\label{tau1tau2}
\end{equation}
\end{theorem}

It follows from Theorem~\ref{PropLAsymp}
that
$%\[
{\rm E}[\Lnr{\tau_1,\tau_2}(\omega)] = 2\uppi \efs{\tau_1, \tau
_2}(\omega)
$ %\]
for all $(\tau_1 , \tau_2)\in(0,1)^2$ and $\omega \in  (0, \uppi
)$, which indicates that an estimator of the scaled spectral density
$2\uppi  \efs{\tau_1, \tau_2}(\omega)$ defined in \eqref{scaled}
could be based on an average of quantities of the form $\Lnryh{n,\tau
_1,\tau_2}(\omega)$. Moreover, the following result, which is an
immediate consequence of Theorem~\ref{PropLAsymp}, yields the
asymptotic distribution of the
copula periodogram kernel.

%co3.1 #&#
\begin{Corol} \label{coroll}
Let $\Omega:= \{\omega_1, \ldots, \omega_\nu\} \subset(0,\uppi)$
denote distinct frequencies and $(\tau_1 , \tau_2) \in(0,1)^2$.
If Assumptions \textup{\ref{(A1)}}--\textup{\ref{(A2)}} hold for every $\tau\in\{\tau_1, \tau_2\}
$, then
%
%e3.16 #&#
\begin{equation}
\label{RBAsymptU} \bigl(\Lnch{n,\tau_1,\tau_2}
\bigl(g_n(\omega_1)\bigr),\ldots, \Lnch{n,\tau
_1,\tau_2}\bigl(g_n(\omega_\nu)
\bigr)\bigr) \stackrel{\mathcal L} {\longrightarrow}\bigl(\Lr{\tau_1,
\tau_2}(\omega _1),\ldots, \Lr{\tau_1,
\tau_2}(\omega_\nu)\bigr),
\end{equation}
where $g_n(\omega)$ is defined in \eqref{gfkt}.
The random variables $\Lr{\tau_1,\tau_2}$ in \eqref{RBAsympt}
associated with distinct frequencies are mutually independent,
%
%e3.17 #&#
\begin{equation}
\Lr{\tau_1,\tau_2}(\omega)\sim%\frac{1}{2} 2
\uppi
\ef{q_{\tau_1}, q_{\tau_2}}(\omega) %\frac{f_{q_{\tau_1} q_{\tau_2}}(\omega_j)}{f(q_{\tau_1})f(q_{\tau_2})}
\chi^2_2 %\text{and}
\quad\quad\mbox{if }
\tau_1=\tau_2, \label{RBtau1tau1}
\end{equation}
%
%if $ \tau_1=\tau_2$, \vspace{-2mm}
and
\[
\Lr{\tau_1,\tau_2}(\omega)\stackrel{d} {=}
\frac{1}{4} (Z_{11}, Z_{12})\pmatrix{ 1&{\rm i}
\cr
-{\rm i}&1 }\pmatrix{
 Z_{21}
\cr
Z_{22} } \quad\quad\mbox{if }
\tau_1\neq\tau_2,
\]
where
$(Z_{11}, Z_{12}, Z_{21}, Z_{22})\pr\sim{\mathcal N}(0, {\bolds
\Sigma} _4(\omega))$
with covariance matrix % \vspace{-3mm}
%e3.18 #&#
\begin{equation}
{\bolds\Sigma} _4(\omega):= 4\uppi\pmatrix{ {\ef{q_{\tau_1}, q_{\tau_1}}(\omega)} &0&\Re{
\ef{q_{\tau_1}, q_{\tau_2}}(\omega)} &\Im{\ef{q_{\tau
_1},
q_{\tau_2}}(\omega)}\vspace*{2pt}
\cr
0& {\ef{q_{\tau_1}, q_{\tau_1}}(\omega)} &-\Im{
\ef{q_{\tau_1}, q_{\tau_2}}(\omega)} &\Re{\ef{q_{\tau_1},
q_{\tau_2}}(\omega)}\vspace*{2pt}
\cr
\Re{\ef{q_{\tau_1}, q_{\tau_2}}(\omega)} &-\Im{
\ef{q_{\tau_1}, q_{\tau_2}}(\omega)} & {\ef{q_{\tau_2},
q_{\tau_2}}(\omega)} &0\vspace*{2pt}
\cr
\Im{\ef{q_{\tau_1}, q_{\tau_2}}(\omega)} & \Re{
\ef{q_{\tau_2}, q_{\tau_2}}(\omega)} & 0 & {\ef{q_{\tau_2},
q_{\tau_2}}(\omega)} }.\label{RBtau1tau2}
\end{equation}
\end{Corol}

In particular,
$%\[
{\rm E}[\Lnc{\tau_1,\tau_2}(\omega)] = 2\uppi\ef{q_{\tau
_1},q_{\tau_2}}(\omega)$. %\]
This indicates that the copula periodogram kernels $\Lnch{n,\tau
_1,\tau_2}$,
rather than the Laplace ones $\Lnryh{n,\tau_1,\tau_2}$, are likely
to be the appropriate tools for
statistical inference about $\ef{q_{\tau_1},q_{\tau_2}}$.
Unfortunately, they are not statistics, since they involve
the unknown marginal distribution $F$ which in practice is unspecified.
This problem is taken care of in the next section.

%s3.2 #&#
\subsection{Asymptotics of rank-based Laplace periodogram kernels}\label{s3.2}

The final result of this section establishes the asymptotic equivalence
of the copula and rank-based Laplace
periodogram kernels $\Lnch{n,\tau_1\tau_2}(\omega)$ and $\Lrh
{n,\tau_1\tau_2}(\omega)$, where the latter do not involve $F$,
hence can be computed
from the data. In particular, the following  results show that $\beyrh
{n,  \tau}$, and $\Lrh{n,\tau_1,\tau_2}(\omega)$ are
asymptotically distribution-free with respect to the marginal
distribution of $Y_t$ in the sense that their asymptotic distributions
only depend on the process $\{ U_t\}_{t \in\IZ}$.

%th3.3 #&#
\begin{theorem}\label{PropCLAsymp}
Let $\Omega:= \{\omega_1, \ldots, \omega_\nu\} \subset(0,\uppi)$
and $T := \{\tau_1, \ldots, \tau_p\} \subset(0,1)$ denote distinct
frequencies and quantile orders, respectively. Let Assumptions
\textup{\ref{(A1)}}--\textup{\ref{(A2)}} be satisfied with \textup{\ref{(A2)}} holding for every $\tau\in T$. Then,
\[
\bigl( \beyrh{n, \tau}\bigl(g_n(\omega)\bigr) \bigr)_{\tau\in T,  \omega
\in\Omega}
\operatorname{\stackrel{\mathcal{L}} {\longrightarrow}}\limits
_{n
\rightarrow\infty} \bigl(
N^{U}_{\tau,\omega} \bigr)_{\tau\in T,
\omega\in\Omega},
\]
where $( N^{U}_{\tau,\omega} )_{\tau\in T,  \omega\in\Omega}$ is
a Gaussian random vector with mean zero and covariance matrix
\begin{eqnarray*}
M_{\tau_{1},\tau_{2}}^{\omega_1, \omega_2} &:= &\Cov\bigl(N^{U}_{\tau_1,
\omega_1},
N^{U}_{\tau_2, \omega_2}\bigr)
\\
&=& \cases{ 4 \uppi \pmatrix{ \Re{
\ef{q_{\tau_1}, q_{\tau_2}}(\omega)} & \Im{\ef{q_{\tau
_1},
q_{\tau_2}}(\omega)}\vspace*{2pt}
\cr
-\Im{\ef{q_{\tau_1}, q_{\tau_2}}(
\omega)} & \Re{\ef{q_{\tau
_1}, q_{\tau_2}}(\omega)} } \vspace*{3pt}& \quad\mbox{if $
\omega_1 = \omega_2 =: \omega$, and}
\cr
\pmatrix{ 0 & 0
\cr
0 & 0 } &\quad \mbox{if $\omega_1 \neq\omega_2$.} }
\end{eqnarray*}
\end{theorem}

At first glance, the fact that replacing the $U_t$'s with their ranks
does not have any impact on the asymptotic distribution of $ \beyrh{n,  \tau}(g_n(\omega))$ seems quite surprising: a
completely different phenomenon indeed typically occurs when estimating
a copula, see e.g.
Genest and Segers \cite{GenestSegers2010}. The explanation for this is that the Bahadur
representation for the vector $(\hat a,\uut{\hat{\mathbf{b}}})$ is (see the proof of
Theorem~\ref{PropCLAsymp}) of the very special form
\[
\sqrt{n}\bigl(\bigl(\hat a,\beyrh{}' \bigr)' -
(q_\tau,0,0)'\bigr) = \bigl(\mathbf{Q}_{n,
\omega}^{U}
\bigr)^{-1}n^{-1/2}\sum_{t=1}^n
\mathbf{c}_{t}(\omega) \bigl(\tau- I\{U_{t} \leq\tau\} +
F\bigl(\hat F_n^{-1}(\tau)\bigr) - \tau \bigr),
\]
where the matrix $\mathbf{Q}_{n, \omega}^{U} := \frac{1}{n}\sum_{t=1}^n \mathbf{c}_{t}(\omega)\mathbf{c}_{t}'(\omega)$ is
diagonal. The additional term\break $F(\hat F_n^{-1}(\tau)) - \tau$ comes
into play because we are using ranks to estimate the unknown marginals.
However, due to the fact that, for Fourier frequencies $\omega$, $\sum_{t=1}^n \cos(\omega t) = \sum_{t=1}^n \sin(\omega t) = 0$, this
effect is not present in the first-order expansion of $\beyrh{}$ and
thus does not influence its asymptotic distribution.

Together with the above result, the continuous mapping theorem then
yields the following result.

%th3.4 #&#
\begin{theorem} \label{PropRBAsymp}
Under the assumptions of Theorem~\ref{PropCLAsymp}
%
%e3.19 #&#
\begin{equation}
\label{RBAsympt} \bigl(\Lrh{n,\tau_1,\tau_2}
\bigl(g_n(\omega_1)\bigr),\ldots, \Lrh{n,
\tau_1,\tau _2}\bigl(g_n(
\omega_\nu)\bigr)\bigr) \stackrel{\mathcal L} {\longrightarrow}\bigl(
\Lr {\tau_1,\tau_2}(\omega_1),\ldots, \Lr{
\tau_1,\tau_2}(\omega_\nu )\bigr),
\end{equation}
where $g_n(\omega)$ and the distribution of the random variables $\Lr
{\tau_1,\tau_2}$ are defined in
\eqref{gfkt} and Corollary~\ref{coroll}, respectively.
\end{theorem}

\begin{pf*}{Proof of Theorem~\ref{PropCLAsymp}} Recall that $\hat F_{n}$
denotes the empirical distribution function of $Y_1,\ldots, Y_n$; let
$\mathbf{e}_1 := (1,0,0)^\tr$, $\bdelta= (\delta_1, \delta_2,
\delta_3)^{\tr}$,
and $U_{t} := F(Y_{t})$. We
introduce the functions
\begin{eqnarray*}
\Zyrh{n,\tau, \omega}(\bdelta) &:=& \sum_{t=1}^n
\bigl(\rho_{\tau
} \bigl(\hat F_{n}(Y_{t}) - \tau-
n^{-1/2} \mathbf{c}_{t}^\tr(\omega )\bdelta\bigr)
- \rho_{\tau}\bigl(\hat F_{n}(Y_{t}) - \tau\bigr)
\bigr),
\\
\Zuh{n,\tau, \omega}(\bdelta) &:=& \sum_{t=1}^n
\bigl(\rho_{\tau
} \bigl( U_{t} - \tau- n^{-1/2}
\mathbf{c}_{t}^\tr(\omega)\bdelta\bigr) -
\rho_{\tau}(U_{t} - \tau) \bigr) - \delta_1
\sqrt{n}\bigl(F\bigl(\hat F_{n}^{-1}(\tau)\bigr) - \tau
\bigr),
\\
\Zu{n,\tau, \omega}(\bdelta) &:=& -\bdelta^\tr \bigl( \bzetau
_{n,\tau,\omega} + \mathbf{e}_1^\tr\sqrt{n}\bigl(F\bigl(
\hat F_{n}^{-1}(\tau)\bigr) - \tau\bigr) \bigr) +
\frac{1}{2} \bdelta^\tr\bQu_{n,\omega} \bdelta,
\end{eqnarray*}
where $
\mathbf{Q}_{n, \omega}^{U} := n^{-1} \sum_{t=1}^n \mathbf
{c}_{t}(\omega) \mathbf{c}_{t}^\tr(\omega) \mbox{{ and }}
\bzetau_{n,\tau,\omega} := n^{-1/2}\sum_{t=1}^n \mathbf
{c}_{t}(\omega)  (\tau- I\{U_{t} \leq\tau\}  )$. %\]
If we can show that the difference $\Zyrh{n,\tau, \omega}(\bdelta
)-\Zu{n,\tau, \omega}(\bdelta) $ is uniformly small in
probability, a slight modification of the arguments developed in the
proof of Theorem~\ref{PropLAsymp} yields a uniform linearization of $
\deyrh{n,\tau,\omega} := \arg\min_\bdelta\Zyrh{n,\tau,\omega
}(\bdelta)$.
More precisely, we show that
%
%e3.20 #&#
\begin{equation}
\label{thm:uniflinearrank} \sup_{\omega\in\mathcal{F}_n}\bigl\|\deyrh{n,\tau,\omega}-\deu {n,
\tau,\omega}\bigr\| = \mathrm{O}_{\mathrm{P}} \bigl( n^{(1/8)
(1-\delta)/(1+\delta)} \log n \bigr),
\end{equation}
where $
\deu{n,\tau,\omega} := \arg\min_\bdelta\Zu{n,\tau, \omega
}(\bdelta) = (\bQu_{n,\omega})^{-1} ( \bzetau_{n,\tau,\omega}
+ \mathbf{e}_1 \sqrt{n}(F(\hat F_{n}^{-1}(\tau)) - \tau)  )
$. %\]
The asymptotic normality of the linearization $\deu{n,\tau,\omega} $
then follows by the same arguments as in Step (2) of the proof of
Theorem~\ref{PropBAsymp}; details are omitted for the sake of brevity.

In order to prove \eqref{thm:uniflinearrank}, we note that
Lemma~\ref{lem:convrate} in   Appendix~\ref{appA} also holds with $\Zyh{n,\tau
, \omega}(\bdelta)$, $\Zx{n,\tau, \omega}(\bdelta)$, $\dex
{n,\tau,\omega}$ and $\deyh{n,\tau,\omega}$
replaced by$\Zyrh{n,\tau, \omega}(\bdelta)$, $\Zu{n,\tau, \omega
}(\bdelta)$, $\deu{n,\tau,\omega}$ and
$\deyrh{n,\tau,\omega}$, respectively.
Therefore, it suffices to establish that, for some $\epsilon>0$,
%
%e3.21 #&#
\begin{equation}
\label{eqn:ZbarRvsZtildeU} \sup_{\omega\in\mathcal{F}_n}\sup_{\|\bdelta- \deu{n,\tau
,\omega}\| \leq\epsilon} \bigl|
\Zyrh{n,\tau, \omega}(\bdelta) - \Zu {n,\tau, \omega}(\bdelta)\bigr| =
\mathrm{O}_{\mathrm{P}} \bigl(n^{(1/4)(1-\delta)/(1+\delta)} (\log n)^2 \bigr).
\end{equation}
Note that $\deu{n,\tau,\omega}$ decomposes into a term containing $
\bzetau_{n,\tau,\omega}$, to which Lemma~\ref{lem:boundelta}
applies, and a term involving $\sqrt{n}(F(\hat F_{n}^{-1}(\tau)) -
\tau)$ which, for every $\tau$, converges in distribution, so that
$\mathrm{P}( \sqrt{n}(F(\hat F_{n}^{-1}(\tau)) - \tau) > A\sqrt {\log n})\rightarrow0$ for any $A>0$. Therefore, there exists a
constant $A$ such that
$
\lim_{n\to\infty}\mathrm{P} (\sup_{\omega\in\mathcal{F}_n}
\|\deu{n,\tau,\omega}\| > A\sqrt{\log n}  ) = 0$.
It follows that, in order to establish (\ref{eqn:ZbarRvsZtildeU}), we
may restrict to a supremum with respect to the set $\|\bdelta\| \leq
2A\sqrt{\log n}$. Knight's identity (Knight \cite{Knight1998}; see page~121
of Koenker \cite{Koenker2005}) yields% the decomposition
\[
\Zyrh{n,\tau, \omega}(\bdelta) = \Zyrh{n,\tau, \omega,1}(\bdelta) + \Zyrh{n,
\tau, \omega ,2}(\bdelta),
\]
where
\[
\Zyrh{n,\tau, \omega,1}(\bdelta) = -\bdelta^\tr n^{-1/2}
\sum_{t=1}^n \mathbf{c}_{t}(
\omega) \bigl(\tau- I\bigl\{U_{t} \leq F\bigl(\hat
F_{n}^{-1}(\tau)\bigr)\bigr\} \bigr),
\]
and
\[
%F_{n,Y}^{-1}(s+\tau))\} - I\{I\{U_{t,n} \leq F_{n,X}(\hat F_{n,Y}^{-1}(
\Zyrh{n,\tau,
\omega,2}(\bdelta) =  \sum_{t=1}^n \int
_0^{n^{-1/2}{\mathbf{c}_{t}^\tr(\omega)\bdelta}} \bigl( I\bigl\{U_{t} \leq F
\bigl(\hat F_{n}^{-1}(s+\tau)\bigr)\bigr\} - I\bigl
\{U_{t} \leq F\bigl(\hat F_{n}^{-1}(\tau)\bigr)
\bigr\} \bigr)\,\mathrm{d}s.
\]
A similar {representation} holds for $\Zuh{n,\tau, \omega}(\bdelta)$.
Now the proof of \eqref{eqn:ZbarRvsZtildeU} is a consequence of
the following two auxiliary results,
which are proved in Sections~\ref{subsub213}--\ref{subsub212}:
%
%e3.22 #&#
\begin{eqnarray}
\label{h3} &&\sup_{\omega\in\mathcal{F}_n}\sup_{\|\bdelta\| \leq A\sqrt
{\log n}} \Biggl|
\Zyrh{n,\tau, \omega,1}(\bdelta) - \bdelta^\tr n^{-1/2} \sum
_{t=1}^n \mathbf{c}_{t}(
\omega) \bigl(\tau- I\{U_{t} \leq\tau\}\bigr)
 \nonumber
\\[-8pt]\\[-8pt]
&&\hphantom{\sup_{\omega\in\mathcal{F}_n}\sup_{\|\bdelta\| \leq A\sqrt
{\log n}} \Biggl|}{}- \delta_1\sqrt{n}\bigl(F\bigl(\hat F_{n}^{-1}(
\tau)\bigr) - \tau\bigr) \Biggr| = \mathrm{O}_{\mathrm{P}} \bigl( n^{ (1/4)  (1-\delta)
/(1+\delta) } (\log
n)^2 \bigr)
\nonumber
\end{eqnarray}
and
%
%e3.23 #&#
\begin{eqnarray}
\label{h4}
&&\sup_{\omega\in\mathcal{F}_n}\sup_{\|\bdelta\| \leq A
\sqrt{\log n}} \Biggl|
\Zyrh{n,\tau, \omega,2}(\bdelta) - \sum_{t=1}^n
\int_0^{n^{-1/2}{\mathbf{c}_{t}^\tr(\omega)\bdelta}}
\bigl( I\{U_{t} \leq s+
\tau\} -I\{U_{t} \leq\tau\} \bigr) \,\mathrm{d}s \Biggr|\qquad\quad\nonumber
\\[-8pt]\\[-8pt]
&&\quad= \mathrm{O}_{\mathrm{P}} \bigl(n^{(1/4) (1-\delta)/(1+\delta)}
(\log n)^2 \bigr).\nonumber
\end{eqnarray}
Note that the combination of \eqref{h3} and \eqref{h4} implies that
$\Zyrh{n,\tau, \omega}$ and $\Zuh{n,\tau, \omega}$ are uniformly
close in probability.
Finally,
we obtain from \eqref{conslemm2} that
%
%e3.24 #&#
\begin{equation}
\label{h4a} \sup_{\omega\in\mathcal{F}_n}\sup_{\|\bdelta\|\leq A\sqrt{\log
n}} \bigl|
\Zuh{n,\tau,\omega}(\bdelta) - \Zu{n,\tau, \omega}(\bdelta )\bigr| = \mathrm{O}_{\mathrm{P}}
\bigl( r_n(\delta)^2 \bigr),
\end{equation}
where we may replace $\Zyh{n,\tau,\omega}(\bdelta)$ with $\Zuh
{n,\tau,\omega}(\bdelta)$ and $\Zy{n, \tau, \omega}(\bdelta)$
with $\Zu{n,\tau, \omega}(\bdelta)$, since $U_{1}, \ldots, U_{n}$
are $\beta$-mixing with the rate from \textup{\ref{(A1)}}, as required, and the
additional term
$\delta_1 \sqrt{n}(F(\hat F_{n}^{-1}(\tau)) - \tau)$ %\]
appears in both $\Zuh{n,\tau,\omega}(\bdelta) $ and $ \Zu{n,\tau,
\omega}(\bdelta)$. Combining \eqref{h3}--\eqref{h4a} yields \eqref
{eqn:ZbarRvsZtildeU}, thus completing the proof of Theorem~\ref{PropCLAsymp}.
\end{pf*}

%s4 #&#
\section{Smoothed periodograms}
\label{sec4}

% \subsection{Smoothed Laplace periodograms.}

We have seen in Section~\ref{sec31} that the Laplace periodogram
kernel, for all $(\tau_1, \tau_2)$, converges in distribution, and
that the expectation of the limit is the \textit{scaled} spectral
density kernel (at $(\tau_1, \tau_2)$)
\[
2\uppi\efs{\tau_1, \tau_2}(\omega) :=2\uppi
\frac{\ef{q_{\tau_1}
q_{\tau_2}}(\omega)}{f(q_{\tau_1})f(q_{\tau_2})} = \frac{1}{f(q_{\tau_1}) f(q_{\tau_2})} \sum_{k=-\infty}^{\infty}
\gamma_k(q_{\tau_1}, q_{\tau_2}) \mathrm{e}^{-\mathrm{i} \omega k}.
\]
%
%where
In practice, however, this is not enough, and consistent estimation is
a minimal requirement. For this purpose, we consider, as in traditional
spectral estimation, smoothed versions of our periodograms, of the form
%
%e4.1 #&#
\begin{equation}
\label{eqn:defDSAE} \efh{n,\tau_1,\tau_2} (
\omega_{j,n}) := \sum_{|k| \leq N_n}
W_n(k)\Lnryh{n,\tau_1,\tau_2} (
\omega_{j+k,n})
\end{equation}
at the Fourier frequencies $\omega_{j,n} = 2 \uppi j / n$,
where
%$(N_n)$ is a sequence of positive integers such that
$N_n \rightarrow\infty$ as $n\to\infty$ is a sequence of positive
integers, and $W_n = \{W_n(j)   \dvt   |j| \leq N_n\}$ is a sequence of
positive weights satisfying
\[
W_n(k) = W_n(-k) \quad\quad\mbox{for all $k$}\quad  \mbox{and} \quad
\sum_{|k| \leq N_n} W_n(k) = 1.
\]

Extending the definition of $\efh{n,\tau_1,\tau_2}$ to the interval
$(0,\uppi)$, we introduce
\[
\bigl\{(0,\uppi) \ni\omega\mapsto\efh{n,\tau_1,\tau_2}(
\omega) | (\tau_1, \tau_2)\in(0,1)^2
\bigr\}
\]
as the \textit{smoothed Laplace periodogram kernel},
where
%
%e4.2 #&#
\begin{equation}
\label{eqn:defDSAE1} \efh{n,\tau_1,\tau_2}(\omega) := \efh{n,
\tau_1,\tau_2} \bigl(g_n(\omega)\bigr),
\end{equation}
and the function $g_n$ is defined in (\ref{gfkt}). In order to show
that $\efh{n,\tau_1,\tau_2}(\omega) $ is a consistent
estimator of the scaled spectral density $\efs{\tau_1, \tau
_2}(\omega)$, we make the following additional assumptions.
% \begin{itemize}
% \item[(N)]
\def\theAssumption{(A3)}

%as3 #&#
\begin{Assumption}\label{(A3)}
$N_n / n \rightarrow0$, and
% \item[(W)]
$\sum_{|k| \leq N_n} W_n^2(k) = \mathrm{O}(1/n)$ as $n\to\infty$.
\end{Assumption}
% \item[(S1)]

\def\theAssumption{(A4)}

%as4 #&#
\begin{Assumption}\label{(A4)}
%For any $\tau_1, \tau_2 \in(0,1)$,
% \item[(S2)]
For any $\tau_1, \tau_2, \tau_3, \tau_4 \in(0,1)$,
\[
\sum_{k_2, k_3, k_4 =-\infty}^{\infty} \bigl|\cum\bigl(I
\{Y_t \leq q_{\tau
_1}\},I\{Y_{t+k_2} \leq
q_{\tau_2}\}, I\{Y_{t+k_3} \leq q_{\tau_3}\}, I
\{Y_{t+k_4} \leq q_{\tau_4}\}\bigr)\bigr| < \infty,
\]
where $\cum(\zeta_1, \ldots, \zeta_r) := \sum(-1)^{p-1} (p-1)!
(\mathrm{E} \prod_{j \in\nu_1} \zeta_j) \cdots(\mathrm{E} \prod_{j \in\nu_p} \zeta_j)$ (with summation extending over all
{partitions $\{\nu_1, \ldots, \nu_p\}$, $p=1,\ldots,r$ of $\{
1,\ldots,r\}$)} denotes the \textit{$r$th order joint cumulant} of
the random vector $(\zeta_1, \ldots, \zeta_r)$ (cf. Brillinger \cite
{Brillinger1975}, page~19).
\end{Assumption}
% \item[(F)]

\def\theAssumption{(A5)}

%as5 #&#
\begin{Assumption}\label{(A5)}
The functions
$\omega\mapsto\ef{q_{\tau_1}, q_{\tau_2}}$ defined in \eqref{sd2}
are continuously differentiable for all $(\tau_1,\tau_2)\in(0,1)^2$.
% \item[(M1)]
\end{Assumption}

Note that an assumption similar to \textup{\ref{(A4)}}, but with the cumulant of
$Y_t$'s instead of the cumulant of the indicators, is made when
consistency of smoothed cross-periodograms is proved,
%implied by 4-th order strong mixing (cf. \citeauthor{Rao1999} \cite{Rao1999}) with an
%appropriate rate,
and that \textup{\ref{(A5)}} follows if \textup{\ref{(A1)}} holds with $\delta> 2$, because this
implies
\[
\sum_{k \in\IZ} |k| \bigl|\gamma_k (
\tau_1, \tau_2)\bigr| < \infty.
\]

%th4.1 #&#
\begin{theorem}\label{41414125544254}
 Let \textup{\ref{(A1)}}--\textup{\ref{(A5)}} hold. Then the smoothed Laplace periodogram
defined in \eqref{eqn:defDSAE} and \eqref{eqn:defDSAE1} is a
consistent estimator for the scaled Laplace spectral density; more precisely,
%
%e4.3 #&#
\begin{equation}
\label{consistrate} \efh{n,\tau_1,\tau_2}(\omega) = 2\uppi
\efs{\tau_1, \tau_2}(\omega ) + \mathrm{O}_{\mathrm{P}}
\bigl(R_n + n^{-1/2} + N_n / n \bigr) = 2\uppi\efs {
\tau_1, \tau_2}(\omega) + \mathrm{o}_{\mathrm{P}}(1),
\end{equation}
where
$R_n = (n^{-1/8} (\log n)^{3/2}) \vee(n^{(1/4) (1-\delta)
/(1+\delta)} (\log n)^{9/4})$.
\end{theorem}

\begin{pf} The proof proceeds in several steps which are sketched
here -- technical details can be found in Appendix~\ref{appB}. We first show
(Section~\ref{secProof41det_1})
that
%
%e4.4 #&#
\begin{equation}\label{lem:ConsistencyA}
\Lnryh{n,\tau_1,\tau_2} (\omega_{j,n}) =
\Lnr{n,\tau_1,\tau_2} (\omega_{j,n})  /
\bigl(f(q_{\tau_1}) f(q_{\tau_2}) \bigr) + \mathrm{O}_{\mathrm{P}}(R_n),
\end{equation}
uniformly in the Fourier frequencies $\omega_{j,n} := 2\uppi j /n$, where
\begin{eqnarray*}
\Lnr{n,\tau_1,\tau_2}(\omega_{j,n}) &:=&
n^{-1} d_n(\tau_1, \omega _{j,n})
d_n(\tau_2, -\omega_{j,n}),
\\
d_n(\tau, \omega_{j,n}) &:=& \sum_{t=1}^n {\rm e}^{{\rm i} \omega
_{j,n} t} \bigl(\tau- I\{Y_t \leq q_{\tau}\}\bigr) = (1, \mathrm{i}) n \bey
{n,\tau,\omega_{j,n}}2^{-1}f(q_\tau) \quad \mbox{and}
\\
n^{1/2} \bey{n, \tau, \omega_{j,n}} &:=& \frac{2}{f(q_{\tau})}
n^{-1/2} \sum_{t=1}^n %
\pmatrix{ \cos(\omega_{j,n} t)
\cr
\sin(\omega_{j,n} t) } %
\bigl(\tau- I
\{Y_t \leq q_{\tau}\}\bigr).
\end{eqnarray*}
As an immediate consequence, we obtain
\[
\efh{n,\tau_1,\tau_2} (\omega_{ j,n }) = \sum
_{|k| \leq N_n} W_n(k) {\Lnr{n,
\tau_1,\tau_2}(\omega_{j+k,n})}/ \bigl(
f(q_{\tau
_1}) f(q_{\tau_2}) \bigr)+ \mathrm{O}_{\mathrm{P}}(R_n).
\]
In Section~\ref{secProof41det_2}, we show that, for any
$\omega_{j,n} = 2\uppi j / n$,
%
%e4.5 #&#
\begin{equation}
\label{lem:ConsistencyB} K_n := \sum_{|k| \leq N_n}
W_n(k) \biggl(\frac{\Lnr{n,\tau_1,\tau
_2} (\omega_{j+k,n})}{f(q_{\tau_1}) f(q_{\tau_2})} - \efs{\tau_1,
\tau_2}(\omega_{j+k,n}) \biggr) = \mathrm{O}_{\mathrm
{P}}(1/
\sqrt{n}).
\end{equation}
Now, let $\omega_{j_n n}$ be a sequence of Fourier frequencies such
that $|\omega_{j_n, n} - \omega| = \mathrm{O}(N_n / n)$ for some $\omega\in
(0,\uppi)$: both for $f \equiv\Re\efs{\tau_1, \tau_2}$ and $f
\equiv\Im\efs{\tau_1, \tau_2}$, we have
\begin{eqnarray*}
  \Biggl|\sum_{|k| \leq N_n} W_n(k) \bigl( f (
\omega_{j_n + k,n}) - f (\omega) \bigr) \Biggr|
&\leq&\sum
_{|k| \leq N_n} W_n(k) \bigl|f'(
\xi_{j_n+k,n})\bigr| \llvert \omega _{j_n + k,n} - \omega\rrvert
\\
&  \leq& C_n \sum_{|k| \leq N_n}
W_n(k) \llvert {2\uppi k}/{n} + \omega _{j_n n} - \omega\rrvert
\\
&\leq& C_n \sum_{|k| \leq N_n} W_n(k)
\llvert {2\uppi k}/{n}\rrvert + C_n \sum_{|k| \leq N_n}
W_n(k) \llvert \omega_{j_n n} - \omega \rrvert
\\
& \leq& C_n \bigl( 2 \uppi{N_n} {/n} + \llvert
\omega_{j_n n} - \omega \rrvert \bigr) \sum_{|k| \leq N_n}
W_n(k) = \mathrm{O}(N_n / n),
\end{eqnarray*}
where $|\xi_{j_n+k,n} - \omega| \leq|\omega- \omega_{j_n+k,n}|$
and $C_n := \sup_{\xi\in\Xi_n} |f'(\xi)|$
is the supremum over
\[
\Xi_n = \bigl\{ \xi | \omega- |\omega-\omega_{j_n,
n}| -
\omega_{N_n, n} \leq\xi\leq\omega+ |\omega-\omega_{j_n,
n}| +
\omega_{N_n, n} \bigr\}.
\]
Note that, since $|\omega-\omega_{j_n, n}| \rightarrow0$ and $\omega
_{N_n, n} = 2\uppi N_n / n \rightarrow0$, $C_n \rightarrow f'(\omega)$,
so that $(C_n)$ is a bounded sequence.
This yields
\[
\Biggl|\sum_{|k| \leq N_n} W_n(k) \bigl(\,\efs{
\tau_1, \tau_2} (\omega_{j_n+k}) - \efs{
\tau_1, \tau_2}(\omega) \bigr) \Biggr| = \mathrm{O}(N_n /
n),
\]
which completes the proof of Theorem~\ref{41414125544254}.
\end{pf}

For a consistent estimation of the (unscaled) Laplace spectral density
$\ef{\tau_1, \tau_2}(\omega)$,
we propose a smoothed version %discrete spectral average estimator
\[
\efth{n,\tau_1,\tau_2}(\omega) := \efth{n,
\tau_1,\tau_2}%(\omega)
\bigl(g_n(\omega)
\bigr),\quad\quad \efth{n,\tau_1,\tau_2} %(\omega)
(
\omega_{j,n}) := \sum_{|k| \leq N_n}
W_n(k) \Lrh{n,\tau_1,\tau_2}(
\omega_{j+k,n})
\]
of the rank-based Laplace periodogram $\Lrh{n,\tau_1,\tau_2}(\omega)$.
%: where $g_n(\omega)$ is the same as in (\ref{eqn:defDSAE1}).
We then have the following result.
%, which is proved in Section~\ref{secProof42}.
%$\utL_{n}^{\tau_1, \tau_2} (\omega_{j+k,n})$ are the Laplace
%periodograms calculated from the standardized ranks, is used.

%th4.2 #&#
\begin{theorem}\label{ConsistenctRankP}
Let Assumptions \textup{\ref{(A1)}}--\textup{\ref{(A5)}} hold. Then the smoothed rank-based Laplace
periodogram $\efth{n,\tau_1,\tau_2}$ is a consistent estimator of
the (unscaled) Laplace spectral density $\ef{q_{\tau_1}, q_{\tau
_2}}$. More precisely,
\[
\efth{n,\tau_1,\tau_2}%(\omega)
(\omega) = 2\uppi
\ef{q_{\tau_1}, q_{\tau_2}}(\omega) + \mathrm{O}_{\mathrm
{P}} \bigl(
n^{(1/8)(1-\delta)/(1+\delta)} (\log n)^{3/2} + {N_n}/{n} \bigr) = 2\uppi
\ef{q_{\tau_1}, q_{\tau_2}}(\omega) + \mathrm{o}_{\mathrm{P}}(1).
\]
\end{theorem}

\begin{pf}
The proof is very similar to that of Theorem~\ref
{41414125544254}. The main difference lies in the asymptotic
representation for the second and third coordinates $n^{1/2} \beuh
{n,\tau,\omega}$ of the quantity $\deu{n,\tau,\omega}$ in \eqref
{thm:uniflinearrank}. Here we use \eqref{thm:uniflinearrank}, which
implies that
\begin{eqnarray*}
&&\sup_{\omega\in\mathcal{F}_n} \Biggl\|n^{1/2} \beuh{n,\tau,\omega} -
2n^{-1/2} \sum_{t=1}^n
\pmatrix{ \cos(\omega t)
\cr
\sin(\omega t)  } %
\bigl(\tau- I\bigl\{F
(Y_t) \leq\tau\bigr\}\bigr) \Biggr\|
\\
&&\quad= \mathrm{O}_{\mathrm{P}}
\bigl(n^{(1/8)(1-\delta)/(1+\delta)} (\log n)^{3/2}\bigr).
\end{eqnarray*}
The rest of the proof follows as in the proof of Theorem~\ref
{41414125544254}, yielding the estimate
\[
\efth{n,\tau_1,\tau_2}%(\omega)
(\omega) = 2\uppi\ef{
\tau_1, \tau_2}(\omega) + \mathrm{O}_{\mathrm{P}}
\bigl(n^{(1/8) (1-\delta)/(1+\delta) } (\log n)^{3/2} + n^{-1/2} +
{N_n}/{n} \bigr).
\]
Finally, the assumptions imply that $n^{-1/2}
= \mathrm{O}(n^{(1/8)(1-\delta)/(1+\delta)} (\log n)^{3/2})$,
which completes the proof
of Theorem~\ref{ConsistenctRankP}.
\end{pf}

Note that Theorem~\ref{41414125544254} solves an important open
problem raised in Li \cite{Li2008,Li2011}, who considered the Laplace periodogram
$\Lnryh{n,\tau_1,\tau_2}$ for $\tau_1=\tau_2$.
This author established the asymptotic unbiasedness of a
smoothed version of the Laplace periodogram, but not its consistency. Moreover,
as pointed out in Theorem~\ref{PropBAsymp} the smoothed version of
$\Lnryh{n,\tau_1,\tau_2}$ is not consistent for the copula
spectral density kernel, which is the main object of interest in this
paper. Theorem~\ref{ConsistenctRankP} shows that the smoothed
rank-based Laplace periodogram
yields a consistent estimate of this quantity.

%s5 #&#
\section{Finite-sample properties}\label{sec5} %\setcounter{equation}{0}
%s5.1 #&#
\subsection{Simulation results}\label{secsim}

In order to illustrate the finite-sample properties of the new
estimates, we present a small
simulation study, where we consider four models. In Models (1) and (2),
the observations are AR(1) processes with $Y_t = -0.3 Y_{t-1} +
\varepsilon_t$, and $\mathcal{N}(0,1)$- and $t_1$-distributed
innovations $\varepsilon_t$. Note that in Model (2) no moments exist,
hence the traditional spectral density is not defined. Model (3) is a
QAR(1) model (cf. Koenker and Xiao \cite{Koenker2006}), that is, a model of the form
$Y_t = \theta_0(U_t) + \theta_1(U_t) Y_{t-1}$, where $(U_t)$ is a
sequence of i.i.d. standard uniform random variables and $\theta
_1$ and $\theta_0$ are functions from $[0,1] $ to $\IR$; more
specifically, we chose $\theta_1(u) = 1.9(u-0.5)$ and $\theta_0(u) =
0.1 \Phi^{-1}(u)$, with $\Phi$ the standard normal distribution
function. Model (4) is the $\operatorname{ARMA}(1,1)$ model $Y_t = -0.8 Y_{t-1} + 1.25
\varepsilon_{t-1} + \varepsilon_t$ with $\varepsilon_t \sim t_3$.
Note that this defines an {\it all-pass} $\operatorname{ARMA}(1,1)$ process where the
observations are uncorrelated, but not independent (cf. e.g., Breidt, Davis and
  Trindade \cite
{BreidtEtAl2001}). All results presented in this section are based on
$5000$ independent replications.

For each of those four models, we generated pseudo-random time series
of lengths $n=100$, $n=500$ and $n=1000$, and computed the Laplace and
rank-based Laplace periodogram
for $\tau_1, \tau_2 \in\{0.05, 0.25, 0.5, 0.75, 0.95\}$. We also
computed the smoothed estimates using Daniell kernels with parameters
$(2,1)$ for $n=100$, $(10,4)$ for $n=500$, and $(10,25)$ for $n=1000$
-- namely, the kernel $W_n^{(m_1, \ldots, m_p)}(j)$ recursively
defined, for parameters $(m_1, \ldots, m_p)$, with $N_n = \sum_{j=1}^p m_j < n/2$, as
\begin{eqnarray*}
W_n^{(m)}(j) & := & (2 m - 1)^{-1} I\bigl\{|j| \leq m
\bigr\},
\\
W_n^{(m_1, \ldots, m_p)}(j) & := & C \bigl(W_n^{(m_1, \ldots, m_{p-1})}
* W_n^{(m_p)}\bigr) (j)
\\
&= &C \sum_{|k| \leq m_p} (2 m_p -
1)^{-1} W_n^{(m_1, \ldots, m_{p-1})}(j-k),
\end{eqnarray*}
where $*$ denotes convolution of two kernels and the constant $C$ is
chosen such that $\sum_{|j| \leq N_n} W_n^{(m_1, \ldots, m_p)}(j) =
1$; the parameters $m_1$ and $m_2$, $N_n = m_1 + m_2$, were chosen by
empirical considerations.

\begin{figure}[t]

\includegraphics{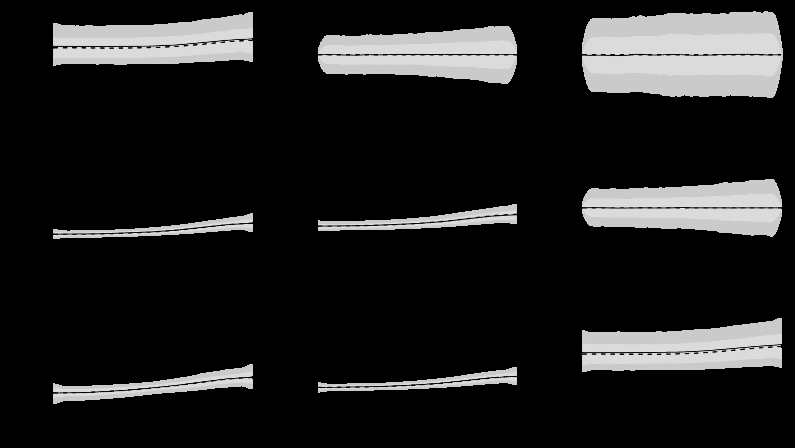}

    \caption{Smoothed Laplace periodograms and (scaled) spectral densities as defined in \protect\eqref{scaled} from 5000 replications of length
     $n=500$ of   an  AR(1) process with  $\mathcal{N}(0,1)$   innovations.
 All curves are plotted against $\omega / (2\uppi)$. Real parts (Imaginary parts) of the periodogram and
 spectrum are presented in subfigures with $\tau_2 \leq \tau_1$ ($\tau_2 > \tau_1$): the dashed line
 represents the actual scaled spectrum [cf.~\protect\eqref{scaled}], the solid line the (pointwise) mean of
 the simulated smoothed Laplace periodograms. The gray areas represent the $0.1$, $0.25$, $0.75$ and
 $0.9$ (pointwise) sample quantiles of the smoothed
 periodograms over the 5000 replications.}\label{fig:LPG-B3-M1}
\end{figure}

\begin{figure}

\includegraphics{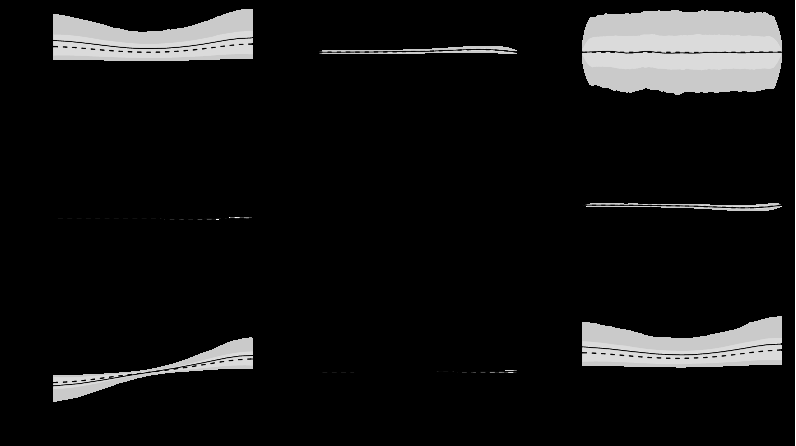}

    \caption{Smoothed Laplace periodograms and (scaled) spectral densities as defined in \protect\eqref{scaled}.
    The process is an  AR(1) process with  $t_1$  innovations and  the sample size is $500$.
    All curves are plotted against $\omega / (2\uppi)$. Real parts (Imaginary parts) of the
    periodogram and spectrum are presented in subfigures with $\tau_2 \leq \tau_1$ ($\tau_2 > \tau_1$):
     the dashed line represents the actual scaled spectrum [cf.~\protect\eqref{scaled}], the solid line the
     (pointwise) mean of the simulated smoothed Laplace periodograms. The gray areas represent the $0.1$,
     $0.25$, $0.75$ and $0.9$ (pointwise) sample quantiles of the smoothed periodograms over the 5000 replications.}\label{fig:LPG-B3-M2}
\end{figure}

\begin{figure}

\includegraphics{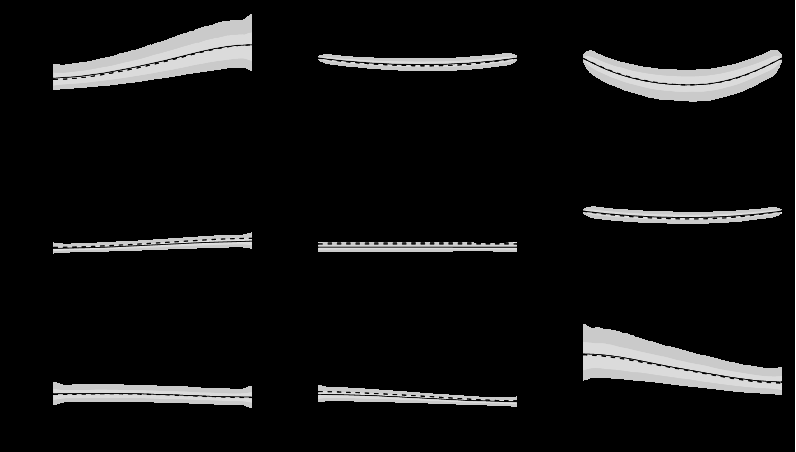}

    \caption{Smoothed Laplace periodograms and (scaled) spectral densities as defined in \protect\eqref{scaled}.
    The process is a  QAR(1) process with  $\theta_1(u) = 1.9(u-0.5)$, $\theta_0(u) = 0.1 \Phi^{-1}(u)$ and the sample size is~$500$.
    All curves are plotted against $\omega / (2\uppi)$. Real parts (Imaginary parts) of the periodogram and
    spectrum are presented in subfigures with $\tau_2 \leq \tau_1$ ($\tau_2 > \tau_1$): the dashed
    line represents the actual scaled spectrum [cf.~\protect\eqref{scaled}], the solid line the
    (pointwise) mean of the simulated smoothed Laplace periodograms. The
    gray areas represent the $0.1$, $0.25$, $0.75$ and $0.9$ (pointwise) sample quantiles
    of the smoothed periodograms over the 5000 replications.}\label{fig:LPG-B3-M3}
\end{figure}

\begin{figure}

\includegraphics{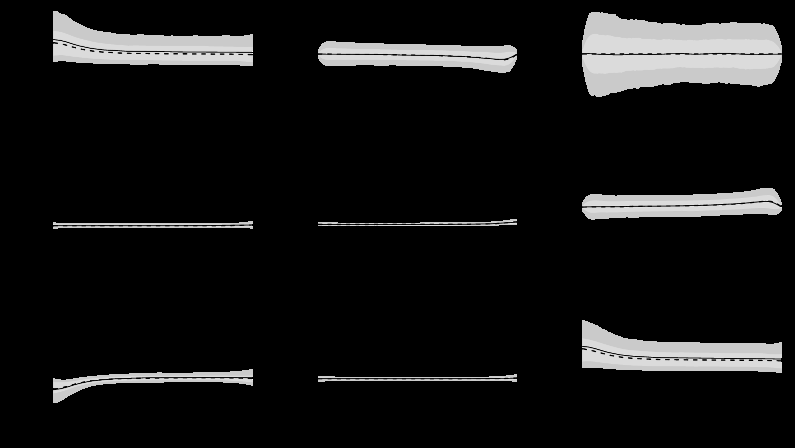}

    \caption{Smoothed Laplace periodograms and (scaled) spectral densities as defined in \protect\eqref{scaled}.
    The process is an $\operatorname{ARMA}(1,1)$ process with  $t_3$  innovations and  the sample size is~$500$.
    All curves are plotted against $\omega / (2\uppi)$. Real parts (Imaginary parts) of the periodogram
    and spectrum are presented in subfigures with $\tau_2 \leq \tau_1$ ($\tau_2 > \tau_1$):
    the dashed line represents the actual scaled spectrum [cf.~\protect\eqref{scaled}], the solid
    line the (pointwise) mean of the simulated smoothed Laplace periodograms. The gray areas
    represent the $0.1$, $0.25$, $0.75$ and $0.9$ (pointwise) sample quantiles of the smoothed periodograms over the 5000 replications.} \label{fig:LPG-B3-M4}
\end{figure}

\begin{figure}

\includegraphics{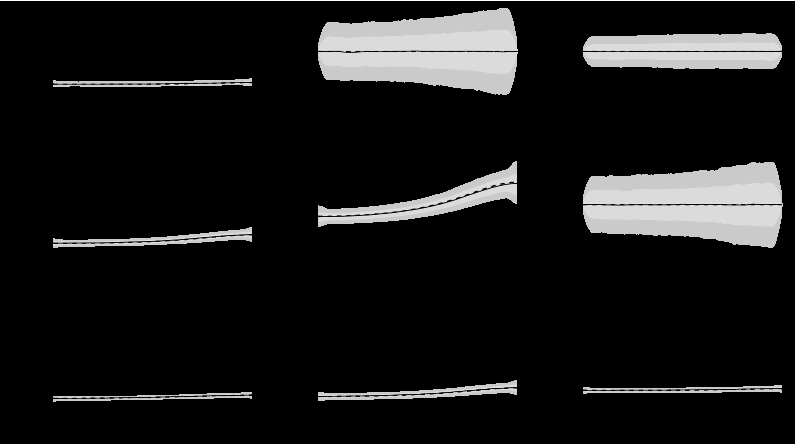}

    \caption{Smoothed rank-based  Laplace periodograms and (unscaled)
    spectral  densities as defined in~\protect\eqref{sd2}.
    The process is an  AR(1) process with  $\mathcal{N}(0,1)$  innovations and the sample size is~$500$.
    All curves are plotted against $\omega / (2\uppi)$. Real parts
    (Imaginary parts) of the periodogram and spectrum are presented in subfigures
    with $\tau_2 \leq \tau_1$ ($\tau_2 > \tau_1$): the dashed line represents the actual scaled
    spectrum [cf.~\protect\eqref{sd2}], the solid line the (pointwise) mean of the simulated smoothed Laplace periodograms.
     The gray areas represent the $0.1$, $0.25$, $0.75$ and $0.9$ (pointwise) sample quantiles of the smoothed periodograms over the
     5000 replications.}\label{fig:LPG-B4-M1}
\end{figure}

\begin{figure}

\includegraphics{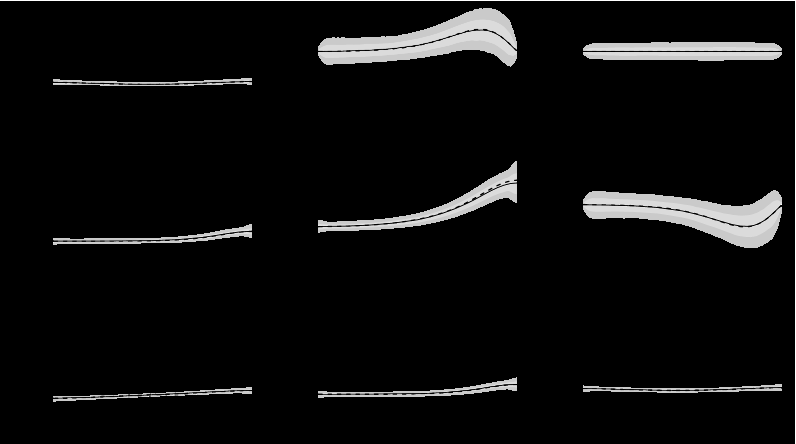}

    \caption{Smoothed rank-based  Laplace periodograms and (unscaled)
    spectral densities as defined in~\protect\eqref{sd2}.
    The process is an  AR(1) process with  $t_1$  innovations and  the sample size is~$500$.
    All curves are plotted against $\omega / (2\uppi)$. Real parts (Imaginary parts)
    of the periodogram and spectrum are presented in subfigures with $\tau_2 \leq \tau_1$ ($\tau_2 > \tau_1$):
    the dashed line represents the actual scaled spectrum [cf.~\protect\eqref{sd2}], the solid line the (pointwise) mean
    of the simulated smoothed Laplace periodograms. The gray areas represent the $0.1$, $0.25$, $0.75$ and $0.9$
    (pointwise) sample quantiles of the smoothed periodograms over the 5000 replications.}\label{fig:LPG-B4-M2}\vspace*{6pt}
\end{figure}

\begin{figure}

\includegraphics{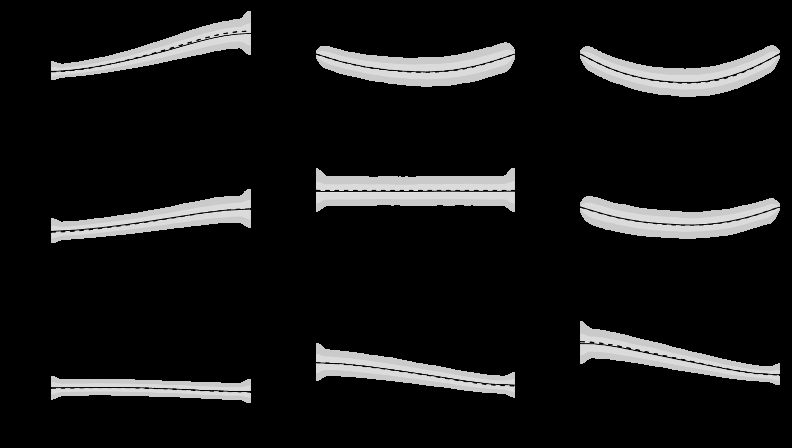}

    \caption{Smoothed rank-based  Laplace periodograms and (unscaled)
    spectral densities as defined in~\protect\eqref{sd2}.
    The process is a  QAR(1) process with  $\theta_1(u) = 1.9(u-0.5)$, $\theta_0(u) = 0.1 \Phi^{-1}(u)$ and the sample size is~$500$.
    All curves are plotted against $\omega / (2\uppi)$. Real parts (Imaginary parts) of the periodogram and
    spectrum are presented in subfigures with $\tau_2 \leq \tau_1$ ($\tau_2 > \tau_1$): the dashed line represents
    the actual scaled spectrum [cf.~\protect\eqref{sd2}], the solid line the (pointwise) mean of the simulated smoothed Laplace
    periodograms. The gray areas represent the $0.1$, $0.25$, $0.75$ and $0.9$ (pointwise) sample quantiles of the smoothed
    periodograms over the 5000 replications.}\label{fig:LPG-B4-M3}
\end{figure}

\begin{figure}

\includegraphics{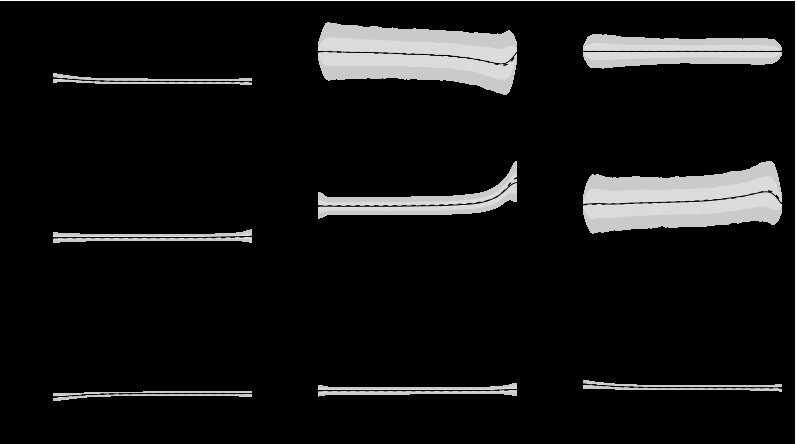}

    \caption{Smoothed rank-based  Laplace periodograms and (unscaled)
    spectral densities as defined in~\protect\eqref{sd2}.
    The process is an $\operatorname{ARMA}(1,1)$ process with  $t_3$  innovations and  the sample size is $500$.
    All curves are plotted against $\omega / (2\uppi)$. Real parts (Imaginary parts) of the periodogram and
    spectrum are presented in subfigures with $\tau_2 \leq \tau_1$ ($\tau_2 > \tau_1$): the dashed line
    represents the actual scaled spectrum [cf.~\protect\eqref{sd2}], the solid line the (pointwise) mean of the
    simulated smoothed Laplace periodograms. The gray areas represent the $0.1$, $0.25$, $0.75$ and $0.9$
    (pointwise) sample quantiles of the smoothed periodograms over the 5000 replications.}\label{fig:LPG-B4-M4}
\end{figure}

From all calculated periodograms, we determine the mean as an
approximation to the expectation of the various
estimates. Each of the following figures subdivides into nine subfigures.
For any combination of $\tau_1$ and $\tau_2$, the imaginary parts of
periodograms and spectra are represented above the diagonal, %whereas
and the real parts
%are represented
below; for $\tau_1 = \tau_2$, those quantities
are real and we represent them on the diagonal.
All curves are plotted against $\omega/ (2\uppi)$.
In all figures, the dashed line represents the ``true'' spectrum
(scaled for Figures~\ref{fig:LPG-B3-M1}--\ref{fig:LPG-B3-M4};
unscaled for Figures~\ref{fig:LPG-B4-M1}--\ref{fig:LPG-B4-M4}) and the
solid line the (pointwise) mean of the simulated smoothed Laplace periodograms.
%In order to illustrate the variability the figures contain some
%additional information.
The gray areas represent the $0.1$, $0.25$, $0.75$ and $0.9$
(pointwise) sample quantiles of the
smoothed periodograms from the $5000$ simulation runs.

%f1 #&#

%f2 #&#

%f3 #&#
%f4 #&#

%f5 #&#
%f6 #&#

%f7 #&#
%f8 #&#

For the sake of brevity, only results for sample size $n=500$
are presented here, but further results, which show a similar behavior,
are available from the authors.

We first discuss the results for the smoothed Laplace periodogram
in the case of an AR(1) process. Figure~\ref{fig:LPG-B3-M1} is with
Gaussian innovations, while the case of $t_1$-distributed innovations
is shown in Figure~\ref{fig:LPG-B3-M2}.
Inspection of these figures reveals that the imaginary component of the
spectrum is vanishing in the
case of Gaussian innovations (see Figure~\ref{fig:LPG-B3-M1}).
This observation reflects the fact that AR processes with Gaussian
innovations are time-reversible. On the other hand,
for $t_1$-distributed innovations, this
phenomenon only takes place for the extreme quantiles ($\tau_1 =
0.05$, $\tau_2 = 0.95$), meaning that $\IP(X_t \leq q_{0.05}, X_{t+k}
\leq q_{0.95})$ is approximately equal to $\IP(X_t \leq q_{0.95},
X_{t+k} \leq q_{0.05})$. This, however, does not hold for $\tau_1 =
0.5$ and $\tau_2 = 0.05$ or $0.95$, which indicates a
time-irreversible impact of extreme values on the central ones.

In Figure~\ref{fig:LPG-B3-M3}, the simulation results for the QAR(1)
process are shown. We see that the (scaled) copula spectrum for $\tau
_1 = \tau_2 = 0.25$ has the shape previously observed in the case of
the AR(1) process, where the autoregressive parameter was negative.
Note that the function $\theta_1(u)$ takes negative values for $u \in
(0,0.5)$. On the other hand, for $\tau_1 = \tau_2 = 0.75$, it has the
shape of the spectral density in the AR(1) case when the autoregressive
parameter is positive, while $\theta_1(u)$ is positive for $u \in
(0.5,1)$. For $\tau_1 = \tau_2 = 0.5$ we observe a flat spectrum,
indicating that the sequence $( I\{Y_t \leq q_{0.5}\} )$ has zero
autocorrelation, which would imply $\IP(X_t \leq q_{0.5}, X_{t+k} \leq
q_{0.5}) = \IP(X_t \leq q_{0.5}) \IP(X_{t+k} \leq q_{0.5})$. The
imaginary part of the spectrum clearly indicates time-irreversibility,
which implies that the QAR(1) process, irrespective of the choice of
$\theta_0$, cannot be a Gaussian process.

The simulation results for the all-pass $\operatorname{ARMA}(1,1)$ process are shown in
Figure~\ref{fig:LPG-B3-M3}. We see here that the statistics proposed
are very able of capturing the serial dependence which (due to
uncorrelatedness) would completely escape the traditional analysis.
Another finding is that, in most cases, the bias is larger for the
estimation of the Laplace spectrum with $\tau_1=\tau_2$: see, for
instance, the diagonals of
Figures~\ref{fig:LPG-B3-M1}--\ref{fig:LPG-B3-M4}.

The corresponding rank-based Laplace periodograms are shown in
Figures~\ref{fig:LPG-B4-M1}--\ref{fig:LPG-B4-M4}, respectively.
The results indicate the same type of time-reversibility features as
observed with the Laplace periodogram.
It is interesting to note that, for the rank-based Laplace periodograms,
the bias appears to be much smaller, and smoothing seems to be
more effective.

Finally, we investigate the quality of the estimates by their mean
squared properties.
In Table~\ref{tab:rmse}, we provide the square roots of the integrated
mean squared errors (MSE). We consider the smoothed rank-based Laplace
periodograms for
sample sizes $n=100$, $500$, and $1000$. Note that, because of
symmetry, we do
not display all combinations. For example, the spectra corresponding to
the quantiles $(0.05,0.05)$
and $(0.95,0.95)$ coincide in the scenario under consideration. In all
cases, we observe, from the point of view of MSE, a reasonable behavior
of the
rank-based Laplace periodograms. It is interesting to
note that the integrated MSE is larger when quantiles in the
neighborhood of
$\tau=0.5$ are involved. For example, the integrated MSE is increasing
from $(0.05,0.05)$ to $(0.05,0.25)$
and $(0.05,0.50)$, then decreasing from $(0.05,0.75)$ to $(0.05,0.95)$.
This phenomenon is closely related to the fact that the empirical
copula has variance zero at the boundaries of the unit cube, see Genest and Segers \cite
{GenestSegers2010} for more details on this fact.

%s5.2 #&#
\subsection{An empirical application: S\&P 500 returns}

The smoothed rank-based Laplace periodogram was computed from the
series of daily return values of the S\&P 500 index
(Jan/2/1963--Dec/31/2009, $n=11\,844$), based on a Daniell kernel with
parameters $(200,100)$, for the same quantile orders as in the previous section.
The results for
the smoothed traditional periodogram are shown in Figure~\ref{fig:pgram-SP500},
and those for the rank-based Laplace periodogram
in Figure~\ref{fig:LPGrb-SP500}, with the same convention as in
Section~\ref{secsim}.
%The autocorrelations of order up to 40 for the S\&P 500 only range
%from $-0.05$ to $0.04$.

The nonlinear features of that series have been stressed by many
authors (see, e.g., Abhyankar, Copeland and
  Wong \cite
{AbhyankarEtAl1997}, Berg, Paparoditis and
  Politis \cite{BergEtAl2010}, Brock, Hsieh and
  LeBaron \cite{BrockEtAl1992}, Hinich and
  Patterson \cite{HinichPatterson1990,HinichPatterson1985}, Hsieh \cite{Hsieh1989},
and Vaidyanathan and
  Krehbiel \cite{VaidyanathanKrehbiel1992}). Those nonlinear features
cannot be detected by classical correlogram-based spectral methods, and
hence do not appear in Figure~\ref{fig:pgram-SP500}, where the
traditional smoothed periodogram is depicted. They do appear, however,
in the plots of Figure~\ref{fig:LPGrb-SP500}. Whereas the picture for
the central quantiles $\tau_1 = \tau_2 = 0.5$ looks quite similar to
that in Figure~\ref{fig:pgram-SP500}, the remaining ones, which
involve at least one extreme quantile, are drastically different,
indicating a marked discrepancy between tail and central dependence
structures. All plots involving at least one extremal quantile yield a
peak at the origin, which possibly corresponds to a long-range memory
for extremal events. Imaginary parts are not zero, suggesting again
time-irreversibility. Such features entirely escape a traditional
spectral analysis.
%By design, the Laplace spectral density kernel may serve as a tool to
%represent the periodicities in nonlinear dependencies, while the
%autocorrelation function and ordinary spectral density fail to do so,
%since they reflect only linear dependencies.
%
%The smoothed Laplace periodogram for $\tau_1 = \tau_2 = 0.5$ looks
%similar to the smoothed ordinary periodogram, but the remaining
%Laplace periodograms differ largely from the one in the middle. All
%the plots involving at least one extremal quantile have in common that
%there seems to be a peak around $\omega= 0$. Possibly, this could be
%interpreted as an indication for long-range dependencies. This feature
%entirely escapes the traditional spectral analysis.
%Assumption (A8) follows.\hfill$\Box$

%%%%%%%%%%%%%%%%%%%%%%%%%%%%%%%%%%%%%%%%%%%%%%%%%%%%%%%%%%%%%%%%%%%%%%%%%%%%%%%%%%%%%%%%%%%%%%%%%%%%%%%%%%%%%%%%%%%%%%%%%%%%%%%%%%%%%%%%%%%%%%%%%%%%%%%%%%%%%%%%%%%%%%%%%%%%%%%%%%%%%%%%%%%%%%%%%%%%%%%%%%%%%%%%%%%%%%%%%%%%%%%%%%%%%%%%%%%%%%%%%%

%sA #&#
\begin{appendix}
%sA #&#
\section{Technical details for the proofs in
Section \texorpdfstring{\protect\ref{321315458}}{3}}\label{appA}

In this section, we give the technical details for the proofs of
Theorems \ref{PropBAsymp} and \ref{PropCLAsymp}. Those proofs rely on
a series of lemmas. Two of them (Lemmas  \ref{lem:bennett} and \ref
{lem:df}) are general results, to be used at several places in both
proofs; their statements and proofs are postponed to Section~\ref{BasicSec}. Two further ones (Lemmas \ref{lem:df2} and \ref
{lem:nummn}) are specific to the proof of \eqref{thm:uniflinearrank}
and Theorem~\ref{PropCLAsymp}: they are presented in Section~\ref{subsub211}.
Finally, Lemmas \ref{lem:convrate} and \ref{lem:boundelta} are
auxiliary results used in the proofs of \eqref{thm:uniflinearworank}
and \eqref{thm:uniflinearrank}; they are regrouped in Section~\ref{subsub121}, along with Lemma~\ref{lem:varzn2}, %(also in Section
which is specific to the proof of \eqref{thm:uniflinearworank}. %

%sA.1 #&#
\subsection{Details for the proof of \texorpdfstring{\protect\eqref{thm:uniflinearworank}}{(3.9)}} \label{details}

Recall that \eqref{thm:uniflinearworank} was obtained by combining
Lemmas \ref{lem:convrate} and \ref{lem:boundelta} with
Equation \eqref{conslemm2}. In Section~\ref{subsub12x}, we
establish \eqref{conslemm2}, thus completing (but for Lemmas \ref{lem:convrate}--\ref{lem:varzn2})
the proof of Theorem~\ref{PropBAsymp}.
In Section~\ref{subsub121}, we state and prove Lemmas \ref
{lem:convrate}--\ref{lem:varzn2}, which completes the proof of (\ref
{thm:uniflinearworank}). The notation of Theorem~\ref{PropBAsymp} is
used throughout this section.

%sA.1.1 #&#
\subsubsection{Proof of \texorpdfstring{\protect\eqref{conslemm2}}{(3.11)}} \label{subsub12x}

In this proof, we use a blocking argument by Yu \cite{Yu1994} -- call
it the \textit{independent blocks argument}. Let $m_n := \lceil
n^{1/(1+\delta)} \log n \rceil$, $\mu_n := \lfloor n/(2 m_n) \rfloor
$, and split the set $\{1,\ldots,n\}$ into $2 \mu_n$ subsets of size
$m_n$ and a ``residual'' subset of size $n - 2 m_n \mu_n$:
%
%eA.1 #&#
\begin{eqnarray}
\label{eqn:blocks} %\begin{split}
S_i & := &\bigl\{ k \in\IN\dvt 2(i-1)
m_n + 1 \leq k \leq(2i-1) m_n\bigr\},\quad\quad i = 1,\ldots,
\mu_n,\nonumber
\\
T_i & := &\bigl\{ k \in\IN\dvt (2i-1) m_n + 1 \leq k
\leq2i m_n\bigr\},\quad\quad i = 1,\ldots, \mu_n,
\\
R_n & := &\{2 m_n \mu_n + 1,\ldots,n\}.\nonumber
\end{eqnarray}
%
%As in \citeauthor{Yu1994} \cite{Yu1994} consider an independent block $m_n$-sequence
%$X_1, \ldots, X_n$. This sequence can be obtained, by first splitting
%$Y_1,\ldots,Y_n$ into the $2 \mu_n + 1$ blocks defined in
%(Y_t)_{t \in S_{\mu_n}}, (Y_t)_{t \in T_{\mu_n}} \text{ and } (Y_t)_{t
%The independent block $m_n$-sequence then is defined by choosing
%$(X_t)_{t \in S_j}$ (respective $(X_t)_{t \in T_j}$ or $(X_t)_{t \in
%R_n}$) to have the same distribution as $(Y_t)_{t \in S_j}$
%(respective $(Y_t)_{t \in T_j}$ or $(Y_t)_{t \in R_n}$) and such that
%$X_1, \ldots, X_n$ is independent of $Y_1, \ldots, Y_n$ and such that
%the $2 \mu_n + 1$ blocks of $X_1, \ldots, X_n$ are independent.
Associated with this partition of $\{1,\ldots,n\}$, consider the partition
\[
(Y_t)_{t\in S_1}, (Y_t)_{t\in T_1} ;
(Y_t)_{t\in S_2},\ldots, (Y_t)_{t\in T_{\mu_n-1}};
(Y_t)_{t\in S_{\mu_n}}, (Y_t)_{t\in T_{\mu
_n}};
(Y_t)_{t\in R_n}
\]
of $\{Y_1,\ldots, Y_n\}$ into $2\mu_n $ blocks of length $m_n$ and a
``residual'' block of length $n-2m_n\mu_n$. The independent block
$m_n$-sequence then is defined as a collection of $2\mu_n$ mutually
independent $m_n$-dimensional random variables $(X_t)_{t\in S_i}$,
$(X_t)_{t\in T_i}$, $i=1,\ldots, \mu_n$, such that $(X_t)_{t\in S_i}
\stackrel{d}{=}(Y_t)_{t\in S_i}$ and $(X_t)_{t\in T_i} \stackrel
{d}{=}(Y_t)_{t\in T_i}$, along with an $(n-2m_n\mu_n)$-dimensional
variable $(X_t)_{t\in R_n}$ independent of all other $(X_t)$'s.\def\thetable{1}

The independent blocks argument will be used to establish results of
the form
\[
{\rm P} \Biggl( \sup_{\theta\in\Theta_n} \Biggl| \sum
_{t=1}^n \theta (t,Y_t) \Biggr| >
\eta_n \Biggr) = \mathrm{o}(1),
\]
where $\Theta_n$ are sets of measurable functions $\theta\dvtx  \IR^{2}
\rightarrow\IR$.%
For the argument, consider the decomposition\vspace*{-1pt}
\begin{eqnarray*}
&&{\rm P} \Biggl( \sup_{\theta\in\Theta_n} \sum_{t=1}^n
\theta(t,Y_t) > \eta_n \Biggr)
\\
&&\quad\leq{\rm P} \Biggl( \sup
_{\theta\in\Theta_n} \Biggl| \sum_{i=1}^{\mu
_n}
\sum_{t \in S_i} \theta(t,Y_t) \Biggr| >
\eta_n/3 \Biggr)
\\
&&\quad\quad{}+ {\rm P} \Biggl( \sup_{\theta\in\Theta_n} \Biggl| \sum
_{i=1}^{\mu
_n} \sum_{t \in T_i}
\theta(t,Y_t) \Biggr| > \eta_n/3 \Biggr) + {\rm P} \biggl( \sup
_{\theta\in\Theta_n} \Biggl| \sum_{t \in R_n}
\theta(t,Y_t) \Biggr| > \eta_n/3 \biggr)
\\
&&\quad=: P^{(1)}_n + P^{(2)}_n +
P^{(3)}_n.
\end{eqnarray*}

%
%Because of $\sup_{\theta\in\Theta_n} \Big| \sum_{t \in R_n}
%that the last of the three probabilities on the right hand side of the
%inequality is $o(1)$.\nepage
\noindent The last probability $P^{(3)}_n$ is zero as soon as
\[
 \mbox{(i)}  \quad  \sup_{\theta\in\Theta_n} \sup_{t=1,\ldots,n} \bigl|
\theta(t,Y_t) \bigr| \leq C_n \qquad\mbox{a.s.}\quad \mbox{and}\quad
m_n C_n < \eta_n/3,
\]
which will be the case in all applications of the independent blocks
argument. The first probability $P^{(1)}_n$ can be handled by applying
Lemma~4.1 from Yu \cite{Yu1994}, by which we have
\[
{\rm P} \Biggl( \sup_{\theta\in\Theta_n} \Biggl| \sum
_{i=1}^{\mu_n} \sum_{t \in S_i}
\theta(t,Y_t) \Biggr| > \eta_n/3 \Biggr) \leq {\rm P} \Biggl(
\sup_{\theta\in\Theta_n} \Biggl| \sum_{i=1}^{\mu_n}
\sum_{t \in S_i} \theta(t,X_t) \Biggr| >
\eta_n/3 \Biggr) + \mathrm{o}(1),
\]
since by the choice of $m_n$ we have $\mu_n \beta(m_n) = \mathrm{o}(1)$. A
similar argument applies to the second probability $P^{(2)}_n$. We
assume that the set $\Theta_n$ consists of finitely many, say $|\Theta
_n|$, elements to further obtain
\[
{\rm P} \Biggl( \sup_{\theta\in\Theta_n} \Biggl| \sum
_{i=1}^{\mu_n} \sum_{t \in S_i}
\theta(t,X_t) \Biggr| > \eta_n/3 \Biggr) \leq|
\Theta_n| \sup_{\theta\in\Theta_n} {\rm P} \Biggl( \Biggl| \sum
_{i=1}^{\mu_n} \sum
_{t \in S_i} \theta(t,X_t) \Biggr| > \eta_n/3
\Biggr),
\]
where the summands $\sum_{t \in S_i} \theta_t(X_t)$, $i=1,\ldots,\mu
_n$ are independent by construction. If we additionally show that
\[
\mbox{(ii)} \quad  \sup_{\theta\in\Theta_n} \sum_{j=1}^{\mu_n} \Var \Biggl(\sum_{t\in S_j} \theta(t,X_t)  \Biggr)
 \leq
V_n^2
\quad\mbox{and}\quad
   \sup_{\theta\in\Theta_n} \sum_{j=1}^{\mu_n} \Var
\Biggl(\sum_{t\in T_j} \theta(t,X_t)  \Biggr) \leq V_n^2,
\]
%

%t1 #&#
\begin{sidewaystable}
\tablewidth=\textwidth
\tabcolsep=0pt
\caption{Root Integrated Mean Square
Errors of smoothed, rank-based Laplace periodograms, for the four
models described in Section~\protect\ref{secsim}, and various series lengths}
\label{tab:rmse}
\begin{tabular*}{\textwidth}{@{\extracolsep{4in minus 4in}}llllllllll@{}}
\hline
& & \multicolumn{8}{l@{}}{$(\tau_1, \tau_2)$} \\[-5pt]
&&\multicolumn{8}{l@{}}{\hrulefill}\\
$Y_t$ & $n$ & $(0.05,0.05)$ & $(0.05,0.25)$ & $(0.05,0.5)$ & $(0.05,0.75)$ &
$(0.05,0.95)$ & $(0.25,0.25)$ & $(0.25,0.5)$ & $(0.5,0.5)$\\
\hline
Model (1) & \hphantom{0}$100$ & 0.0212 & 0.0408 & 0.0459 & 0.0401 & 0.0219 & 0.0651 & 0.0837
& 0.0876 \\
& \hphantom{0}$500$ & 0.0085 & 0.0185 & 0.0215 & 0.0189 & 0.0099 & 0.0347
& 0.0429 & 0.0474 \\
& $1000$ & 0.0054 & 0.0117 & 0.0137 & 0.0121 & 0.0064 & 0.0225 & 0.0275
& 0.0310 \\[3pt]
Model (2)& \hphantom{0}$100$ & 0.0223 & 0.0418 & 0.0462 & 0.0405 & 0.0234 & 0.0672 & 0.0852
& 0.0929 \\
 & \hphantom{0}$500$ & 0.0091 & 0.0188 & 0.0213 & 0.0188 & 0.0110 & 0.0353
& 0.0441 & 0.0506 \\
& $1000$ & 0.0059 & 0.0120 & 0.0135 & 0.0120 & 0.0072 & 0.0228 & 0.0282
& 0.0330 \\[3pt]
Model (3)& \hphantom{0}$100$ & 0.0207 & 0.0398 & 0.0452 & 0.0386 & 0.0214 & 0.0652 & 0.0830
& 0.0873 \\
 & \hphantom{0}$500$ & 0.0084 & 0.0184 & 0.0213 & 0.0186 & 0.0098 & 0.0349
& 0.0428 & 0.0471 \\
& $1000$ & 0.0053 & 0.0115 & 0.0135 & 0.0119 & 0.0064 & 0.0227 & 0.0277
& 0.0309 \\[3pt]
Model (4)& \hphantom{0}$100$ & 0.0220 & 0.0412 & 0.0453 & 0.0398 & 0.0226 & 0.0654 & 0.0834
& 0.0873 \\
 & \hphantom{0}$500$ & 0.0097 & 0.0191 & 0.0214 & 0.0190 & 0.0108 & 0.0344
& 0.0422 & 0.0465 \\
& $1000$ & 0.0064 & 0.0122 & 0.0135 & 0.0121 & 0.0071 & 0.0226 & 0.0271
& 0.0306 \\
\hline
\end{tabular*}
\end{sidewaystable}

\def\thefigure{9}
%f9 #&#
\begin{figure}\vspace*{20pt}%

\includegraphics{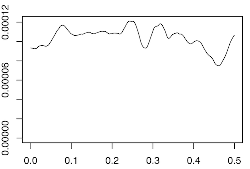}

    \caption{Smoothed traditional periodogram, S\&P~500 returns
    curve is plotted against $\omega / (2\uppi)$.}\label{fig:pgram-SP500}
%\end{figure}
\vspace*{50pt}

\def\thefigure{10}%

%f10 #&#
%\begin{figure}[p]

\includegraphics{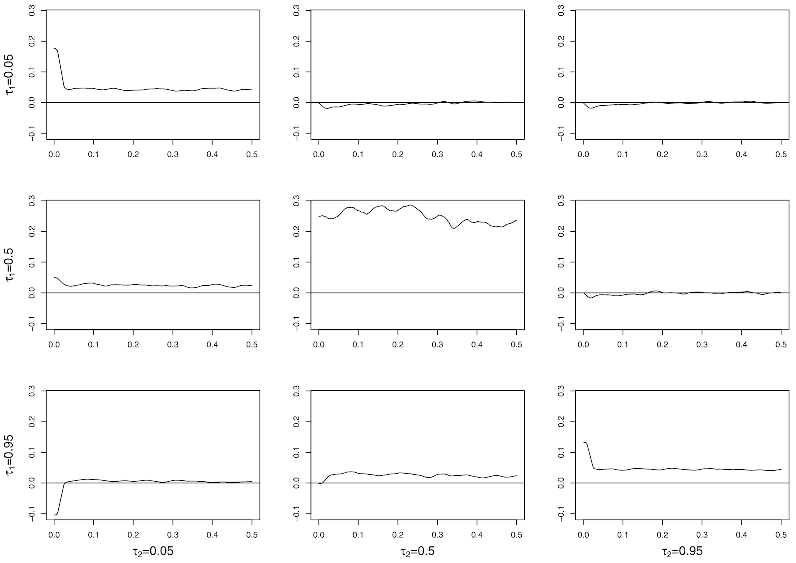}

    \caption{Smoothed rank-based Laplace periodograms, S\&P~500 returns.
    All curves are plotted against $\omega / (2\uppi)$. Real parts (Imaginary parts)
    of the periodogram and spectrum are presented in subfigures with $\tau_2 \leq \tau_1$ ($\tau_2 > \tau_1$).}   \label{fig:LPGrb-SP500}
\end{figure}%

\noindent the version of Bennett's inequality given in Lemma~\ref{lem:bennett}
can be applied, so that, under (i) and (ii),
\begin{eqnarray*}
{\rm P} \Biggl( \sup_{\theta\in\Theta_n} \Biggl| \sum
_{i=1}^{\mu_n} \sum_{t \in S_i}
\theta(t,X_t) \Biggr| > \eta_n/3 \Biggr) &\leq&{\rm P} \Biggl(
\sup_{\theta\in\Theta_n} \Biggl| \sum_{i=1}^{\mu_n}
\sum_{t \in S_i} \bigl(\theta(t,X_t) - {\rm
E}\bigl[\theta(t,X_t)\bigr] \bigr) \Biggr| > \lambda_n \Biggr)
\\
&\leq&2|\Theta_n|\exp \biggl(-\frac{\log2}{4} \biggl(
\frac{\lambda
_n^2}{2 V_n^2} \wedge\frac{\lambda_n}{m_nC_n} \biggr) \biggr),
\end{eqnarray*}
where $\lambda_n := \eta_n/3 - n \sup_{\theta\in\Theta_n} \sup_{t=1,\ldots,n} | {\rm E}[\theta(t,X_t)] |$.
Exactly the same argument can be used to handle the probability
$P^{(2)}_n$. Hence, we obtain
%
%eA.2 #&#
\begin{eqnarray}
\label{def:en} &\displaystyle{\rm P} \Biggl( \sup_{\theta\in\Theta_n} \Biggl| \sum
_{t=1}^n \theta (t,Y_t) \Biggr| >
\eta_n \Biggr) \leq E_n + \mathrm{o}(1),\nonumber&
\\
[-8pt]\\[-8pt]
&\displaystyle E_n := 4|
\Theta_n|\exp \biggl(-\frac{\log2}{4} \biggl(\frac{\lambda
_n^2}{2 V_n^2}
\wedge\frac{\lambda_n}{m_nC_n} \biggr) \biggr).\nonumber&
\end{eqnarray}
An application of the independent block argument for finite $\Theta_n$
thus boils down to establishing (i)--(ii) discussed above and ensuring
that $E_n = \mathrm{o}(1)$.

Regarding the proof of \eqref{conslemm2} note that, it is obviously
possible to construct $N = \mathrm{o}(n^5)$ points $d_{1},\ldots,d_{N}$ (dependence
on $n$ is not reflected in
the notation) such that, for every $\bdelta$ with $\|\bdelta\| \leq
A\sqrt{\log n}$, there exists an index $j(\delta)$ for which $\|
\bdelta-d_{j(\delta)}\| \leq n^{-3/2}$.
Define
\[
K_{n}(\bdelta;\tau,\omega) := \sum_{t=1}^{n}
\biggl( \int_0^{n^{-1/2} \mathbf{c}^\tr_{t}(\omega) \bdelta} \bigl(I\{Y_t
\leq s + q_{\tau}\} - I\{Y_t \leq q_{\tau}\}\bigr)
\,\mathrm{d}s - {f(q_{\tau
})}(2n)^{-1}\bigl(\mathbf{c}_{t}^\tr(
\omega) \bdelta\bigr)^2 \biggr)
\]
and note, by direct calculation, that, for $n\geq n_0$ with $n_0\in\IN
$ independent of $\tau$ and $\omega$,
\[
\sup_{\omega\in\mathcal{F}_n } \bigl| K_{n}(a;\tau,\omega) -
K_{n}(b;\tau,\omega) \bigr| \leq1.5 \sqrt{n}\|a-b\|.
\]
By applying Knight's identity, we therefore have
\[
\sup_{\omega\in\mathcal{F}_n} \sup_{\|\bdelta\| \leq A\sqrt{\log n}} \bigl| \Zyh{n,\tau,
\omega}(\bdelta) - \Zy{n, \tau, \omega}(\bdelta) \bigr| = \sup_{\theta\in\Theta_n}
\Biggl| \sum_{t=1}^n \theta(t,
Y_t) \Biggr| + \mathrm{O}_{\rm P}\bigl(n^{-1}\bigr),
\]
where
\begin{eqnarray*}
\Theta_n &:=& \biggl\{ \theta(t,y) := \int_0^{n^{-1/2} \mathbf{c}^\tr_{t}(\omega) d_j}
\bigl(I\{y \leq s + q_{\tau}\} - I\{y \leq q_{\tau}\}\bigr)
\,\mathrm{d}s - {f(q_{\tau})}(2n)^{-1}\bigl(\mathbf{c}_{t}^\tr(
\omega) d_j\bigr)^2 \Big|
\\
&&\hphantom{\biggl\{}\omega\in\mathcal{F}_n, j=1,\ldots,N \biggr\}.
\end{eqnarray*}

In order to show that $ \sup_{\theta\in\Theta_n}  | \sum_{t=1}^n
\theta(t, Y_t)  | = \mathrm{O}_{\rm P}(r_n(\delta)^2)$, we apply
the independent blocks argument with $\Theta_n$ defined above and
$\eta_n := D r_n(\delta)^2$ for a suitable constant $D$.

Due to the fact that $n^{(1-\delta)/(2+2\delta)} (\log n)^{3/2} \ll
r_n(\delta)^2$ and that, by Lemma~\ref{lem:varzn2},
\[
\sup_{\theta\in\Theta_n} \sup_{t=1,\ldots,n} \bigl| \theta(t,
Y_t) \bigr| \leq C n^{-1/2} (\log n)^{1/2} =:
C_n,
\]
almost surely, (i) in the independent blocks argument follows.

Next, a direct calculation shows that (ii) in the independent blocks
argument holds with $V_n^2 := C n^{-1/2} (\log n)^2$.

Finally, let us complete the independent blocks argument by
establishing that for $E_n$ defined in \eqref{def:en} we have $E_n =
\mathrm{o}(1)$. Observe that the bounds in Lemma~\ref{lem:varzn2} imply
\[
\sup_{\theta\in\Theta_n} \sup_{t=1,\ldots,n} {\rm E}\bigl[\bigl|
\theta(t, X_t)\bigr|\bigr] \leq C \log(n)^3 n^{-3/2} =
\mathrm{o}\bigl(n^{-1}r_n(\delta)^2\bigr).
\]
Thus, we find that for sufficiently large $n$
\[
\lambda_n := D \Bigl(r_n(\delta)^2/3 - n
\sup_{\theta\in\Theta_n} \sup_{t=1,\ldots,n} {\rm E}\bigl[\bigl|
\theta(t, X_t)\bigr|\bigr] \Bigr) \leq D r_n(
\delta)^2/6.
\]
Noting that $|\Theta_n| = N n = \mathrm{o}(n^6)$ direct calculations yield $E_n
= \mathrm{o}(1)$ for $D$ in the definition of $\eta_n$ being large enough. This
completes the application of the independent blocks argument and shows
that $ \sup_{\theta\in\Theta_n}  | \sum_{t=1}^n \theta(t,
Y_t)  | = \mathrm{O}_{\rm P}(r_n(\delta)^2)$.

Summing up, except for Lemma~\ref{lem:varzn2} which is taken care of in the next
section, we have proven \eqref{conslemm2}. If we now prove Lemmas \ref{lem:convrate}
and \ref{lem:boundelta}, \eqref{conslemm1} and \eqref{thm:uniflinearworank}, hence
Theorem~\ref{PropBAsymp}, follow.
The purpose of Section~\ref{subsub121} below is to complete the proof of
Theorem~\ref{PropBAsymp} by establishing the missing Lemmas \ref{lem:convrate}--\ref{lem:varzn2}. %\qed

%sA.1.2 #&#
\subsubsection{Three auxiliary lemmas} \label{subsub121}
We now state and prove the three lemmas that have been used in the
proof of Theorem~\ref{PropBAsymp}.
Lemma~\ref{lem:convrate} generalizes ideas from Pollard \cite{Pollard1991}.

%leA.1 #&#
\begin{lemma}\label{lem:convrate}
Let $B_{a_n}(\mathbf{x})$ denote the closed ball (in $\mathbb{R}^3$)
with center $\mathbf{x}$ and radius $a_n > 0$.
Assume that, for some sequence of real numbers $a_n=\mathrm{o}(1)$,
\[
\Delta_n := \sup_{\omega\in\mathcal{F}_n} \sup_{\bdelta\in
B_{a_n}(\dey{n,\tau,\omega})}
\bigl| \Zyh{n,\tau, \omega}(\bdelta) - \Zy{n, \tau, \omega}(\bdelta )\bigr| =
\mathrm{o}_{\mathrm{P}}\bigl(a_n^2\bigr).
\]
Then, $
\sup_{\omega\in\mathcal{F}_n} | \deyh{n,\tau,\omega} - \dey
{n,\tau,\omega} | = \mathrm{o}_{\mathrm{P}}(a_n)
$.
\end{lemma}

\begin{pf}
Let $r_{n,\tau,\omega}(\bdelta) := \Zyh{n,\tau,\omega}(\bdelta)
- \Zy{n, \tau, \omega}(\bdelta)$. Simple algebra and the explicit
form (\ref{deltaquer}) of $\dey{n, \tau, \omega}$
%(recall that $\dex{n, \tau, \omega} = (\bQx_{n,\tau,\omega})^{-1}
yield
%
%eA.3 #&#
\begin{equation}
\label{h1} \Zyh{n,\tau, \omega}(\bdelta) = \tfrac{1}{2}(\bdelta-
\dey{n,\tau,\omega})^\tr\bQy_{n,\tau
,\omega} (\bdelta- \dey{n,\tau,
\omega}) - \tfrac{1}{2} (\dey{n,\tau,\omega})^\tr
\bQy_{n,\tau,\omega} \dey{n,\tau,\omega} + r_{n,\tau,\omega}(\bdelta).
\end{equation}
Any $\bdelta\in\IR^3 \setminus B_{a_n}(\dey{n,\tau,\omega})$
with distance $l_n := \|\bdelta-\dey{n,\tau,\omega}\| > a_n$ to
$\dey{n,\tau,\omega}$ can be represented as $
\bdelta= \dey{n,\tau,\omega} + l_{n,\tau,\omega} \mathbf
{d}_{n,\tau,\omega}$,
where
$
\mathbf{d}_{n,\tau,\omega} := l_{n,\tau,\omega}^{-1}(\bdelta-\dey
{n,\tau,\omega})$.

The point $\bdelta_{n,\tau,\omega}^* = \dey{n,\tau,\omega} + a_n
\mathbf{d}_{n,\tau,\omega}$ lies on the boundary of the ball
$B_{a_n}(\dey{n,\tau,\omega})$. The convexity of $\Zyh{n,\tau,
\omega}(\bdelta)$ therefore implies
\begin{eqnarray*}
&&{a_n} {l_{n,\tau,\omega}^{-1}} \Zyh{n,\tau, \omega}(
\bdelta) + \bigl(1-{a_n} {l_{n,\tau,\omega}^{-1}} \bigr)
\Zyh{n,\tau, \omega }(\dey{n,\tau,\omega})
\\
&&\quad\geq \Zyh{n,\tau, \omega}\bigl(\bdelta_{n,\tau,\omega}^*\bigr) = \Zy{n, \tau,
\omega}\bigl(\bdelta_{n,\tau,\omega}^*\bigr) + r_{n,\tau
,\omega}\bigl(
\bdelta_{n,\tau,\omega}^*\bigr)
\\
&&\quad\geq \tfrac{1}{2}\mathbf{d}_{n,\tau,\omega}^\tr
\bQy_{n,\tau
,\omega} \mathbf{d}_{n,\tau,\omega} a_n^2 -
\tfrac{1}{2}(\dey {n,\tau,\omega})^\tr\bQy_{n,\tau,\omega} \dey{n,
\tau,\omega} - \Delta_n
\\
&&\quad\geq \tfrac{1}{2}\mathbf{d}_{n,\tau,\omega}^\tr
\bQy_{n,\tau
,\omega} \mathbf{d}_{n,\tau,\omega} a_n^2 +
\Zyh{n,\tau, \omega }(\dey{n,\tau,\omega}) - 2 \Delta_n.
\end{eqnarray*}
Rearranging and taking the infimum over $\omega$ and $\bolds\delta$,
we obtain
%
%eA.4 #&#
\begin{eqnarray}
\label{eqn:ConvergMinima}&& \inf_{\omega\in\mathcal{F}_n} \inf_{\bdelta: |\bdelta- \dex
{n,\tau,\omega}| > a_n}
\bigl( \Zyh{n,\tau, \omega}(\bdelta) - \Zyh{n,\tau, \omega}(\dey{n,\tau,\omega})
\bigr)\nonumber
\\[-8pt]\\[-8pt]
&&\quad\geq\inf_{\omega\in\mathcal{F}_n} \inf_{\bdelta: |\bdelta-
\dey{n,\tau,\omega}| > a_n}
{l_{n,\tau,\omega}} {a_n^{-1}} \biggl(\frac{1}{2}
\mathbf{d}_{n,\tau,\omega}^\tr\bQy_{n,\tau,\omega}
\mathbf{d}_{n,\tau,\omega} a_n^2- 2 \Delta_n
\biggr).\nonumber
\end{eqnarray}
By assumption,
the smallest eigenvalue of $\bQy_{n,\tau,\omega}$ is bounded away
from zero uniformly in $\omega\in\mathcal{F}_n$, for $n$
sufficiently large. Hence, $2 \Delta_n < \frac{1}{2}\mathbf
{d}_{n,\tau,\omega}^\tr\bQy_{n,\tau,\omega} \mathbf{d}_{n,\tau
,\omega} a_n^2$ with probability tending to one,
the right-hand side in (\ref{eqn:ConvergMinima}) is strictly positive,
and the minimum of $\Zyh{n,\tau, \omega}(\bdelta)$ cannot be
attained at any $\bdelta$ with $|\bdelta- \dey{n,\tau,\omega}| >
a_n$.
\end{pf}

%leA.2 #&#
\begin{lemma} \label{lem:boundelta}
Let \textup{\ref{(A1)}} hold, and $\dey{n,\tau,\omega}$ be defined as in \eqref{deltaquer}.
Then, for any $\tau\in(0,1)$ for which $f(q_{\tau}) > 0$, there
exists a constant $A$ such that
\[
\lim_{n\to\infty} \mathrm{P} \Bigl(\sup_{\omega\in\mathcal{F}_n} \|
\dey{n,\tau ,\omega}\| > A\sqrt{\log n} \Bigr) = 0. % \sup_{\omega\in\mathcal{F}_n} \|\tilde\bdelta_{n,\tau,\omega}\|_
\]
\end{lemma}

\begin{pf}
Denote by $\|\mathbf{x}\|_{\infty}$ the sup-norm of $\mathbf{x}$. Since, for
${\bf x}\in\mathbb{R}^3$, $\sqrt{3} \|{\bf x}\|_{\infty} \geq\|
{\bf x}\|$, it suffices to prove that
\[
\lim_{n\to\infty} \mathrm{P} \Bigl(\sup_{\omega\in\mathcal{F}_n} \|
\dey{n,\tau ,\omega}\|_{\infty} > 3^{-1/2} A\sqrt{\log n} \Bigr) =
0.
\]
Next note that $\sqrt{n} \sup_{\omega\in\mathcal{F}_n} \|\dey
{n,\tau,\omega}\|_{\infty} = \sup_{\theta\in\Theta_n} | \sum_{t=1}^{n} \theta(t,Y_t) |$,
where
\[
\Theta_n := \bigl\{ \theta(t,y) := f(q_{\tau})^{-1}
c_{t,k}(\omega) \bigl(\tau- I\{y \leq q_{\tau}\}\bigr) |
k=1,2,3, \omega\in\mathcal{F}_n\bigr\},
\]
with $(c_{t,1}(\omega), c_{t,2}(\omega), c_{t,3}(\omega)) := (1,
\cos(\omega t), \sin(\omega t))$.

We apply the independent blocks argument described in Section~\ref{subsub12x}, with $\Theta_n$ defined above and $\eta_n := 3^{-1/2} A
n^{1/2} (\log n)^{1/2}$ with a suitably chosen constant $A$. To this
end, remark that (i) in the independent blocks argument holds for $A$
large enough, because we have, almost everywhere,
\[
\sup_{\theta\in\Theta_n} \sup_{t=1,\ldots,\mu_n m_n} \bigl|\theta
(t,Y_t)\bigr| \leq\frac{2}{f(q_{\tau})} =: C_n
\]
which implies,
\[
\sup_{\theta\in\Theta_n} \Biggl| \sum_{t \in R_n}
\theta(t,Y_t) \Biggr| \leq\frac{2 m_n}{f(q_{\tau})}\quad\quad \mbox{a.e.}
\]
Regarding (ii) from the independent blocks argument note that for all
$\theta\in\Theta_n$
\begin{eqnarray*}
\Var \biggl( \sum_{t \in S_i}
\theta(t,X_t) \biggr) & =& \sum_{s \in S_i}
\sum_{t \in S_i} \mathrm{E}\bigl[\theta(s,X_s)
\theta(t,X_t)\bigr]
\\
& =& (\bQy_{n,\tau,\omega})^{-2} \sum_{|\iota| < m_n}
\gamma _{\iota}(\tau, \tau) \sum_{j=2(i-1) m_n+1+(0\vee\iota
)}^{(2i-1)m_n + (\iota\wedge0)}
c_{j+\iota,k}(\omega) c_{j,k}(\omega)'. %\end{split}
\end{eqnarray*}
Since $| c_{t,k}(\omega) | \leq1$ and
\[
\sum_{\iota=-\infty}^{\infty} \bigl|\gamma_{\iota}(
\tau,\tau)\bigr| \leq 1+ C_1 \mathop{ \sum_{\iota=-\infty}}_{ \iota\neq0}^{\infty}
\iota ^{-\delta} =: C < \infty,
\]
we have
\[
\sum_{i=1}^{\mu_n} \Var \biggl(\sum
_{t \in S_i} \theta(t,X_t) \biggr) \leq4 C
\bigl(f(q_{\tau})\bigr)^{-2} n =: V_n^2.
\]

Direct calculations show that $E_n$ defined in \eqref{def:en} of the
independent blocks argument satisfies $ E_n = \mathrm{o}(1)$. This completes the
independent blocks argument and concludes the proof.
\end{pf}

%leA.3 #&#
\begin{lemma} \label{lem:varzn2}
 For the Fourier frequencies $\omega\in\mathcal{F}_n$, let
%
%eA.5 #&#
\begin{equation}
\label{lem:varzn2:eqn:Ht} H_t(\bdelta;\tau,\omega) := \int
_0^{n^{-1/2} \mathbf{c}^\tr
_{t}(\omega) \bdelta} \bigl(I\{X_t \leq s +
q_{\tau}\} - I\{X_t \leq q_{\tau}\}\bigr)\,
\mathrm{d}s
\end{equation}
and define
%
%eA.6 #&#
\begin{equation}
\label{def:wtn} W_{t,n}(\omega,\bdelta) := H_t(\bdelta;
\tau,\omega) - {f(q_{\tau
})}(2n)^{-1}\bigl(
\mathbf{c}_{t}^\tr(\omega) \bdelta\bigr)^2.
\end{equation}
Then, for some finite constant $C$ (independent of $t,t_1,t_2$) and $n$
large enough,
%
%eA.7 #&#
\begin{eqnarray} \label{W1}
%{eqnarray}
\sup_{\omega\in\mathcal{F}_n} \sup_t \bigl|
\mathrm{E}\bigl[W_{t,n}(\omega ,\bdelta)\bigr]\bigr| %&\leq&
&\leq& C\|
\bdelta\|^3 n^{-3/2},\nonumber  \\[-8pt]\\[-8pt]
\sup
_{\omega\in\mathcal{F}_n} \sup_t \bigl|W_{t,n}(\omega,
\bdelta)\bigr| %&\leq&
&\leq& C \bigl(n^{-1/2} \|\bdelta\| + n^{-1}
\|\bdelta\|^2\bigr) \nonumber% \label{W2},
\end{eqnarray}
almost surely, and
%
%eA.8 #&#
\begin{equation}
\sup_{\omega\in\mathcal{F}_n} \bigl|\mathrm{E}\bigl[W_{t_1,n}(\omega ,
\bdelta)W_{t_2,n}(\omega,\bdelta)\bigr]\bigr| %&\leq&
\leq C \bigl(\|
\bdelta\|^{4}\vee1\bigr) \bigl(n^{-3/2}I\{t_1=t_2
\} + n^{-2}I\{ t_1\neq t_2\}
\bigr)\label{W3}.
\end{equation}
\end{lemma}

\begin{pf}
%Assumption (A2) guarantees that all calculations below still hold when
%distributions depend on $n$.
First, note that
%
%eA.9 #&#
\begin{eqnarray}
\label{eqn:EHt} \mathrm{E} \bigl[H_t(\bdelta; \tau, \omega) \bigr] &
= & \mathrm{E} \biggl[ \int_0^{n^{-1/2} \mathbf
{c}^\tr_{t}(\omega) \bdelta} \bigl(I
\{X_t \leq u + q_{\tau}\} - I\{ X_t \leq
q_{\tau}\} \bigr) \biggr] \,\mathrm{d} u\nonumber
\\[-8pt]\\[-8pt]
& = & \int_0^{n^{-1/2} \mathbf{c}^\tr_{t}(\omega)
\bdelta} \bigl( f(q_{\tau}) u
+ r_4(u,\tau) \bigr) \,\mathrm {d}u = \frac{f(q_{\tau})}{2n} \bigl(
\mathbf{c}_{t}^\tr(\omega) \bdelta \bigr)^2 +
r_1(\tau,\omega)
\nonumber,
\end{eqnarray}
where $|r_4(u,\tau)| \leq C_3 u^2$, hence $|r_1(\omega,\tau)| \leq
C_4 \|\bdelta\|^3n^{-3/2}$. Next, observe that
%
%eA.10 #&#
%eA.11 #&#
\begin{eqnarray}
&&\mathrm{E} \bigl[H_t(\bdelta; \tau, \omega)^2 \bigr]\nonumber
\\
&&\quad = \mathrm{E} \biggl[ \int_0^{n^{-1/2} \mathbf{c}^\tr
_{t}(\omega) \bdelta}
\int_0^{n^{-1/2} \mathbf{c}^\tr_{t}(\omega)
\bdelta} \bigl(I\{X_t \leq u +
q_{\tau}\} - I\{X_t \leq q_{\tau}\} \bigr)
\nonumber
\\
&&\quad\quad\hphantom{\mathrm{E} \biggl[ \int_0^{n^{-1/2} \mathbf{c}^\tr
_{t}(\omega) \bdelta}
\int_0^{n^{-1/2} \mathbf{c}^\tr_{t}(\omega)
\bdelta} \bigl(}{} \times \bigl(I\{X_t \leq v + q_{\tau}\} - I
\{X_t \leq q_{\tau}\} \bigr) \,\mathrm {d}u \,\mathrm{d}v
\biggr]
\nonumber
\\
&&\quad = \mathrm{E} \biggl[ \int_0^{n^{-1/2} \mathbf{c}^\tr_{t}(\omega)
\bdelta} \int
_0^{n^{-1/2} \mathbf{c}^\tr_{t}(\omega) \bdelta} \bigl(I\bigl\{X_t \leq(u
\wedge v) + q_{\tau}\bigr\} - I\bigl\{X_t \leq(u \wedge0) +
q_{\tau}\bigr\}
\nonumber
\\
& &\hphantom{\mathrm{E} \biggl[ \int_0^{n^{-1/2} \mathbf{c}^\tr_{t}(\omega)
\bdelta} \int
_0^{n^{-1/2} \mathbf{c}^\tr_{t}(\omega) \bdelta} \bigl(}\quad\quad{} - I\bigl\{X_t \leq(v \wedge0) + q_{\tau}\bigr\} + I
\{X_t \leq q_{\tau}\}\bigr) \,\mathrm{d}u \,\mathrm{d}v \biggr]
\nonumber
\\
&&\quad= \int_0^{n^{-1/2} \mathbf{c}^\tr_{t}(\omega) \bdelta} \int_0^{n^{-1/2} \mathbf{c}^\tr_{t}(\omega) \bdelta}
(u \wedge v - u \wedge0 - v \wedge0 ) f(q_{\tau}) + r_2(u,v,
\tau) \,\mathrm {d}u \,\mathrm{d}v\quad \label{eqn:EHtsq1}
\\
& &\quad= 3^{-1}{n^{-3/2}} f(q_{\tau}) \bigl\llvert
\mathbf{c}^\tr_{t}(\omega ) \bdelta\bigr\rrvert
^3 + r_3(\omega,\tau), \label{eqn:EHtsq2} %\end{split}
\end{eqnarray}
where $|r_2(u,v,\tau)| \leq C_1 (u^2 + v^2)$, hence $|r_3(\omega,\tau
)| \leq C_2 \|\bdelta\|^4n^{-2} $. Equality (\ref{eqn:EHtsq1})
follows via a Taylor expansion,
% &   F((u \wedge v) + q_{\tau}) - F((u \wedge0) + q_{\tau}) -
%F((v \wedge0) + q_{\tau}) + F(q_{\tau}) \\
% & = F(q_{\tau}) + (u \wedge v) f(q_{\tau}) + \frac{1}{2} (u \wedge
%v)^2 f'(\xi_1(u,v)) \\
% &  -\Big(F(q_{\tau}) + (u \wedge0) f(q_{\tau}) + \frac{1}{2} (u
% &  -\Big(F(q_{\tau}) + (v \wedge0) f(q_{\tau}) + \frac{1}{2} (v
% & = (u \wedge v - u \wedge0 - v \wedge0) f(q_{\tau}) + r_2(u,v,
(\ref{eqn:EHtsq2}) from the fact that
$\int_0^x\int_0^x  (u \wedge v - u \wedge0 - v \wedge0 )
\,\mathrm{d}u \,\mathrm{d}v = \frac{1}{3} |x|^3$. Similarly, for $t_1
\neq t_2$, but from the same block (otherwise $H_{t_1}$ and $H_{t_2}$
are independent and the previously derived approximation of their
expectations can be used for the proof),
%eA.12 #&#
\begin{eqnarray} \label{eqn:EHt1Ht2}
&&\mathrm{E} \bigl[H_{t_1}(\bdelta; \tau, \omega) H_{t_2}(
\bdelta; \tau, \omega) \bigr]\nonumber
\\
&&\quad = \mathrm{E} \biggl[ \int_0^{n^{-1/2} \mathbf{c}_{t_1}^\tr(\omega
) \bdelta}
\int_0^{n^{-1/2} \mathbf{c}_{t_2}^\tr(\omega) \bdelta} \bigl(I\{X_{t_1} \leq u +
q_{\tau}\} - I\{X_{t_1} \leq q_{\tau}\} \bigr)
\nonumber
\\
&&\hphantom{\mathrm{E} \biggl[ \int_0^{n^{-1/2} \mathbf{c}_{t_1}^\tr(\omega
) \bdelta}
\int_0^{n^{-1/2} \mathbf{c}_{t_2}^\tr(\omega) \bdelta} \bigl(}\quad\quad{} \times \bigl(I\{X_{t_2} \leq v + q_{\tau}\} - I
\{X_{t_2} \leq q_{\tau}\} \bigr)\, \mathrm{d}u\, \mathrm{d}v
\biggr]\quad\quad
\nonumber
\\
%[-8pt]\\[-8pt]
% \end{split}
% \begin{equation}
&&\quad = \int_0^{n^{-1/2} \mathbf{c}_{t_1}^\tr(\omega) \bdelta}
\int_0^{n^{-1/2} \mathbf{c}_{t_2}^\tr(\omega) \bdelta} F_{t_2-t_1}(u+q_{n,\tau},
v+q_{\tau}) - F_{t_2-t_1}(q_{\tau}, v+q_{\tau}
\nonumber
)
\\
&&\hphantom{\int_0^{n^{-1/2} \mathbf{c}_{t_1}^\tr(\omega) \bdelta}
\int_0^{n^{-1/2} \mathbf{c}_{t_2}^\tr(\omega) \bdelta}}\quad\quad{} - F_{t_2-t_1}(u+q_{\tau}, q_{\tau}) +
F_{t_2-t_1}(q_{n,\tau}, q_{\tau
}) \,\mathrm{d}u \,\mathrm{d}v
\nonumber
\\
&&\quad =\int_0^{n^{-1/2} \mathbf{c}_{t_1}^\tr(\omega) \bdelta} \int_0^{n^{-1/2} \mathbf{c}_{t_1}^\tr(\omega) \bdelta}
r_6(u,v,\tau) \,\mathrm{d}u \,\mathrm{d}v%\\
% & =&
= r_7(\omega,\tau),%\nonumber
\end{eqnarray}
where $|r_6(u,v,\tau)| \leq C_6 (u^2 + v^2)$, hence $|r_7(u,v,\tau)|
\leq C_7 \|\bdelta\|^4 n^{-2}$;
equality (\ref{eqn:EHt1Ht2}) follows via a Taylor expansion and some
straightforward algebra.
%From (\ref{eqn:EHtsq1}), (\ref{eqn:EHtsq2}) and (\ref{eqn:EHt1Ht2}),
%we obtain
%f_{n,X}(q_{n,\tau}) \sum_{t=1}^n \left| \mathbf{c}^\tr_{t}(
%with $|R_n| \leq C_3 \|\bdelta\|^4m_nn^{-1}$ for sufficiently large
%$n$, and
%For the sum in the above statement note that
% \leq\sum_{t=1}^n \left| \|\mathbf{c}_{t}(\omega_{j,n})\|_{\infty} \
% \leq3^{3/2} \sum_{t=1}^n \|\bdelta\|^3 = 3^{3/2} n\|\bdelta\|^3.\]
This completes the proof.
\end{pf}

%sA.2 #&#
\subsection{Details for the proof of (\texorpdfstring{\protect\ref{thm:uniflinearrank}}{3.20})}

We now turn to the proof of Theorem~\ref{PropCLAsymp}. Sections~\ref{subsub213}--\ref{subsub212} contain the proofs of \eqref{h3}
and  \eqref{h4}, which are basic in establishing that theorem. %(
Some auxiliary results used in the proofs are collected in Section~\ref{subsub211} under the form of Lemmas \ref{lem:df2} and \ref{lem:nummn}.
Denote by $\hat F_{n}$ the empirical distribution function of
$Y_1,\ldots,Y_n$. Throughout this section, the notation from
Section~\ref{s3.2} is used. %the proof of Theorem~\ref{PropCLAsymp} is used.

%sA.2.1 #&#
\subsubsection{Proof of \texorpdfstring{\protect\eqref{h3}}{(3.22)}} \label{subsub213}

Plugging into \eqref{h3} the definition of $\Zyrh{n,\tau, \omega
,1}(\bdelta)$, it remains to show that [recall that $c_{t,1}(\omega)
= 1$]
%
%eA.13 #&#
\begin{eqnarray}\label{erg1}
&&\max_{k=2,3}\sup_{\omega\in\mathcal{F}_n} \Biggl| n^{-1/2}
\sum_{t=1}^n c_{t,k}(\omega)
\bigl(I\bigl\{U_{t} \leq F\bigl(\hat F_{n}^{-1}(
\tau)\bigr)\bigr\} - I\{U_{t} \leq\tau\} \bigr) \Biggr|\nonumber
\\[-8pt]\\[-8pt]
 &&\quad= \mathrm{O}_{\mathrm{P}}
\bigl( n^{-1/4}m_n^{1/2}\log n\bigr) \nonumber
\end{eqnarray}
 and
%
%eA.14 #&#
\begin{eqnarray}\label{erg2}
&& \Biggl| n^{-1/2} \sum_{t=1}^n
\bigl(I\bigl\{U_{t} \leq F\bigl(\hat F_{n}^{-1}(
\tau )\bigr)\bigr\} - I\{U_{t} \leq\tau\}\bigr) - \sqrt{n}\bigl(F
\bigl(\hat F_{n}^{-1}(\tau)\bigr) - \tau\bigr) \Biggr|
\nonumber
\\[-8pt]\\[-8pt]
&&\quad = \mathrm{O}_{\mathrm{P}}\bigl( n^{-1/4} m_n^{1/2}
\log n\bigr). \nonumber
\end{eqnarray}
First, consider \eqref{erg1}. Since, by Lemma~\ref{lem:df2}, $|F(\hat
F_{n}^{-1}(\tau)) - \tau| = \mathrm{O}_{\mathrm{P}}(n^{-1/2}\sqrt{\log n})$,
we obtain
%
%eA.15 #&#
\begin{eqnarray}\label{secondterm}
&& \sup_{\omega\in\mathcal{F}_n} \Biggl| n^{-1/2} \sum
_{t=1}^n c_{t,k}(\omega) \bigl(I\bigl
\{U_{t} \leq F\bigl(\hat F_{n}^{-1}(\tau)\bigr)
\bigr\} - I\{U_{t} \leq\tau\}\bigr) \Biggr|
\nonumber
\\
&& \quad\leq\sup_{\omega\in\mathcal{F}_n} n^{-1/2}\sup_{|x-\tau|\leq n^{-1/2}\log n}
\Biggl|\sum_{t=1}^n c_{t,k}(\omega)
\bigl(I\{U_{t} \leq x\} - I\{U_{t} \leq\tau\} - (x-\tau)
\bigr) \Biggr|\quad\quad
 \\
&& \quad\quad{}+ \sup_{\omega\in\mathcal{F}_n} n^{-1}\log n \Biggl|\sum
_{t=1}^n c_{t,k}(\omega) \Biggr|\nonumber
%+ \sup_{\omega\in\mathcal{F}_n} n^{-1/2}\sup_{|x-\tau|\leq n^{-1/2}
\end{eqnarray}
for $k=2,3$, with probability tending to one. The second term in (\ref
{secondterm}) vanishes, because, for all $\omega\in\mathcal{F}_n$,
$\sum_{t=1}^n \cos(\omega t) = \sum_{t=1}^n \sin(\omega t) = 0$. In
order to bound the first term, cover the set $\mathcal{Z} := \{
u\dvt|u-\tau|\leq n^{-1/2}\log n\}$ with $N < n$ balls of radius $1/n$
and centers $u_1,\ldots,u_N \in\mathcal{Z}$, and define $
\IG_{n,\omega,k}(u) := n^{-1/2}\sum_{t=1}^n c_{t,k}(\omega) (I\{
U_{t} \leq u\} - u)
$. %\]
Then, almost surely,
\begin{eqnarray*}
&& \sup_j \sup_{\omega\in\mathcal{F}_n}\sup
_{|u-u_j|\leq n^{-1}} \bigl|\IG_{n,\omega,k}(u) - \IG_{n,\omega,k}(u_j)
\bigr|
\\
&&\quad \leq\sup_{u \in\mathcal{Z}}n^{-1/2}\sum
_{t=1}^n \bigl( I\bigl\{ U_{t} \leq u +
2n^{-1}\bigr\} - I\bigl\{U_{t} \leq u - 2n^{-1}
\bigr\} + 4n^{-1} \bigr) + \mathrm{O}\bigl(n^{-1/2}\bigr)
\\
&& \quad\leq\sqrt{n} \sup_{j=1,\ldots,N} \bigl|\hat F_{n,U}
\bigl(u_j+2n^{-1}\bigr) - \hat F_{n,U}
\bigl(u_j - 2n^{-1}\bigr) - 4n^{-1} \bigr| + \mathrm{O}
\bigl(n^{-1/2}\bigr) %= O_{\mathrm{P}}(n^{-1/2})
,
\end{eqnarray*}
where the latter bound, in view of Lemma~\ref{lem:df}, is $\mathrm{O}_{\mathrm
{P}}(n^{(1-\delta)/(2+2\delta)} \log n)$.
Thus, % we have obtained the estimate
\[
\sup_j \sup_{\omega\in\mathcal{F}_n}\sup
_{|u-u_j|\leq n^{-1}} \bigl|\IG_{n,\omega,k}(u) - \IG_{n,\omega,k}(u_j)
\bigr| = \mathrm{O}_{\mathrm{P}}\bigl(n^{(1-\delta)/(2+2\delta)} \log n\bigr),\quad\quad k=2,3,
\]
and therefore
%
%eA.16 #&#
\begin{eqnarray} \label{maxest}
&&\max_{k=2,3}\sup_{\omega\in\mathcal{F}_n} \Biggl| n^{-1/2}
\sum_{t=1}^n c_{t,k}(\omega)
\bigl(I\bigl\{U_{t} \leq F\bigl(F_{n}^{-1}(\tau)
\bigr)\bigr\} - I\{U_{t} \leq\tau\}\bigr) \Biggr|
\nonumber
\\[-8pt]\\[-8pt]
&&\quad\leq \max_{k=2,3}\sup_{j=1,\ldots,N} \sup
_{\omega\in\mathcal{F}_n} \bigl|\IG_{n,\omega,k}(u_j) -
\IG_{n,\omega,k}(\tau) \bigr| +\mathrm{ O}_{\mathrm
{P}}\bigl(n^{(1-\delta)/(2+2\delta)} \log n
\bigr).\nonumber
\end{eqnarray}
Now, by construction, $\max_j |u_j-\tau|\leq n^{-1/2}\log n$.

Moreover,
\[
\max_{k=2,3}\sup_{j=1,\ldots,N} \sup
_{\omega\in\mathcal{F}_n} \bigl|\IG_{n,\omega,k}(u_j) -
\IG_{n,\omega,k}(\tau) \bigr| = \sup_{\theta\in\Theta_n} \Biggl| \sum
_{t=1}^n \theta(t,U_t) \Biggr|,
\]
where
\begin{eqnarray*}
\Theta_n &:=& \bigl\{\theta(t,u) := n^{-1/2}
c_{t,k}(\omega) \bigl(I\{u \leq u_j\} - I\{u \leq\tau\} -
(u_j - \tau ) \bigr) |
\\
&&\hphantom{\bigl\{} \omega\in\mathcal{F}_n, j=1,\ldots,N, k=2,3\bigr\}.
\end{eqnarray*}

Apply the independent blocks argument with $\eta_n := \tilde C
n^{-1/2}\sqrt{\log n}(n^{1/2}m_n\log n)^{1/2}$, where $\tilde C$ is a
large enough constant, and $\Theta_n$ defined above.
Direct calculations show that
\[
\sup_{\theta\in\Theta_n} \bigl |\theta(t, U_t)\bigr | \leq2 n^{-1/2} =: C_n
\qquad\mbox{a.s.},
\]
 which yields (i) from the independent blocks argument, since $m_n
C_n \sim m_nn^{-1/2}\log n \ll\eta_n$.
Additionally, for some finite constant $C$ independent of $\theta\in
\Theta_n$
$\mathrm{E}|\theta(t,U_t)|^2 \leq C n^{-3/2} \log n$,
and
$\mathrm{E}[\theta(t_1,U_{t_1}) \theta(t_2,U_{t_2})] \leq C
n^{-2}(\log n)^2$,
and thus
\begin{eqnarray*}
\sup_{\theta\in\Theta_n} \sum_{j=1}^{\mu_n}
\Var \biggl(\sum_{t
\in S_j} \theta(t, U_t)
\biggr) &\leq&\bar C n^{-1/2} \log n =: V_n^2,
\\ \sup
_{\theta\in\Theta_n} \sum_{j=1}^{\mu_n}
\Var \biggl(\sum_{t \in T_j} \theta(t, U_t)
\biggr) &\leq& V_n^2.
\end{eqnarray*}
Hence, (ii) from the independent blocks argument holds and the fact
that $E_n = \mathrm{o}(1)$ with $E_n$ defined in \eqref{def:en} follows by a
simple calculation. The independent blocks argument thus yields
\[
\max_{k=2,3}\sup_{j=1,\ldots,N} \sup
_{\omega\in\mathcal{F}_n} \bigl|\IG_{n,\omega,k}(u_j) -
\IG_{n,\omega,k}(\tau) \bigr| = \mathrm{O}_{\mathrm
{P}}\bigl(n^{-1/4}m_n
\log n\bigr) = \mathrm{O}_{\mathrm{P}}\bigl(n^{(1-\delta)/(4+4\delta)} \log n\bigr).
\]
Together with (\ref{maxest}), this establishes (\ref{erg1}). Turning
to \eqref{erg2}, Lemmas \ref{lem:df2} and \ref{lem:df} yield
\begin{eqnarray*}
&& \Biggl| n^{-1/2} \sum_{t=1}^n
\bigl(I\bigl\{U_{t} \leq F\bigl(\hat F_{n}^{-1}(
\tau)\bigr)\bigr\} - I\{U_{t} \leq\tau\} - \bigl(F\bigl(\hat
F_{n}^{-1}(\tau )\bigr) - \tau\bigr) \bigr) \Biggr|
\\
&& \quad\leq\sup_{|u - \tau|\leq n^{-1/2}\log n} \Biggl| n^{-1/2} \sum
_{t=1}^n \bigl(I\{U_{t} \leq u\} - I
\{U_{t} \leq\tau\} - (u - \tau) \bigr) \Biggr|
\\
&&\quad = {n}^{1/2} \sup_{|u - \tau|\leq n^{-1/2}\log
n} \bigl|\hat F_{n,U}(u)
- \hat F_{n,U}(\tau) - (u-\tau) \bigr|
\\
&&\quad = \mathrm{O}_{\mathrm{P}}\bigl(n^{-1/2}\bigl(m_n\vee
n^{1/4}\bigr)\log n\bigr)\leq \mathrm{O}_{\mathrm{P}}\bigl(n^{-1/4}m_n^{1/2}
\log n\bigr).
\end{eqnarray*} %\upqed\qed

%sA.2.2 #&#
\subsubsection{Proof of % the estimate
\texorpdfstring{\protect\eqref{h4}}{(3.23)}} \label{subsub212}

%With the notation $U_{t,n} = F_{n,X}(X_{t,n}) \sim{\mathcal U} [0,1]$,
Observe the decomposition
\begin{eqnarray*}
&&\Zyrh{n,\tau, \omega,2}(\bdelta) - \sum_{t=1}^n
\int_0^{n^{-1/2}{\mathbf{c}_{t}^\tr(\omega)\bdelta}} \bigl( I\{U_{t} \leq s+
\tau\} - I\{U_{t} \leq\tau\} \bigr)\,\mathrm{d}s
\\
&&\quad= \sum_{t=1}^n \int
_0^{n^{-1/2}{\mathbf{c}_{t}^\tr(\omega)\bdelta
}} \bigl( I\bigl\{U_{t} \leq F
\bigl(\hat F_{n}^{-1}(s+\tau)\bigr)\bigr\} - I\bigl
\{U_{t} \leq F\bigl(\hat F_{n}^{-1}(\tau)\bigr)
\bigr\} \\
&&\hphantom{\sum_{t=1}^n \int
_0^{n^{-1/2}{\mathbf{c}_{t}^\tr(\omega)\bdelta
}} \bigl(}\quad\quad{}- I\{U_{t} \leq s+\tau\}
 + I\{U_{t} \leq\tau\} \bigr) \,\mathrm{d}s
\\
&&\quad= \int_{-2\|\bdelta\|}^{2\|\bdelta\|} n^{-1/2}\sum
_{t=1}^n \bigl( I\bigl\{U_{t} \leq F
\bigl(\hat F_{n}^{-1}\bigl(n^{-1/2}s+\tau\bigr)\bigr)
\bigr\} - I\bigl\{U_{t} \leq F\bigl(\hat F_{n}^{-1}(
\tau)\bigr)\bigr\}
\\
&&\hphantom{\int_{-2\|\bdelta\|}^{2\|\bdelta\|} n^{-1/2}\sum
_{t=1}^n \bigl(}\quad\quad{} - I\bigl\{U_{t} \leq n^{-1/2}s+\tau\bigr\} + I
\{U_{t} \leq\tau\} \bigr)
\\
&&\hphantom{\int_{-2\|\bdelta\|}^{2\|\bdelta\|} n^{-1/2}\sum
_{t=1}^n}\quad\quad{}\times\bigl(I\bigl\{0 \leq s \leq
\mathbf{c}_{t}^\tr(\omega)\bdelta\bigr\} - I\bigl\{0 \geq s
\geq\mathbf{c}_{t}^\tr(\omega)\bdelta\bigr\} \bigr)
\,\mathrm{d}s
\\
&&\quad= A_n^{(1)} - A_n^{(2)} -
A_n^{(3)}+ A_n^{(4)}, \quad\quad\mbox{say,}
\end{eqnarray*}
where
\begin{eqnarray*}
% & \Zxrh{n,\tau, \omega,2}(\bdelta) - \sum_{t=1}^n\int_0^{n^{-1/2}{
A_n^{(1)}
& :=& \int_{-2 \|\bdelta\|}^{2 \|\bdelta\|} \bigl(\IS _{n,\omega,\bdelta}^{(+)}
\bigl(F \bigl(\hat F_{n}^{-1}\bigl(n^{-1/2} s + \tau
\bigr)\bigr), n^{-1/2} s + \tau;s\bigr) - \IS_{n,\omega,\bdelta}^{(+)}
\bigl(F \bigl(\hat F_{n}^{-1}(\tau)\bigr), \tau;s\bigr)
\bigr) \,\mathrm{d}s ,
\\
A_n^{(2)} & :=& \int_{-2 \|\bdelta\|}^{2 \|\bdelta\|}
n^{-1/2} \sum_{t=1}^n \bigl[
\bigl(F\bigl(\hat F_{n}^{-1} \bigl(n^{-1/2} s + \tau
\bigr)\bigr) - \bigl(n^{-1/2} s + \tau\bigr)\bigr) - \bigl(F \bigl(\hat
F_{n}^{-1}(\tau)\bigr) - \tau\bigr) \bigr]
\\
&&\hphantom{\int_{-2 \|\bdelta\|}^{2 \|\bdelta\|}
n^{-1/2} \sum_{t=1}^n}{} \times I\bigl\{0 \leq s \leq\mathbf{c}_{t}^\tr (\omega)
\bdelta\bigr\} \,\mathrm{d}s,
\\
A_n^{(3)} & :=& \int_{-2 \|\bdelta\|}^{2 \|\bdelta\|}
\bigl(\IS _{n,\omega,\bdelta}^{(-)}\bigl(F \bigl(\hat F_{n}^{-1}
\bigl(n^{-1/2} s + \tau\bigr)\bigr), n^{-1/2} s + \tau;s\bigr) -
\IS_{n,\omega,\bdelta}^{(-)}\bigl(F \bigl(\hat F_{n}^{-1}(
\tau)\bigr), \tau;s\bigr) \bigr) \,\mathrm{d}s ,
\\
A_n^{(4)} & :=& \int_{-2 \|\bdelta\|}^{2 \|\bdelta\|}
n^{-1/2} \sum_{t=1}^n \bigl[
\bigl(F\bigl(\hat F_{n}^{-1} \bigl(n^{-1/2} s + \tau
\bigr)\bigr) - \bigl(n^{-1/2} s + \tau\bigr)\bigr) - \bigl(F \bigl(\hat
F_{n}^{-1}(\tau)\bigr) - \tau\bigr) \bigr]
\\
&&\hphantom{\int_{-2 \|\bdelta\|}^{2 \|\bdelta\|}
n^{-1/2} \sum_{t=1}^n}{} \times I\bigl\{0 \geq s \geq\mathbf{c}_{t}^\tr (\omega)
\bdelta\bigr\} \,\mathrm{d}s, %\end{split}
\end{eqnarray*}
and %we have introduced the notation
\begin{eqnarray*}
\IS_{n,\omega,\bdelta}^{(+)}(u,v;s) & :=& n^{-1/2}\sum
_{t=1}^n \bigl( I\{U_{t} \leq u\} - I
\{U_{t} \leq v\} - (u-v) \bigr) I\bigl\{0 \leq s \leq
\mathbf{c}_{t}^\tr(\omega)\bdelta\bigr\},
\\
\IS_{n,\omega,\bdelta}^{(-)}(u,v;s) & :=
&n^{-1/2}\sum_{t=1}^n \bigl( I
\{U_{t} \leq u\} - I\{U_{t} \leq v\} - (u-v) \bigr) I\bigl
\{0 \geq s \geq\mathbf{c}_{t}^\tr(\omega)\bdelta\bigr\}.
\end{eqnarray*}
%
%Since ... we have, with probability tending to one,
%&\sup_{\omega\in\mathcal{F}_n}\sup_{\|\bdelta\| \leq A\sqrt{\log n}}
%d_n} \Big(|\IS_{n,\omega,\bdelta}^{(+)}(x,\tau)| + |\IS_{n,\omega,
%&\ +2 \sup_{|x-q_{\tau,m_n}|\leq2\|\bdelta\|n^{-1/2}}\sqrt{n}\Big|
First note that, in view of Lemma~\ref{lem:df2},
\begin{eqnarray*}
\bigl|A_n^{(2)} \bigr| &\leq&4 \| \bdelta\| \sqrt{n} \sup
_{|u - \tau| \leq2 \|
\bdelta\| / \sqrt{n}} \bigl| F\bigl(\hat F_{n}^{-1} (u)\bigr)
- u - \bigl(F \bigl(\hat F_{n}^{-1}(\tau)\bigr) - \tau\bigr)
\bigr|
\\
&=& \mathrm{O}_{\mathrm{P}}\bigl(\rho_n\bigl(2 (\log n)^{1/2}
n^{-1/2}, \delta\bigr) \sqrt {n\log n}\bigr)
\\
&=& \mathrm{O}_{\mathrm{P}}\bigl(
\bigl(n^{-1/4} (\log n)^{5/4}\bigr) \vee\bigl(n^{(1-\delta
)/(2+2\delta)} (
\log n)^{3/2}\bigr)\bigr)
\\
&=& \mathrm{O}_{\mathrm{P}}\bigl(n^{-1/4} m_n^{1/2} \log n
\bigr).
\end{eqnarray*}
A similar bound can be obtained for $A_n^{(4)} $. Next, for
sufficiently large $n$, still in view of Lemma~\ref{lem:df2},
\begin{eqnarray*}
&& \int_{-2 \|\bdelta\|}^{2 \|\bdelta\|} \bigl|
\IS_{n,\omega
,\bdelta}^{(+)}\bigl(F \bigl(\hat F_{n}^{-1}
\bigl(n^{-1/2} s + \tau\bigr)\bigr), n^{-1/2} s + \tau;s\bigr) \bigr|\,
\mathrm{d}s
\\
&&\quad \leq\int_{-2 \|\bdelta\|}^{2 \|\bdelta\|} \sup_{v: |v-\tau| \leq2 \|\bdelta\| / \sqrt{n}}
\bigl| \IS _{n,\omega,\bdelta}^{(+)}\bigl(F \bigl(\hat F_{n}^{-1}(v)
\bigr), v;s\bigr) \bigr|\, \mathrm {d}s
\\
&&\quad \leq\int_{-2 \|\bdelta\|}^{2 \|\bdelta\|} \sup_{v: |v-\tau| \leq2 \|\bdelta\| / \sqrt{n}}
\sup_{u: |u-v|
\leq n^{-1/2}\log n } \bigl| \IS_{n,\omega,\bdelta}^{(+)}(u, v;s) \bigr|\,
\mathrm{d}s
\\
&&\quad \leq4 \|\bdelta\| \sup_{s : |s| \leq2 \|
\bdelta\|} \sup_{v: |v-\tau| \leq2 \|\bdelta\| / \sqrt{n}}
\sup_{u: |u-v| \leq n^{-1/2}\log n } \bigl| \IS_{n,\omega,\bdelta
}^{(+)}(u, v;s) \bigr|.
\end{eqnarray*}
Similar inequalities hold for $\int_{-2 \|\bdelta\|}^{2 \|\bdelta\|}
| \IS_{n,\omega,\bdelta}^{(+)}(F (\hat F_{n}^{-1}(\tau)),\tau
;s)  | \,\mathrm{d}s$. Let us show that
%
%eA.17 #&#
\begin{eqnarray}
\label{eq:sup1} &&\sup_{\omega\in\mathcal{F}_n} \sup_{\bdelta: \|\bdelta\| \leq
A \sqrt{\log n}} \sup
_{s : |s| \leq2 \|\bdelta\|} \mathop{\sup_{(u,v): |v-\tau| \leq2 \|\bdelta\| / \sqrt
{n}}}_{|u-v| \leq n^{-1/2}\log n } \bigl|
\IS_{n,\omega,\bdelta
}^{(+)}(u, v;s) \bigr|\nonumber
\\[-8pt]\\[-8pt]
&&\quad = \mathrm{O}_{\mathrm{P}}
\bigl(n^{-1/4}m_n^{1/2}\log n\bigr).\nonumber
\end{eqnarray}

For any $C>0$, we have $I\{0 \leq s \leq\mathbf{c}_t^\tr\bdelta\} =
I\{0 \leq Cs \leq C\mathbf{c}_t^\tr\bdelta\}$.
Thus, it is sufficient to consider vectors $\bdelta$ satisfying $\|
\bdelta\| = 1$. Since, by definition, $\|\mathbf{c}_{t}(\omega)\| =
\sqrt{2}$, it also is sufficient to consider values of $s$ in the
interval $[0,\sqrt2]$. Finally, note that if
\[
I\bigl\{0 \leq s_1 \leq\mathbf{c}_t^\tr
\bdelta_1\bigr\} = I\bigl\{0 \leq s_2 \leq
\mathbf{c}_t^\tr\bdelta_2\bigr\}\quad\quad \mbox{for all
} t=1,\ldots,n,
\]
then also $\IS_{n,\omega,\bdelta_1}^{(+)}(u,v;s_1)=\IS_{n,\omega
,\bdelta_2}^{(+)}(u,v;s_2)$. We thus can rewrite (\ref{eq:sup1}) as
%
%eA.18 #&#
\begin{equation}
\label{eq:sup2} G_n := \sup_{T \in\mathcal{M}_n} \mathop{\sup
_{ (u,v): |v-\tau| \leq2
\|\bdelta\| / \sqrt{n}}}_{ |u-v| \leq n^{-1/2}\log n } \bigl|\bar\IS _n^{(+)}(u,v;T)\bigr|
= \mathrm{O}_{\mathrm{P}}\bigl(n^{-1/4}m_n^{1/2}\log n
\bigr),
\end{equation}
where
%
%eA.19 #&#
\begin{equation}
\label{mn} \mathcal{M}_n:= \bigl\{T = \bigl\{t \in\{1,\ldots,n\}\dvt 0
\leq s \leq\mathbf {c}_t^\tr\bdelta\bigr\} | \omega\in
\mathcal{F}_n, s \in[0,\sqrt 2], \|\bdelta\| = 1 \bigr\}
\end{equation}
and
\[
\bar\IS_n^{(+)}(u,v;T)  :=  n^{-1/2} \sum
_{t\in T} \bigl( I\{U_{t} \leq u\} - u - \bigl(I
\{U_{t} \leq v\} - v\bigr) \bigr).
\]
In order to prove \eqref{eq:sup1} (equivalently, \eqref{eq:sup2}),
define the set
\[
\mathcal{Z}_n := \bigl\{ (u,v) \in\IR^2  \dvt |u-v| \leq
n^{-1/2}\log n , |v-\tau| \leq2 A n^{-1/2}\sqrt{\log n} \bigr\}
\]
and cover it with $N<n^2$ balls of radius $1/n$ with centers
$z_1,\ldots,z_N \in\mathcal{Z}_n$.
For any $(u,v)$ in $\mathcal{Z}_n$ there exists %an index
a $j$ such that $\|(u,v) - (z_{1j}, z_{2j}) \| \leq1/n$ and,
%that for this index we have (
letting $z_j:=(z_{1j},z_{2j})$, we have, almost surely,
\begin{eqnarray*}
\rho(u,v,z_j) &:= & \bigl| \bar\IS_n^{(+)}(u,v;T) -
\bar\IS _n^{(+)}(z_{1j}, z_{2j};T)\bigr|
\\
& \leq& n^{-1/2} \sum_{t=1}^n
\bigl(I\bigl\{|U_{t} - z_{1j}| \leq n^{-1}\bigr\}
+ I\bigl\{|U_{t} - z_{2j}| \leq n^{-1}\bigr\} +
|u-z_{1j}|+ |v-z_{2j}| \bigr)
\\
& \leq&2n^{-1/2} + n^{-1/2} \sum_{t=1}^n
\bigl(I\bigl\{|U_{t} - z_{1j}| \leq n^{-1}\bigr
\} + I\bigl\{|U_{t} - z_{2j}| \leq n^{-1}\bigr\}
\bigr)
\\
& =& 2n^{-1/2} + n^{-1/2} \sum_{t=1}^n
\bigl(I\bigl\{U_{t} \leq z_{1j} + n^{-1}\bigr\}
- I\bigl\{U_{t} < z_{1j} - n^{-1}\bigr\}
\\
&&\hphantom{2n^{-1/2} + n^{-1/2} \sum_{t=1}^n
\bigl(}{} + I\bigl\{U_{t} \leq z_{2j} + n^{-1}\bigr\} -
I\bigl\{U_{t} < z_{2j} - n^{-1}\bigr\} \bigr)
\\
& \leq& n^{1/2} \bigl(\hat F_{n,U}\bigl(z_{1j} +
2n^{-1}\bigr) - \bigl(z_{1j} + 2n^{-1}\bigr) -
\bigl(\hat F_{n,U}\bigl(z_{1j} - 2n^{-1}\bigr) -
\bigl(z_{1j} - 2n^{-1}\bigr) \bigr)
\\
&&\hphantom{n^{1/2} \bigl(}{}+ \hat F_{n,U}\bigl(z_{2j} + 2n^{-1}\bigr) -
\bigl(z_{2j} + 2n^{-1}\bigr)
\\
&&\hphantom{n^{1/2} \bigl(}{}- \bigl( \hat F_{n,U}
\bigl(z_{2j} - 2n^{-1}\bigr) - \bigl(z_{2j} -
2n^{-1}\bigr) \bigr) \bigr)
+ \mathrm{O}\bigl(n^{-1/2}\bigr),
\end{eqnarray*}
where $\hat F_{n,U}$ denotes the empirical distribution function of
$U_{1},\ldots, U_{n}$.
From Lemma~\ref{lem:df},
\begin{eqnarray*}
&&\sup_{z_1, \ldots, z_N}   \mathop{\sup_{ (u,v) \in[0,1]^2}}_{ \|
z_j - (u,v)\| < n^{-1}}
\bigl| \rho(u,v,z_j) \bigr|
\\
&&\quad\leq n^{1/2} \sup_{z_j \in\mathcal{Z}}
\bigl| \hat F_{n,U}\bigl(z_{1j} + 2n^{-1}\bigr) - \hat
F_{n,U}\bigl(z_{1j}-2n^{-1}\bigr) -
4n^{-1} \bigr|
\\
&&\quad\quad{} + n^{1/2} \sup_{z_j \in\mathcal{Z}} \bigl| \hat F_{n,U}
\bigl(z_{2j} + 2n^{-1}\bigr) - \hat F_{n,U}
\bigl(z_{2j}-2n^{-1}\bigr) - 4n^{-1} \bigr| + \mathrm{O}
\bigl(n^{-1/2}\bigr)
\\
&&\quad = \mathrm{O}_{\rm P} \bigl({m_n}n^{-1/2}\log n \bigr).
\end{eqnarray*}
%
%where the latter rate follows from Lemma~\ref{lem:df}.
With this, we have, for $G_n$ defined in \eqref{eq:sup2},
\[
G_n \leq\sup_{T \in\mathcal{M}_n} \sup_{z_1,\ldots,z_N}
\bigl|\bar\IS _n^{(+)}(z_{1j}, z_{2j} ;T)\bigr|
+ \mathrm{O}_{\mathrm{P}} \bigl( {m_n}n^{-1/2}\log n \bigr).
\]
Note that
\[
\sup_{T \in\mathcal{M}_n} \sup_{z_1,\ldots,z_N} \bigl|\bar\IS
_n^{(+)}(z_{1j}, z_{2j} ;T)\bigr| = \sup
_{\theta\in\Theta_n} \Biggl| \sum_{t=1}^n
\theta(t,U_t) \Biggr|,
\]
where
\begin{eqnarray*}
\Theta_n &:= &\bigl\{ \theta(t,w) := n^{-1/2} I\{t\in T\}
\bigl( I\{w \leq u\} - u - \bigl(I\{w \leq v\} - v\bigr) \bigr) |
\\
&&\hphantom{\bigl\{} (u,v) =
z_1,\ldots,z_N, T \in\mathcal{M}_n \bigr
\}.
\end{eqnarray*}

We apply the independent blocks argument with $\Theta_n$ defined above
and $\eta_n := n^{-1/4} m_n \log n$; note that $|\mathcal{M}_n| \leq
(n+1)^4$ by Lemma~\ref{lem:nummn} and $N < n^2$ by construction.

Simple computations yield (recall that $z_j \in\mathcal{Z}$)
%
%eA.20 #&#
%eA.21 #&#
\begin{eqnarray}
%z_{1j} - (I\{U_{t} \leq z_{2j}\} - z_{2j})|
\sup_{\theta\in\Theta_n } \sup_{t=1,\ldots,n}\bigl|
\theta(t,U_t)\bigr|&\leq& 2n^{-1/2}, \label{eqn:V1}
\\
\sup_{\theta\in\Theta_n } \sum_{j=1}^{\mu_n}
\Var \biggl( \sum_{t
\in S_j}\theta(t,U_t)
\biggr) &\leq& C n^{-1/2}\log n =: V_n^2,\nonumber\label{eqn:V3}
\\[-8pt]\\[-8pt]
 \sup
_{\theta\in\Theta_n } \sum_{j=1}^{\mu_n}
\Var \biggl( \sum_{t \in
T_j}\theta(t,U_t)
\biggr) &\leq& V_n^2 .\nonumber
\end{eqnarray}
Thus (i) from the independent blocks argument follows from \eqref
{eqn:V1} since
\[
n^{-1/4} m_n^{1/2} \log n \gg m_nn^{-1/2}.
\]
Moreover, \eqref{eqn:V3} yields (ii), again from the independent
blocks argument. Finally, verify $E_n = \mathrm{o}(1)$ with $E_n$ defined in
\eqref{def:en} by direct calculation to conclude
\[
\sup_{T \in\mathcal{M}_n} \sup_{z_1,\ldots,z_N} \bigl|\bar\IS
_n^{(+)}(z_{1j}, z_{2j} ;T)\bigr| =
\mathrm{O}_{\rm P}\bigl(n^{-1/4} m_n \log n\bigr).
\]
A similar result can be derived for $\IS_{n,\omega,\bdelta}^{(-)}$.
This completes the proof. %\qed

%sA.2.3 #&#
\subsubsection{Two auxiliary lemmas} \label{subsub211}

We now state and prove Lemmas \ref{lem:df2} and \ref{lem:nummn} that
have been used in Sections~\ref{subsub213} and \ref{subsub212}.

%leA.4 #&#
\begin{lemma} \label{lem:df2}
\textup{(i)}
Assume that, for any $\gamma>0$ such that $[\alpha-\gamma,\beta
-\gamma] \subset(0,1)$,
\[
\inf_{u\in[\alpha-\gamma,\beta+\gamma]} f\bigl(F^{-1}(u)\bigr) > 0.
\]
Then,
$\sup_{u \in[\alpha,\beta]} |F(\hat F_{n}^{-1}(u)) - u| =
\mathrm{O}_{\mathrm{P}}(n^{-1/2}\sqrt{\log n})$.

 \textup{(ii)} Define $\rho_n(a_n,\delta):= (\frac{a_n +
n^{1/(1+\delta)} a_n^2 \log n}{n} \log n  )^{1/2} \vee(n^{-\delta
/(1+\delta)} \log n)$.
If $\rho_n(a_n,\delta) $ is $ \mathrm{o}(a_n)$, then
\[
\sup_{u,v \in[\alpha,\beta],|u-v|\leq a_n} \bigl|F\bigl(\hat F_{n}^{-1}(u)
\bigr) - F\bigl(\hat F_{n}^{-1}(v)\bigr) - (u - v)\bigr|
=\mathrm{O}_{\mathrm{P}} \bigl(\rho_n(2a_n,\delta) \bigr).
\]
\end{lemma}

\begin{pf}
Elementary analytic
considerations show that, for any nondecreasing function $g$, $\sup_{w\in[u,v]}|g(w)-w| \leq a_n$ implies $\sup_{w\in
[u+2a_n,v-2a_n]}|g^{-1}(w)-w| \leq a_n$. This, for $g(w) = \hat
F_{n}(F^{-1}(w))$, $u = \alpha-\delta$, and $v = \beta+\delta$,
along with Lemma~\ref{lem:df}, yields part (i) of the lemma.
Turning to part (ii), by Lemma~\ref{lem:df}, for any bounded $\mathcal
{Y}\subset\IR$,
\[
\sup_{y \in\mathcal{Y}}\sup_{|x| \leq a_n} \bigl|\hat
F_{n}(y+x) - \hat F_{n}(y) - F(x+y) + F(y)\bigr| =
\mathrm{O}_{\mathrm{P}}\bigl(\rho_n(a_n,\delta)\bigr).
\]
Since, for any $A\subset[0,1]$, $\sup_{u,v\in
A}|F^{-1}(u)-F^{-1}(v)|\leq C_A|u-v|$ for
some positive constant~$C_A$,
%= \Big(\inf_{p\in A} f_{m_n}(F_{m_n}^{-1}(p))\Big)^{-1}$
%which yields
\[
\sup_{u,v \in[\alpha-\gamma,\beta+\gamma], |u-v| \leq a_n} \bigl|\hat F_{n}\bigl(F^{-1}(u)
\bigr) - \hat F_{n}\bigl(F^{-1}(v)\bigr) - (u-v)\bigr| =
\mathrm{O}_{\mathrm{P}}\bigl(\rho_n(a_n,\delta)\bigr).
\]
We now apply Lemma~3.5 from Wendler \cite{Wendler2011}, with $F(w) = \hat
F_{n}(F^{-1}(w))$, $l = a_n$,  $c = D \rho_n(a_n,\delta)$,
$C_1 = \hat F_{n}(F^{-1}(\alpha-\gamma))$, $C_2 = \hat
F_{n}(F^{-1}(\beta+ \gamma))$. By assumption,
$l + 2c = a_n + 2 D \rho_n(a_n,\delta) \leq2 a_n
$ %\]
for sufficiently large $n$. By Lemma~\ref{lem:df}, we have $C_1 =
\alpha+ \delta+ \mathrm{o}_{\mathrm{P}}(1)$, $C_2 = \beta- \delta+
\mathrm{o}_{\mathrm{P}}(1)$ and, for any strictly increasing continuous
function $G$, $(F\circ G^{-1})^{-1} = G\circ F^{-1}$ (see Exercise 3 in
Chapter~1 of Shorack and Wellner \cite{ShorWell1986}); moreover, by part (i) of the
present lemma, $F(\hat F_{n}^{-1}(u))$ is uniformly close to $u$ for
large $n$. Hence,
\[
%&&
\sup_{u,v \in[\alpha,\beta],|u-v|\leq2a_n} \bigl|F\bigl(\hat
F_{n}^{-1}(u)\bigr) - F_{n}\bigl(\hat
F_{n}^{-1}(v)\bigr) - (u - v)\bigr| %\\
%&&
> D
\rho_n(2a_n,\delta)%\Big\}
\]
implies
\[
%&\subset&
% \Big\{
\sup_{u,v \in[\alpha-\delta,\beta+\delta], |u-v| \leq a_n} \bigl|\hat
F_{n}\bigl(F^{-1}(u)\bigr) - \hat F_{n}
\bigl(F^{-1}(v)\bigr) - (u-v)\bigr| %\\
%&&
> D
\rho_n(a_n,\delta)%\Big\},
.
\]
Part (ii) of the lemma follows on letting $D$ tend to infinity.
\end{pf}

%leA.5 #&#
\begin{lemma} \label{lem:nummn}
The cardinality of the set $
\mathcal{M}_n$ defined in \textup{\eqref{mn}} is at most $(n+1)^4$.
\end{lemma}

\begin{pf} Fix a Fourier frequency $\omega_{j,n} = 2\uppi j/n \in
\mathcal{F}_n$ and note that
\[
\mathbf{c}_{t}(\omega_{j,n})^\tr\bdelta=
\delta_1 + \delta_2 \cos (\omega_{j,n} t) +
\delta_3 \sin(\omega_{j,n} t) = \delta_1 +
\sqrt{\delta_2^2+\delta_3^2}
\cos\bigl(\omega_{j,n} t + \phi(\delta_2,
\delta_3)\bigr),
\]
where $\phi(\delta_2,\delta_3) \in[0,2\uppi]$ denotes a phase shift.
Moreover, for any $v \in[0,1]$, noting that the mapping $x \mapsto
\cos(\omega_{j,n} x + \phi)$ is ${n}/{j}$-periodic,
\begin{eqnarray*}
&& \Bigl\{ t \in\{1,\ldots,n\} \big| 0 \leq v\leq\delta_1 + \sqrt
{\delta_2^2+\delta_3^2}\cos(
\omega_{j,n}t + \phi) \Bigr\}
\\
&&\quad = \biggl\{\frac{nk}{j} + w \big| w \in[C_{1,\phi
,v,\bdelta}-C_{0,\phi,v,\bdelta},C_{1,\phi,v,\bdelta}+C_{0,\phi
,v,\bdelta}],k
= 0,\ldots,n \biggr\} \cap\{1,\ldots,n\},
\end{eqnarray*}
where $C_{0,\phi,v,\bdelta}\in[0,n/2j]$ and $ C_{1,\phi,v,\bdelta}
\in[0,n/j]$ denote two real-valued constants (depending on $\phi
,v,\bdelta$). Now, we have
\begin{eqnarray*}
&&\biggl\{\frac{nk}{j} + v \big| v \in[a_1,b_1],k=0,1,\ldots,n
\biggr\} \cap \{1,\ldots,n\}
\\
&&\quad= \biggl\{\frac{nk}{j} + v \big| v
\in[a_2,b_2],k=0,1,\ldots,n \biggr\} \cap\{1,\ldots,n\} ,
\end{eqnarray*}
provided that $\lceil ja_1\rceil= \lceil ja_2\rceil, \lceil
jb_1\rceil= \lceil jb_2\rceil$, where $\lceil a\rceil$ denotes the
smallest integer larger or equal to $a$. The argument above holds for
any Fourier frequency. In particular, it implies that
\begin{eqnarray*}
\mathcal{M}_n&\subset& \biggl\{T = \biggl\{t \in\{1,\ldots,n\} \cap
\biggl\{ \frac{kn}{j} + v \Big| v \in\biggl[\frac{a-b}{j},
\frac{a+b}{j}\biggr] \biggr\} \biggr\} \big|
\\
&&\hphantom{\biggl\{} b = 0,\ldots,\lceil n/2\rceil, a,k=0,\ldots,n, j=1,\ldots,n \biggr\},
\end{eqnarray*}
a collection of sets that contains at most $(n+1)^4$ elements. This
completes the proof.
\end{pf}

\subsection{Two basic lemmas}\label{BasicSec}

Finally, we state and prove here Lemmas \ref{lem:bennett} and \ref
{lem:df}, which have been used at several places in this Appendix.

%leA.6 #&#
\begin{lemma} \label{lem:bennett}
Denote by $X_1,\ldots,X_{\mu_n m_n}$ a sequence of $\mu_n$ independent
blocks of $m_n$ random variables such that
$%\[
\sup_{i=1,\ldots,\mu_n m_n} |X_{i}| \leq C_n$ %
a.s., and
\[
\sum_{j=1}^{\mu_n} \Var \Biggl(\sum
_{i=m_n(j-1)+1}^{m_n j} X_i \Biggr) \leq
V_n^2.
\]
Then, for all $\lambda_n >0$,
$%\[
\mathrm{P} (  |\sum_{j=1}^n X_{j} | > \lambda_n  )
\leq2\exp (-\frac{\log2}{4} (\frac{\lambda_n^2}{2 V_n^2}
\wedge\frac{\lambda_n}{m_nC_n} ) )$. %\]
In particular, for $D>0$,
$%\[
\mathrm{P} (  |\sum_{j=1}^n X_{j} | > 6 \max(D V_n\sqrt {\log n},D^2 m_n C_n \log n)  ) \leq4 n^{-D^2}$.%\]
\end{lemma}

\begin{pf}
Defining the random variables $
U_{k} := \sum_{j= m_n(k-1) +1}^{m_n k} X_{j},   k=1,\ldots, \mu_n
$, %\]
note that $U_{1},U_{2},\ldots, U_{\mu_n}$ are independent, that $|U_{j}|
\leq m_n C_n$ a.s. and that $\Var (\sum_j U_{j}  ) \leq
V_n^2$. Applying Bennett's inequality (see Pollard \cite{Pollard1984}) yields
\begin{eqnarray*}
&&\mathrm{P} \Biggl( \Biggl|\sum_{j=1}^n
X_{j} \Biggr| > \lambda_n \Biggr)
\\
 &&\quad\leq2 \exp \biggl(-
\frac{V_n^2}{m_n^2C_n^2}h \biggl(\frac{m_nC_n\lambda
_n}{2 V_n^2} \biggr) \biggr) \leq2 \exp \biggl(-
\frac{1}{4}\frac{\lambda_n}{m_nC_n}\log \biggl(1+\frac{m_nC_n\lambda_n}{2 V_n^2} \biggr)
\biggr)
\\
&&\quad\leq2 \exp \biggl(-\frac{\log2}{4}\frac{\lambda_n}{m_nC_n} \biggl(
\frac{m_nC_n\lambda_n}{2V_n^2} \wedge1 \biggr) \biggr) = 2 \exp \biggl(-\frac{\log2}{2}
\biggl(\frac{\lambda
_n^2}{4V_n^2}\wedge\frac{\lambda_n}{2m_nC_n} \biggr) \biggr),
\end{eqnarray*}
where $h(x) := (1+x)\log(1+x) - x \geq\frac{1}{2} x\log(1+x) \geq
\frac{\log(2)}{2} x (x \wedge1)$. The second assertion follows by
direct calculation. \end{pf}

%leA.7 #&#
\begin{lemma} \label{lem:df}
Let Assumptions \textup{\ref{(A1)}} and \textup{\ref{(A2)}} hold.
\begin{enumerate}[(ii)]
\item[(i)] Let $\mathcal{Y} \subset\IR$ be a bounded set, $D>1$,
and $0 \leq a_n = \mathrm{o}(1)$. Then,
\[
\sup_{y \in\mathcal{Y}}\sup_{|x| \leq a_n} \bigl|\hat
F_{n}(y+x) - \hat F_{n}(y) - F(x+y) + F(y)\bigr| =
\mathrm{O}_{\mathrm{P}}\bigl(\rho_n(a_n,\delta)\bigr),
\]
where $\rho_n(a_n,\delta):= (\frac{a_n + n^{1/(1+\delta)} a_n^2
\log n}{n} \log n  )^{1/2} \vee(n^{-\delta/(1+\delta)} \log
n)$.

%In particular, if $a_n = O(n^{-\alpha} (\log n)^{\gamma})$, $\alpha
%
\item[(ii)] For any bounded $\mathcal{Y} \subset\mathbb{R}$,
$
\sup_{x \in\mathcal{Y}} |\hat F_{n}(x) - F(x)| = \mathrm{O}_{\mathrm
{P}}(n^{-1/2}\sqrt{\log n})
$.
\end{enumerate}
\end{lemma}

\begin{pf}
The bounded set $Z:=\{(x,y)\in\IR^2|y\in\mathcal{Y}, |x|\leq a_n \}
$ can be covered with $N = \mathrm{O}(n^2)$ spheres of radius $\frac{1}{2}
n^{-1}$ and centers $(z_{1j},z_{2j}) \in Z, j=1,\ldots,N$. A Taylor
expansion yields
\begin{eqnarray*}
&&\sup_{\|(x,y)-(z_{1j},z_{2j})\| \leq1/2 n} \bigl| \hat F_{n}(y+x) - \hat
F_{n}(y) - F(x+y) + F(y)
\\
&&  \hphantom{\sup_{\|(x,y)-(z_{1j},z_{2j})\| \leq1/2 n} \bigl|}{}- \bigl(\hat F_{n}(z_{1j}+z_{2j}) - \hat
F_{n}(z_{2j}) - F(z_{1j}+z_{2j}) +
F(z_{2j})\bigr) \bigr|
\\
&&\quad\leq n^{-1} \sum_{t=1}^n
\bigl(I\bigl\{|Y_{t} - z_{2j}|\leq n^{-1}\bigr\}
+ I\bigl\{\bigl|Y_{t} - (z_{1j} + z_{2j})\bigr|\leq
n^{-1}\bigr\}\bigr) + C n^{-1},
\end{eqnarray*}
where the constant $C$ does not dependent on $t$ and $j$. Therefore,
\begin{eqnarray*}
&&\sup_{y \in\mathcal{Y}}\sup_{|x| \leq a_n} \bigl|\hat
F_{n}(y+x) - \hat F_{n}(y) - F(x+y) + F(y)\bigr|
\\
&&\quad \leq \sup
_{\theta\in\Theta_{1,n}} \Biggl| \sum_{t=1}^n
\theta(t,Y_t) \Biggr| + \sup_{\theta\in\Theta_{2,n}} \Biggl| \sum
_{t=1}^n \theta(t,Y_t) \Biggr|,
\end{eqnarray*}
where
\begin{eqnarray*}
\Theta_{1,n} &:=& \bigl\{ \theta(t,y) := n^{-1} \bigl(I\{y
\leq z_{1j}+z_{2j}\} - I\{y \leq z_{2j}\}\bigr) -
F(z_{1j}+z_{2j}) + F(z_{2j}) |
\\
&&\hphantom{\bigl\{} j=1,\ldots,N
\bigr\},
\end{eqnarray*}
and
\begin{eqnarray*}
\Theta_{2,n} &:=& \bigl\{ \theta(t,y) := n^{-1} \bigl(I\bigl
\{|y - z_{2j}|\leq n^{-1}\bigr\} + I\bigl\{\bigl|y -
(z_{1j} + z_{2j})\bigr|\leq n^{-1}\bigr\}\bigr) + C
n^{-1} |
\\
&&\hphantom{\bigl\{} j=1,\ldots,N \bigr\}.
\end{eqnarray*}

We proceed to bound the suprema over $\Theta_{1,n}$ and $\Theta
_{2,n}$ by applying the independent blocks argument with $\eta_n := D
\rho_n(a_n, \delta)$ and a suitable constant $D$. Begin with $\Theta
_{1,n}$. We have ${\rm E} \theta(t,X_t) = 0$ for all $\theta\in
\Theta_{1,n}$,
$
\sup_{\theta\in\Theta_{1,n}}\sup_t |\theta(t,X_t)| \leq2n^{-1}$,
and
\begin{eqnarray*}
&&\sup_y \sum_{j=1}^{\mu_n}
\Var \biggl( \sum_{t \in S_j} I\{X_t \leq y+x
\} - I\{X_t \leq y\} - F(x+y) + F(y) \biggr)
\\
&&\quad\leq C_2
\mu_n m_n \bigl(m_n |x|^2 + |x|
\bigr) =: V_{1,n}^2
\end{eqnarray*}
for some finite constant $C_2$ independent of $x$, and $m_n := \lceil
n^{1/(1+\delta)} \log n \rceil$, defined as within the independent
blocks argument (see Section~\ref{subsub12x}). The same bound holds
with $S_j$ replaced by~$T_j$. This implies
\[
\sup_{\theta\in\Theta_{1,n}} \sum_{j=1}^{\mu_n}
\Var \biggl( \sum_{t \in S_j} \theta(t,X_t)
\biggr)\leq\frac{C_2(m_na_n^2+a_n)}{n},
\]
and
\[
\sup_{\theta\in\Theta_{1,n}} \sum_{j=1}^{\mu_n}
\Var \biggl( \sum_{t \in T_j} \theta(t,X_t)
\biggr)\leq\frac{C_2(m_na_n^2+a_n)}{n}.
\]
A simple calculation [observe that $n \rho_n(a_n,\delta) \geq
n^{1/(1+\delta)} \log n \sim m_n$] shows that this implies (i) and
(ii) from the independent blocks argument and that for $E_n$ defined
in \eqref{def:en} we have $E_n = \mathrm{o}(1)$. Thus $\sup_{\theta\in\Theta
_{1,n}}  | \sum_{t=1}^n \theta(t,Y_t)  | = \mathrm{o}_{\rm P}(\eta_n)$.

Next, apply the independent blocks argument with $\Theta_{n,2}$.
Observe that
\[
\sup_{\theta\in\Theta_{1,n}}\sup_t \bigl|\theta(t,X_t)\bigr|
\leq (C+2)n^{-1} \quad\quad\mbox{a.s.}
\]
This yields (i) from the independent blocks argument. Furthermore, we have
\[
\sup_{\theta\in\Theta_{2,n}} \sum_{j=1}^{\mu_n}
\Var \biggl( \sum_{t \in S_j} \theta(t,X_t)
\biggr)\leq C'n^{-2},\quad\quad \sup_{\theta
\in\Theta_{2,n}} \sum
_{j=1}^{\mu_n} \Var \biggl( \sum
_{t \in T_j} \theta(t,X_t) \biggr)\leq
C'n^{-2}
\]
for a constant $C'$ and the same bound holds with $S_j$ replaced by
$T_j$. Thus (ii) from the independent blocks argument is established.
Based on this and the fact that
\[
\sup_{\theta\in\Theta_{2,n}}\sup_t \bigl|\E\bigl[
\theta(t,X_t)\bigr]\bigr| = \mathrm{O}\bigl(n^{-2}\bigr),
\]
some simple calculations show that for $E_n$ defined in \eqref{def:en}
we have $E_n = \mathrm{o}(1)$. This completes the independent blocks argument
for $\Theta_{2,n}$. Combining the results obtained so far establishes
the first part of this lemma. The second part follows from similar arguments.
% Applying the independent blocks argument twice, first with $
%:= D \rho_n(a_n, \delta)$, we first note that, almost surely,
%2 m_n,  \text{and}  \sup_{\theta\in\Theta_{2,n}} \sum_{t
%
%Due to the fact that $n \rho_n(a_n,\delta) \geq n^{1/(1+\delta)} \log
%n \sim m_n$, where $a_n \asymp b_n$ denotes that $a_n$ and $b_n$ are
%of the same order, (i) in the independent blocks argument holds in
%both cases, when $D$ is chosen large enough.
%Choosing $D>2$, we obtain Part (i) of the Lemma, by verifying that
%$E_n = o(1)$ in both cases of applying the independent blocks argument.
%Part (ii) of the Lemma follows along the same lines.
\end{pf}

%
%Let Assumptions (A1) and (A2) hold. Then, for any bounded $\mathcal{Y}
%with $r_n(\eta_n,m_n)$ defined in \eqref{rn}.
%}
%
%$%\[
%|I\{Y_t \leq s\} - I\{X_{t,n} \leq s\}| \leq I\{|D_{t,n}| \geq\eta_n
%%\]
%Thus,
%The first term on the right-hand side of (\ref{equationtheabove}) is
%a.s. non negative; its expected value is bounded by $\kappa_n$, and
%thus the term is of order $O_{\mathrm{P}}(\kappa_n)$. As for the
%second term, we have (recall that, by assumption, $F_{n,X}$ has
%uniformly bounded derivative)
%&&\sup_{s\in\mathcal{Y}} \big|\hat F_{n,X}(s+2\eta_n) - \hat
%F_{n,X}(s-2\eta_n)\big|
%&&  \leq\sup_{s\in\mathcal{Y}} \big|\hat F_{n,X}(s+2\eta_n) -
%&& = O_{\mathrm{P}}(r_n(2\eta_n,m_n)) + O(\eta_n),
%where the last identity follows from Lemma~\ref{lem:df}. This
%completes the proof. \hfill$\Box$

%sA #&#
\section{Technical details for the proof of Theorem~\texorpdfstring{\protect\ref{41414125544254}}{4.1}}\label{appB}

The proof of Theorem~\ref{41414125544254} in Section~\ref{sec4} is
relying on Equations \eqref{lem:ConsistencyA} and \eqref
{lem:ConsistencyB}, which we establish in Sections~\ref{secProof41det_1} and \ref{secProof41det_2}, respectively.
%sA.1 #&#
\subsection{Proof of %the representation
\texorpdfstring{\protect\eqref{lem:ConsistencyA}}{(4.4)}} \label{secProof41det_1}

Putting
\begin{eqnarray*}
4n^{-1} \tilde{\Delta}_n & := & (\beyh{n,
\tau_1, \omega_{j,n}} - \bey{n, \tau_1,
\omega_{j,n}})^{\prime} \pmatrix{ 1 & \mathrm{i}
\cr
-\mathrm{i} &
1 } \bey{n, \tau_2, \omega_{j,n}}
\\
&&{}+ (\bey{n, \tau_1, \omega_{j,n}})^{\prime}
\pmatrix{ 1 & \mathrm{i}
\cr
-\mathrm{i} & 1 } (\beyh{n, \tau_2,
\omega_{j,n}} - \bey{n, \tau_2, \omega_{j,n}})
\\
&&{}+ (\beyh{n, \tau_1, \omega_{j,n}} - \bey{n,
\tau_1, \omega _{j,n}})^{\prime} \pmatrix{1 &
\mathrm{i}
\cr
-\mathrm{i} & 1 } (\beyh{n, \tau_2,
\omega_{j,n}} - \bey{n, \tau_2, \omega_{j,n}}),
\end{eqnarray*}
we obtain, from the definition of the Laplace periodogram,
\begin{eqnarray*}
\Lnr{n,\tau_1,\tau_2} (\omega_{j,n}) & :=&
\frac{n}{4} (\beyh{n, \tau_1, \omega_{j,n}})^\tr
\pmatrix{ 1 & \mathrm{i}
\cr
-\mathrm{i} & 1 }\beyh{n, \tau_2,
\omega_{j,n}}
\\
& =& \frac{n}{4} (\bey{n, \tau_1, \omega_{j,n}})^\tr
\pmatrix{ 1 & \mathrm{i}
\cr
-\mathrm{i} & 1 } \bey{n, \tau_2,
\omega_{j,n}} + \tilde{\Delta}_n
\\
& = &\frac{1}{f(q_{\tau_1}) f(q_{\tau_2})} \bigl( n^{-1} d_n(
\tau_1, \omega_{j,n}) d_n(\tau_2, -
\omega_{j,n}) \bigr) + \tilde{\Delta}_n.
\end{eqnarray*}
By \eqref{thm:uniflinearworank}, for $\tau\in\{\tau_1, \tau_2\}$,
\[
n^{1/2} \sup_{\omega_{j,n} \in\mathcal{F}_n} \|\beyh{n, \tau,
\omega_{j,n}} - \bey{n, \tau, \omega_{j,n}}\| =
\mathrm{O}_{\mathrm{P}} \bigl( n^{(1/8)(1-\delta)/(1+\delta)} (\log n)^{7/4}\bigr),
\]
while Lemma~\ref{lem:boundelta} implies that
\[
n^{1/2} \sup_{\omega_{j,n} \in\mathcal{F}_n} \|\bey{n, \tau,
\omega_{j,n}} \| = \mathrm{O}_{\mathrm{P}}(\sqrt{\log n}),
\]
so that
$\| \tilde{\Delta}_n \| = \mathrm{O}_{\mathrm{P}}( n
\|\beyh{n, \tau, \omega_{j,n}} - \bey{n, \tau, \omega_{j,n}}\|
\cdot\|\bey{n, \tau, \omega_{j,n}} \|) = \mathrm{O}_{\mathrm{P}}(R_n)$. %\]

%sA.2 #&#
\subsection{Proof of %the representation
\texorpdfstring{\protect\eqref{lem:ConsistencyB}}{(4.5)}} \label{secProof41det_2}

%We prove $L_2$-convergence.
Note that $\Lnr{n,\tau_1,\tau_2} (\omega_{j,n})$ is the
cross-periodogram of the bivariate time series
%
%eA.1 #&#
\begin{equation}
\label{bivariatets} \bigl(\tau_1 - I\{Y_t \leq
q_{\tau
_1}\}, \tau_2 - I\{Y_t \leq
q_{\tau_2}\} \bigr).
\end{equation}
Let $\omega_{j,n},  \omega_{k,n} \in(0,\uppi)$ be two sequences of
Fourier frequencies.
Corollary~7.2.2 in Brillinger \cite{Brillinger1975} implies that
%
%eA.2 #&#
\begin{equation}
\label{eqn:VarLtilde} \Var\bigl(\Lnr{n,\tau_1,\tau_2}(
\omega_{j,n})\bigr) = \ef{1,1}(\omega _{j,n}) \ef{2,2} (
\omega_{j,n}) + \frac{2\uppi}{n} \ef {1,2,1,2}(\omega_{j,n},
-\omega_{j,n}, -\omega_{k,n}) + \mathrm{O}(1/n)\quad
\end{equation}
and, for $\omega_{j,n} \neq\omega_{kn}$,
%
%eA.3 #&#
\begin{equation}
\label{eqn:CovLtilde} \Cov \bigl( \Lnr{n,\tau_1,\tau_2}(
\omega_{j,n}), \Lnr{n,\tau_1,\tau_2}(
\omega_{k,n}) \bigr) = \frac{2\uppi}{n} \ef{1,2,1,2}(
\omega_{j,n}, -\omega_{j,n}, -\omega _{k,n}) + \mathrm{O}
\bigl(1/n^2\bigr),
\end{equation}
where $\ef{1,1}$, $\ef{2,2}$ and $\ef{1,2,1,2}$ are the spectra and
cumulant spectra of the bivariate time series (\ref{bivariatets}),
which exist by Assumption~\textup{\ref{(A4)}}. Note that the orders $\mathrm{O}(1/n)$ and
$\mathrm{O}(1/n^2)$ of the remainders in (\ref{eqn:VarLtilde}) and (\ref
{eqn:CovLtilde}) hold uniformly with respect to $j$ and $k$. The
aforementioned spectra coincide with the Laplace spectra $\ef{\tau_1,
\tau_1}$, and $\ef{\tau_2, \tau_2}$ and the cumulant spectra are
also bounded (see Brillinger \cite{Brillinger1975}, page~26). Therefore,
\[
\Cov \bigl( \Lnr{n,\tau_1,\tau_2}(\omega_{j,n}),
\Lnr{n,\tau_1,\tau_2}(\omega_{k,n}) \bigr) =
\cases{ \ef{\tau_1, \tau_1}(\omega_{j,n}) \ef{
\tau_2, \tau_2} (\omega _{j,n}) +
\bar{R}_n, \vspace*{2pt}&\quad $\omega_{j,n} = \omega_{k,n},$
\cr
\bar{R}_n ,&\quad $\omega_{j,n} \neq\omega_{k,n}$, }
\]
where $\bar{R}_n = \mathrm{O}(1/n)$ does not depend on $j$ and $k$.
The assertion follows by the fact that the variance and the bias of the
random variable $K_n$ in \eqref{lem:ConsistencyB}
both are of the order $\mathrm{O}(1/n)$. For the variance, note that
\begin{eqnarray*}
% &  \Var\Big(\sum_{|k| \leq N_n} W_n(k) \Big(\frac{\tilde{L}_n^{
\Var(K_n) & =& \frac{1}{f^2(q_{\tau_1}) f^2(q_{\tau_2})}
\\
&&{}\times
\biggl[\sum_{|k| \leq N_n} W_n^2(k)
\Var \bigl( \Lnr{n,\tau_1,\tau_2}(\omega_{j+k,n})
\bigr)
\\
&& \hphantom{{}\times\biggl[}{}+ \sum_{|k_1| \leq N_n} W_n(k_1)
\mathop{\sum_{ |k_2| \leq
N_n }}_{ k_2 \neq k_1}
W_n(k_2) \Cov \bigl(\Lnr{n,\tau_1,
\tau_2}(\omega_{j+k_1,n}),\overline{\Lnr {n,\tau_1,
\tau_2} (\omega_{j+k_2,n})} \bigr) \biggr]
\\
& =& \mathrm{O}(1/n) ,
\end{eqnarray*}
due to the second part of Assumption \ref{(A3)} and \eqref{eqn:CovLtilde}.
As for the bias, $\mathrm{E}[K_n]= \mathrm{O} (1/n)$
due to the fact that
$\mathrm{E} \Lnr{n,\tau_1,\tau_2}(\omega_{j+k,n}) = \sum_{k=-\infty}^{\infty} \gamma_k(q_{\tau_1}, q_{\tau_2}) \mathrm
{e}^{-\mathrm{i} \omega_{j+k,n} k} + \mathrm{O}(1/n)$
uniformly with respect to the frequencies (see Theorem~4.3.2 in Brillinger \cite
{Brillinger1975}).%\qed

\end{appendix}

% zodis ''Acknowledgments" paliekamas pagal autoriu
\section*{Acknowledgements}
This work has been supported by the Sonderforschungsbereich
``Statistical modelling of nonlinear dynamic processes'' (SFB 823) of
the Deutsche Forschungsgemeinschaft.

Marc Hallin acknowledges the support of a Discovery Grant of the
Australian Research Council, the IAP Research Network Grant P7/06 of
the Belgian federal Government (Belgian Science
Policy), and a Humboldt-Forschungspreis of the Alexander-von-Humboldt-Stiftung.

Tobias Kley is supported by a Ph.D. Grant of the Ruhr-Universit\"at
Bochum and by the Ruhr-Universit\"at Research School funded by
Germany's Excellence Initiative (DFG GSC 98/1).

The constructive remarks and comments by an anonymous referee and a
Associate Editor are gratefully acknowledged.
%suskaldyti doi

% imsref loaded by arune.pranskunaite, 2014-02-28 16:22:37

\printhistory

\end{document}